\documentclass[12pt]{article}

\usepackage[left=1.05in,right=1.05in,top=1.1in,bottom=1.1in]{geometry}
\usepackage[pdftex]{graphicx}
\usepackage{amsfonts,amssymb,amsmath,amsthm}
\usepackage{epstopdf}
\usepackage{longtable}
\usepackage[normalem]{ulem}
\usepackage{graphicx,multirow,float}
\usepackage{hyperref,verbatim,float,stackrel}
\usepackage{bbm}
\usepackage[usenames, dvipsnames]{color}
\usepackage[export]{adjustbox}
\usepackage[font={small}]{caption}
\usepackage{enumitem}

\numberwithin{theorem}{section}

\def\R#1{(\ref{#1})}

\providecommand{\keywords}[1]{\textbf{\small{Keywords:}} {\small{#1}}}

\def\D{\,\mathrm{d}}
\def\I{\mathrm{i}}
\def\E{\mathrm{e}}

\DeclareMathOperator*{\argmin}{argmin}

\DefineNamedColor{named}{Purple}{cmyk}{0.45,0.86,0,0}

\title{Reconstruction of optical vector-fields with applications in endoscopic imaging\footnote{This work has been submitted to the IEEE for possible publication. Copyright may be transferred without notice, after which this version may no longer be accessible.}} 

\author{Milana Gataric\thanks{Center for Mathematical Imaging in Healthcare, Statistical Laboratory, University of Cambridge; m.gataric{@}statslab.cam.ac.uk} 
\and
George S.~D.~Gordon\thanks{Department of Engineering, University of Cambridge; gsdg2@cam.ac.uk} 
\and 
Francesco Renna\thanks{Instituto de Telecomunica\c{c}\~{o}es, Faculdade de Ci\^{e}ncias da Universidade do Porto; frarenna@dcc.fc.up.pt} 
\and 
Alberto Gil C.~P.~Ramos\thanks{Nokia Bell Labs, Cambridge; albertogilramos{@}cantab.net} 
\and 
Maria P.~Alcolea\thanks{Wellcome Trust-MRC Cambridge Stem Cell Institute, University of Cambridge; mpa28{@}cam.ac.uk} 
\and 
Sarah E.~Bohndiek\thanks{Department of Physics and Cancer Research UK Cambridge Institute, University of Cambridge; seb53{@}cam.ac.uk}
}

\begin{document}

\maketitle

\begin{abstract}
We introduce a framework for the reconstruction of the amplitude, phase and polarisation of an optical vector-field using calibration measurements acquired by an imaging device with an unknown linear transformation.
By incorporating effective regularisation terms, this new approach is able to recover an optical vector-field with respect to an arbitrary representation system, which may be different from the one used in calibration.
In particular, it enables the recovery of an optical vector-field with respect to a Fourier basis, which is shown to yield indicative features of increased scattering associated with tissue abnormalities.
We demonstrate the effectiveness of our approach using synthetic holographic images as well as biological tissue samples in an experimental setting where measurements of an optical vector-field are acquired by a fibre endoscope, and observe that indeed the recovered Fourier coefficients are useful in distinguishing healthy tissues from lesions in early stages of oesophageal cancer.
\end{abstract}

\keywords{inverse problems, image reconstruction, blind deconvolution, calibration, Fourier features, optical phase and polarisation, optical fibres, endoscope, oesophageal cancer}

\section{Introduction}

Recently, there has been a significant interest in developing new types of optical fibre endoscopes for medical imaging applications \cite{Cizmar2012,Loterie2015a,Thompson2011,Carpenter2014a,Ba2016}.  Typically, these new endoscopes aim to be thinner, and therefore less invasive, and/or use different properties of light than conventional white light endoscopes making them more sensitive for detecting diseases such as cancer \cite{Waterhouse2018}. A full optical vector-field reflected from a tissue consists of amplitude, phase and polarisation information. Phase and polarisation have recently shown promise as diagnostic indicators, but are discarded by conventional white light endoscopes which record amplitude information only.  Phase is highly sensitive to surface scattering that arises due to microstructural tissue changes in early cancer, creating distorted reflected wavefronts \cite{Drezek2003,Wang2011a, Wang2011b, Su2015}.  This effect has been utilised in phase contrast and quantitative phase microscopy to predict recurrence of prostate cancer \cite{Sridharan2015}.  Similarly, polarisation information can indicate the formation of dense collagen networks \cite{Arifler2007}, and the concentration of other polarisation-sensitive compounds, such as glucose, linked with early cancer  \cite{Kunnen2015, Alali2015}.  This has found use in the diagnosis of colon \cite{Pierangelo2011, Ahmad2015} and gastric cancers \cite{Wang2014}. Currently, there are no commercial phase and polarisation endoscopes but a number of prototype devices have been demonstrated \cite{Qi2016,Cizmar2012, Ploschner2015, Thompson2011, Warren2016}.

To achieve phase and polarisation imaging in fibre endoscopes, it is necessary to characterise the underlying transformation of the optical fibre. In realistic clinical settings, this transformation changes frequently due to bending and temperature fluctuations and it is therefore important that this characterisation is efficient and accurate. 
For the characterisation, typically a set of known fields that form some kind of a basis are input into one end of the fibre and the resulting outputs are recorded at the other end, a procedure termed \emph{calibration}. The task then becomes to recover a representation of the optical field reflected from a tissue given the calibration measurements and the samples of the output field measured by an imaging sensor outside of the fibre.

In this paper, we investigate the following questions:
(i) is there a particularly useful representation of the full optical field reflected from a tissue that can be used for detecting optical aberrations associated with early cancer, and
(ii) how can such a representation be recovered by an efficient and reliable algorithm from raw endoscopic measurements, namely from the calibration measurements and the samples of the output optical field?

To address these questions, we show 
that a Fourier representation recovered directly from the raw measurements has the statistical power to distinguish healthy tissues from tumours, and we provide a general reconstruction framework that can perform such recovery efficiently and stably.

More concretely,
after reviewing related previous works in Section \ref{ss:previous_work}, in Section~\ref{s:framework} we introduce a general reconstruction framework for the recovery of a full two-dimens\-ion\-al complex vector-field, where different regularisation terms are permitted and the bases used for image representation and device calibration are allowed to be different and/or non-orthogonal.
In Section~\ref{s:fourierfeatures} we demonstrate that
it is possible to extract informative features for detecting cancer from
images of simulated tissue samples
by projecting them onto a Fourier basis and observing the decay of their respective Fourier coefficients.
Finally, in Section~\ref{s:numerics}, we apply our new approach to experimental data acquired using a custom-built fibre endoscope \cite{Gordon2018} and recover synthetic holographic images as well as images of mouse oesophageal tissue containing small tumours (lesions). In particular, by recovering images of a biological tissue with respect to a Fourier basis using $\ell_1$-regularisation, we observe that the corresponding Fourier coefficients are indicative of differences between lesions and healthy tissues and demonstrate their potential for medical diagnostic applications. We conclude with a discussion of our results
and directions for future research
in Section \ref{s:discussion}.

\subsection{Relation to previous work}
\label{ss:previous_work}

Before precisely formulating our new approach in the next section, it is useful to give a brief overview of some existing reconstruction techniques used in the applications of interest and differentiate them from the reconstruction framework developed in this paper. 

In imaging through optical fibres or other scattering media, typical recovery procedures use the same, finite-dimensional basis for calibration and image representation in conjunction with standard inversion techniques. They start by discretizing the mathematical operator of the fibre as a mapping between pixels at different ends of the fibre $\mathbf{A}:\mathbf{x}\in\mathbb{C}^p \mapsto\mathbf{y}\in\mathbb{C}^n$, leading to a transmission matrix $\mathbf{A}\in\mathbb{C}^{n\times p}$ which is then characterised through calibration. The individual calibration inputs are arranged into the columns of matrix $\mathbf{X}_{cal}\in\mathbb{C}^{p\times m}$ and the corresponding calibration outputs into the columns of matrix $\mathbf{Y}_{cal}\in\mathbb{C}^{n\times m}$. Most existing systems use a full orthogonal basis of the discretized input space as the calibration inputs, e.g.~a set of tilted plane waves (a Fourier basis) \cite{Cizmar2009} or a Hadamard basis generated using a phase-only spatial light modulator \cite{Popoff2010}. The orthogonality of such bases ensures that matrix $\mathbf{X}_{cal}$ is unitary. Then, by assuming that $\mathbf{A}$ is also unitary, images can be recovered using \emph{phase conjugation}. In this approach a (generalised) inverse of the transmission matrix is calculated as $\mathbf{X}_{cal}\mathbf{Y}_{cal}^{*}$, where $\cdot^{*}$ denotes the conjugate transpose, and a representation of $\mathbf{x}$ with respect to the calibration inputs is recovered as $\mathbf{X}_{cal}\mathbf{Y}_{cal}^{*}\mathbf{y}$  \cite{Popoff2010,Cizmar2012,Cizmar2011}. Although simple and straightforward to compute, the unitary assumptions in this approach are typically violated in practice \cite{Carpenter2014}.  In the context of imaging through scattering media, the inversion of the transmission matrix was also performed through alternative approaches to phase-conjugation such as least-squares (Moore-Penrose pseudoinverse) or Tikhonov regularisation \cite{Popoff2010a,Popoff2011a}. In particular, $\ell_1$-regularisation has not been previously explored in these applications.

When compared with these conventional techniques, we emphasise that our new framework is able to recover a representation of the unknown optical field with respect to any particular infinite-dimensional basis which is allowed to be different from the one used for calibration, directly from the raw measurements.
If an image representation with respect to a particular basis (such as Fourier) is desired, 
alternatively to our new approach one could in principle use the
conventional techniques to recover an approximation to such a representation as we now describe. One could calibrate the fibre with respect to a Fourier basis and use standard recovery techniques to reconstruct images with respect to the same basis. However, in high resolution imaging, calibration with respect to a Fourier basis may become prohibitively slow in practice and it may be preferable to use different, more efficient systems for calibration, as we do in this paper. Another possible approach would be to first recover the image with respect to the calibration basis and then approximate its Fourier coefficients in a post-processing step. However, as a two-stage procedure, such approach is inherently less efficient and  suffers from greater error than the approach proposed in this paper, which can recover Fourier coefficients directly from the raw measurements.

In the earlier work \cite{Gordon2018}, a fibre endoscope was developed to produce images of phase and polarisation for early cancer detection. In this case, a set of calibration inputs was chosen so as to greatly speed up experimental characterisation measurement time by enabling parallelized calibration that exploits the localised confinement of light of the underlying fibre structure. However, this input basis is non-orthogonal and so phase-conjugation cannot be naively applied. Using the framework presented in this paper, we are able to reconstruct phase and polarisation images with respect to a diagnostically relevant representation system, while simultaneously preserving the benefits of efficient experimental calibration achieved with a system tailored to the fibre structure. Moreover, this new approach significantly decreases the reconstruction time compared to the previously implemented reconstruction technique \cite{Gordon2018}, providing an important advance towards real-time image reconstruction.

Finally, we mention that changing representation systems between image recovery and sampling has previously been applied to inverse problems arising in various image and signal processing applications (see \cite{Adcock2014,Adcock2016} and references therein). Typically, it is assumed that the imaging device of a known linear transformation provides image samples with respect to a specified sampling system, while the aim is to recover a  representation of the image with respect to a different system chosen so that a good approximation of the image is obtained or the number of required samples is decreased. As in this paper, the systems considered are modelled by Riesz bases or frames of infinite-dimensional function-spaces. By contrast, in this paper the imaging device produces pixel samples of a transformed image where the underlying transformation is unknown and is characterised through a calibration procedure.

\section{Reconstruction framework}\label{s:framework}

In this section we introduce our reconstruction framework. We start by presenting an infinite-dimensional imaging model in Section~\ref{ss:model}. We then consider a simplified scalar-valued setting in Section~\ref{ss:scalars}, where we derive a linear system of equations and its regularised solution while providing flexibility in choosing different systems for calibration and image representation. We then extend our framework to the vector-field case in Section~\ref{ss:vectors}.

\subsection{Imaging model and reconstruction problem}\label{ss:model}

In imaging through fibres or other scattering media, an input optical vector-field $\mathbf{F}$ is related to its corresponding output optical vector-field $\tilde{\mathbf{F}}$ through an integral transformation with a spatially-varying kernel $\mathbf{G}$, also called Green's function or point-spread function. Specifically, such transformation can be modelled by
\begin{equation}
\label{eq:model}
\tilde{\mathbf{F}}(\mathbf{y}) = \int_S  \mathbf{G}(\mathbf{y},\mathbf{x}) \mathbf{F}(\mathbf{x}) \D \mathbf{x},
\end{equation}
where $\mathbf{F}:S\rightarrow \mathbb{C}^2$ is a complex-valued vector-field representing the unknown optical field on the input plane $S\subseteq\mathbb{R}^2$, $\tilde{\mathbf{F}}:\mathbb{R}^2\rightarrow \mathbb{C}^2$ is a complex-valued vector-field on the output plane which can be sampled, and  $\mathbf{G}:\mathbb{R}^2\times\mathbb{R}^2\rightarrow \mathbb{C}^{2\times2}$ is some unknown bounded matrix-valued function\footnote{In general, the kernel  $\mathbf{G}$ is also time-dependent as it depends on many factors such as bending of the fibre and temperature. In this paper, we account for significant measurement noise but only consider imaging at a single time point; c.f.~Section~\ref{s:discussion}.} \cite{Rotter2017}. In particular, we consider the input field $\mathbf{F}$ to be an object with infinite resolution, and thus, we model $\mathbf{F}$ as an element of an infinite-dimensional function-space, such as the $\mathcal{L}^2$-space of square-integrable vector-valued functions.

In endoscopic imaging, we are interested in capturing a full optical field (i.e., amplitude, phase and polarisation) reflected from a human tissue inside the body, which is also called a wavefront, and which in this paper, we refer to as an image. In the terminology above, $\mathbf{F}$ denotes an image, which is observed indirectly at the input imaging plane $S$ located at the end of the fibre placed inside the body, termed as the \emph{distal} facet of the fibre. The fibre then transports light from the \emph{distal} facet $S$ inside the body to the \emph{proximal} facet outside the body where the imaging sensor directly observes $\tilde{\mathbf{F}}$ at the output imaging plane. The question then becomes how to recover the unknown, infinite-dimensional optical field $\mathbf{F}$ from the acquired samples of $\tilde{\mathbf{F}}$.


More concretely, given the pointwise measurements of the output vector-field $\tilde{\mathbf{F}}$ collected at the imaging sensor
\begin{equation}\label{eq:samples}
 \tilde{\mathbf{F}}(\mathbf{y}_n), \quad n=1,\ldots,N, 
\end{equation}
where $\mathbf{y}_n\in\mathbb{R}^2$ and $N\in\mathbb{N}$ is the resolution of the imaging sensor, the goal is  to recover the unknown function $\mathbf{F}$ via equation \R{eq:model}. It is important to note that these measurements will also contain noise introduced by the measurement procedure.

This linear inverse problem is especially challenging because both the spatially-varying kernel $\mathbf{G}$ as well as the eigenfunctions associated with the underlying integral transform \R{eq:model} are unknown.  Such eigenfunctions are termed \emph{modes} of the fibre and their analytic form is available only for some limited fibres such as parabolic graded index multimode fibres \cite{Snyder1983}.

To be able to recover $\mathbf{F}$ from finitely many samples of $\tilde{\mathbf{F}}$ in scenarios where neither $\mathbf{G}$ nor the eigenfunctions are known, one strategy that may be considered is to employ a calibration procedure. Concretely, it is possible to design calibration input fields $\mathbf{E}_m$, $m=1,\ldots,M$, and to measure the corresponding output fields $\tilde{\mathbf{E}}_m$, which in line with the notation above are vector-valued functions related through the infinite-dimensional model given in \R{eq:model}.
The advantage of calibration is that we now have access not only to the data given in \R{eq:samples} but also to the calibration data
\begin{equation}\label{eq:calibrationdata}
 \mathbf{E}_m, \quad \tilde{\mathbf{E}}_m(\mathbf{y}_n), \quad m=1,\ldots,M, \ n=1,\ldots, N,
\end{equation}
which forms additional information with which to recover $\mathbf{F}$.

It is noted that while the output fields $\tilde{\mathbf{E}}_m$ are sampled at an output imaging sensor of resolution $N$, the calibration input fields $\mathbf{E}_m$ can be evaluated on a discretised grid whose resolution does not depend on any physical limitation imposed by the fibre or by the sensor collecting the transmitted image; it only depends on the resolution of the sensors used for calibration, which may be much larger than $M$. Therefore, as for the input $\mathbf{F}$, we model the inputs $\mathbf{E}_m$ as elements of an infinite-dimensional function-space. 
Thus, 
the representation of $\mathbf{F}$
as well as the device calibration can be considered with respect to a wide class of infinite-dimensional bases or over-complete systems that may not be orthogonal.

\subsection{Reconstruction of scalar-fields}\label{ss:scalars}


We begin approaching the general problem of recovering the complex vector-field $\mathbf{F}$, described in Section \ref{ss:model}, by first solving a simpler but related problem, which once solved will provide us with the methodology and insights necessary to tackle the problem in its full generality in the next Section \ref{ss:vectors}. Specifically, we assume in this subsection that $F, \tilde{F}, E_m, \tilde{E}_m$ are scalar valued functions that take values in $\mathbb{C}$ rather than $\mathbb{C}^2$, and accordingly $G$ takes values in $\mathbb{C}$ rather than $\mathbb{C}^{2\times2}$. We highlight this difference by writing these quantities with non-bold symbols.

Our approach starts by considering all fields on the input imaging plane $S$ as elements of the same function-space $\mathcal{F}$, such as the $\mathcal{L}^2$-space of square-integrable scalar-valued functions supported on $S$, with inner product defined as $\langle E,H \rangle := \int_{S} E(\mathbf{x}) H^*(\mathbf{x}) \D \mathbf{x}$, for $E,H\in\mathcal{F}$. 
We proceed by recovering $F\in\mathcal{F}$ at resolution $K\in\mathbb{N}$ in terms of some desired representation system $\{H_k\}_{k=1}^K$ in $\mathcal{F}$, using only the available data in \R{eq:samples} and \R{eq:calibrationdata}. Specifically, we aim to estimate the coefficients $\mathbf{f} = \begin{bmatrix} f_1,\ldots, f_K\end{bmatrix}^{\top}\in\mathbb{C}^K$ of the $K$-term approximation of $F$ given as
\begin{equation}\label{eq:Kterm}
F_K(\mathbf{x}) := \sum_{k=1}^K f_k H_k(\mathbf{x}), \quad \mathbf{x}\in S.
\end{equation}

Before turning to the computation of $f_k$ in \eqref{eq:Kterm},
it is insightful to work through special cases of $S$ and $\{H_k\}_{k=1}^K$ that are particularly useful in practice. For instance, if the aim is to recover a Fourier representation of $F$ and $S:=[-1/2,1/2]^2\subseteq\mathbb{R}^2$ is the unit square, then $\{H_k\}_{k=1}^K$ is the $K$-dimensional Fourier basis $\{e^{2\pi\I \mathbf{k}\cdot\mathbf{x}}\}_{\mathbf{k}\in I_{K}}$ where $I_K:=\{\mathbf{k}=(k_1,k_2)\in\mathbb{Z}^2:k_1,k_2=-\sqrt{K}/2,\ldots,\sqrt{K}/2-1\}$, $\mathbf{k}\cdot\mathbf{x}:=k_1x_1+k_2x_2$, $\mathbf{x}:=(x_1,x_2)\in S$, and \eqref{eq:Kterm} specialises to
\begin{equation}\label{eq:fourierceof}
F_K(\mathbf{x}) := \sum_{\mathbf{k}\in I_K} f_{\mathbf{k}} e^{2\pi\I \mathbf{k}\cdot\mathbf{x}}, \quad f_{\mathbf{k}}:=\int_{S} F(\mathbf{x})e^{-2\pi\I \mathbf{k}\cdot\mathbf{x}}\D \mathbf{x}.
\end{equation}
More generally, $\{H_k\}_{k=1}^K$ may contain the first $K$ elements of a Riesz basis in $\mathcal{F}$, such as B-spline wavelets for example, with its corresponding biorthogonal sequence denoted by $\{\breve{H}_k\}_{k=1}^K$, in which case \eqref{eq:Kterm} specialises to  $F_K(\mathbf{x}) = \sum_{k=1}^{K}\langle F,\breve{H}_k \rangle H_k$, i.e., $f_k = \langle F,\breve{H}_k \rangle$. Moreover, since for the reconstruction framework we do not require an explicit form of the coefficients $f_k$, the notion of basis can be further relaxed to over-complete representation systems such as those of over-complete frames \cite{Christensen2002}.

Returning to the key issue of approximating the coefficients $f_k$ in the general $K$-term approximation $F_K$ \eqref{eq:Kterm} from the given measurements \R{eq:samples}--\R{eq:calibrationdata}, we write each $H_k$ in terms of the calibration functions $\left\{E_m\right\}_{m=1}^M$ as 
\begin{equation}\label{eq:Happrox}
H_k(\mathbf{x}) = \sum_{m=1}^M h_{m,k} E_m(\mathbf{x}) + \delta_k(\mathbf{x}),  \quad \mathbf{x}\in S,
\end{equation}
for some coefficients $h_{m,k}\in\mathbb{C}$, whose computation we discuss below, and for some error term $\delta_k$. Since \R{eq:model} is a linear transformation of $F$, by 
substituting $F$ with $F_K+(F-F_K)$ in \R{eq:model} and writing $F_K$ in terms of \eqref{eq:Kterm} and \eqref{eq:Happrox}, we have 
\begin{equation}\label{eq:firstobservation}
\tilde{F}(\cdot) = \sum_{k=1}^K \sum_{m=1}^M f_k h_{m,k} \tilde E_m(\cdot) + \left( \sum_{k=1}^K f_k \int_S G(\cdot,\mathbf{x}) \delta_k(\mathbf{x}) \D \mathbf{x} + \int_S G(\cdot,\mathbf{x}) (F(\mathbf{x})-F_K(\mathbf{x})) \D \mathbf{x}\right).
\end{equation}
By evaluating  equation \eqref{eq:firstobservation} at the measurement points $\{\mathbf{y}_n\}_{n=1}^N$, we obtain the following linear system 
\begin{equation}\label{eq:linsys_CB}
 \mathbf{g} = \mathbf{E} \mathbf{H} \mathbf{f} + \boldsymbol{\varepsilon},
\end{equation}
where $\mathbf{g}\in\mathbb{C}^{N}$ is the vector having its $n$-th entry equal to $\tilde F(\mathbf{y}_{n})$,
$\mathbf{E}\in\mathbb{C}^{N\times M}$ is the matrix having its $(n,m)$-th entry equal to $\tilde E_m (\mathbf{y}_n)$,
$\mathbf{H}\in\mathbb{C}^{M\times K}$ is the matrix 
with its $(m,k)$-th entry equal to $h_{m,k}$
and
$\boldsymbol{\varepsilon}\in\mathbb{C}^N$ is an error term containing  the last two terms in the right-hand-side of \R{eq:firstobservation}. In addition, the error term $\boldsymbol{\varepsilon}\in\mathbb{C}^N$ can be seen also as encapsulating measurement error noise incurred when measuring $\tilde{F}(\mathbf{y}_n)$ and $\tilde{E}_m(\mathbf{y}_n)$ in \eqref{eq:samples} and \eqref{eq:calibrationdata}, respectively.
We then opt to define the solution of \R{eq:linsys_CB} as the minimisation problem
\begin{equation}\label{eq:reg}
\bar{\mathbf{f}}  := \argmin_{\mathbf{f}\in\mathbb{C}^K } \left\{ \| \mathbf{g} - \mathbf{E} \mathbf{H} \mathbf{f} \|_2 + \lambda \mathcal{R}(\mathbf{f}) \right\},
\end{equation}
where $\|\cdot\|_2$ denotes the Euclidean norm on $\mathbb{C}^N$, while the regularisation term $\mathcal{R}$ and its parameter $\lambda\geq0$ are described below. Once the coefficients $\bar{\mathbf{f}} = \begin{bmatrix}\bar f_1,\ldots,\bar f_K\end{bmatrix}^{\top} \in \mathbb{C}^K$ are computed through \eqref{eq:reg}, then in line with \eqref{eq:Kterm} we define the reconstruction of $F$ as the approximation given by
\begin{equation}\label{eq:recK}
\bar F_K(\mathbf{x}) := \sum_{k=1}^{K} \bar f_k H_k(\mathbf{x}), \quad \mathbf{x}\in S.
\end{equation}

To obtain the explicit solution defined in \eqref{eq:recK}, it remains to describe the procedure for computing the coefficients of matrix $\mathbf{H}$ and to define the regularisation term $\mathcal{R}$.

First, observe that if the same system is used for calibration and reconstruction, then $\mathbf{H}=\mathbf{I}$. 
Otherwise, we can estimate $\mathbf{H}$ as follows.
Using \R{eq:Happrox}, we write
$\langle H_k, E_{m'} \rangle = \sum_{m=1}^M h_{m,k} \langle E_m,E_{m'} \rangle + \langle \delta_k, E_{m'} \rangle$, $m'=1,\ldots,M$, and therefore, provided $\langle \delta_k, E_{m'} \rangle\approx0$, we have
$$ \mathbf{H} 
\approx
\begin{bmatrix}
 \langle E_1,E_1 \rangle & \ldots & \langle E_M,E_1 \rangle \\
 \vdots & & \vdots  \\
 \langle E_1,E_M \rangle & \ldots & \langle E_M,E_M \rangle
\end{bmatrix}^{-1}
\begin{bmatrix}
 \langle H_1,E_1 \rangle & \ldots & \langle H_K,E_1 \rangle \\
 \vdots & & \vdots  \\
 \langle H_1,E_M \rangle & \ldots & \langle H_K,E_M \rangle
\end{bmatrix}. 
$$
The matrix featuring this linear system is known as the Gram matrix of $\{E_m\}_{m=1}^M$, which takes the form of the identity matrix when $\{E_m\}_{m=1}^M$ are orthonormal. 
We note that the accuracy of such estimation of matrix $\mathbf{H}$ and its condition number depend on the gap between the function-spaces spanned by $\{H_k\}_{k=1}^K$ and $\{E_m\}_{m=1}^M$ as well as on the conditioning of the Gram matrix. In general, a good estimation requires that $\{E_m\}_{m=1}^M$ forms a good approximation for $\{H_k\}_{k=1}^K$.

We now discuss the choice of the regularisation term $\mathcal{R}$ in \eqref{eq:reg}. If $\mathcal{R}$ is absent from \R{eq:reg}, i.e.~if $\lambda=0$, then the solution of this minimisation problem is equivalent to the least-squares solution $\bar{\mathbf{f}}:=((\mathbf{E} \mathbf{H})^* \mathbf{E} \mathbf{H} )^{-1} (\mathbf{E} \mathbf{H})^* \mathbf{g}$.  If the regularisation term is given by $\mathcal{R}(\mathbf{f}):=\|\mathbf{f}\|_2$, then \eqref{eq:reg} is known as Tikhonov regularisation and its solution is given by $\bar{\mathbf{f}}:=((\mathbf{E} \mathbf{H})^* \mathbf{E} \mathbf{H} + \lambda \mathbf{I} )^{-1} (\mathbf{E} \mathbf{H})^* \mathbf{g}$. However, if $\mathbf{E} \mathbf{H}$ is badly conditioned, $\boldsymbol{\varepsilon}>0$ and, additionally, it is known a priori that only a few elements of $\{H_k\}_{k=1}^K$ are sufficient to represent $F$ well, then $\mathcal{R}(\mathbf{f}):=\|\mathbf{f}\|_0$ is an appropriate choice of the regularisation term. This is termed as the $\ell_0$-regularisation, where the $\ell_0$-norm of $\mathbf{f}$, i.e.~$\|\mathbf{f}\|_0$, is defined as the number of non-zero coordinates in $\mathbf{f}$. The $\ell_0$-regularisation bypasses the ill-conditioning by imposing sparsity in the solution $\bar F_K$ with respect to $\{H_k\}_{k=1}^K$.  In practice, solving the minimisation problem with such a non-convex $\ell_0$-term is  computationally difficult, so typically an $\ell_1$-relaxation is considered instead. The corresponding relaxed minimisation problem can then be solved by fast iterative algorithms \cite{Berg2008,Becker2011}. In addition to the choice of the form of the regularisation term $\mathcal{R}$, the strength of the regularisation is controlled by the parameter $\lambda$, which can be chosen by cross-validation techniques \cite{Doostan2011}.


We conclude this subsection by a  discussion on the accuracy and robustness of the solution defined in \eqref{eq:recK}.

The reconstruction error can be quantified by the magnitude of $F-\bar{F}_K=(F-F_K)+(F_K-\bar{F}_K)$, which depends on several factors. The magnitude of the first term $F-F_K$ depends of how well $F$ can be represented by its $K$-term approximation with respect to $\{H_k\}_{k=1}^K$, and thus it is expected to decrease with increasing $K$. On the other hand, the magnitude of the second term $F_K-\bar{F}_K$ depends on the conditioning of the matrix $\mathbf{E}\mathbf{H}$ and the error term $\boldsymbol{\varepsilon}$ in \R{eq:linsys_CB}, and thus, it is expected to increase with increasing $K$ when $M$ and $N$ are fixed. In other words, if the resolution $K$ at which we reconstruct is increased, we also need to increase $M$ and $N$. However, it may be possible to attain higher resolutions if some form of regularisation is used when solving  \R{eq:linsys_CB}.

As previously noted, the error term $\boldsymbol{\varepsilon}$ contains the measurement error as well as the last two terms of the right hand side in \eqref{eq:firstobservation}, which can be disregarded provided  $F - F_K$ and $\delta_k$ are small or they lie in the span of those eigenfunctions corresponding to a small singular value. Thus, for small $\boldsymbol{\varepsilon}$ it is required that $\{H_k\}_{k=1}^K$ and $\{E_m\}_{m=1}^M$ form a good approximation for $F$  or for the eigenfunctions with large singular values. 
However, if the singular values of the underlying integral operator accumulate at zero, the conditioning of the matrix $\mathbf{E}$ may  become worse if the span of $\{E_m\}_{m=1}^M$ includes too many eigenfunctions including those corresponding to a small singular value. 
Loosely speaking, the calibration functions $\{E_m\}_{m=1}^M$ should form a good representation for the span of the eigenfunctions, excluding those eigenfunctions corresponding to small singular values if they exist. However, as we do not have access to the true eigenfunctions,  we do not have control over the ill-conditioning that we are introducing by using a particular choice of $\{E_m\}_{m=1}^M$. Thus, the use of regularisation in solving \R{eq:linsys_CB} becomes crucial in order to obtain a robust solution.

\subsection{Reconstruction of vector-fields}\label{ss:vectors}

We now extend the scalar-field reconstruction framework developed in Subection \ref{ss:scalars} to the more general vector-field problem presented in Subection \ref{ss:model}. To begin with, let
$$
\mathbf{F}:=
\begin{bmatrix}
 F^{h} \\
 F^{v}
\end{bmatrix} \overset{\R{eq:model}}{\mapsto}
\tilde{\mathbf{F}}:=
\begin{bmatrix}
\tilde{F}^{h} \\
 \tilde{F}^{v}
\end{bmatrix}
$$
be the complex-vector-valued functions related as in equation \R{eq:model},
where the superscripts $h$ and $v$ correspond to the horizontal and vertical polarisations of the optical vector-field, respectively. The goal is to recover both polarisations $F^h$ and $F^v$, which are scalar-valued functions, by using the measurements \R{eq:samples} and \R{eq:calibrationdata}. Since each of $\tilde{F}^{h}$ and $\tilde{F}^{v}$ depends on both $F^v$ and $F^h$, we cannot consider the reconstructions of  $F^v$ and $F^v$ as two independent problems. However, as it will be demonstrated below, we can still use the reconstruction framework introduced earlier in Section \ref{ss:scalars}, as long as the calibration inputs can form a representation system for vector-valued functions such as $\mathbf{F}$. To ensure that this is the case,
rather than straightforwardly sampling the vector-valued calibration inputs from \R{eq:calibrationdata}, i.e.
\begin{equation*}
\mathbf{E}_m:=
\begin{bmatrix}
 E^{h}_m \\
 E^{v}_m
\end{bmatrix}
, \quad
\tilde{\mathbf{E}}_m(\mathbf{y}_n):=
\begin{bmatrix}
\tilde{E}^{h}_m(\mathbf{y}_n) \\
 \tilde{E}^{v}_m(\mathbf{y}_n)
\end{bmatrix}
\end{equation*}
we instead sample the following two related forms
\begin{align}\label{eq:cal_funs_vf}
\mathbf{A}_m \!:=\! 
\begin{bmatrix}
 E_m^{h} \\
 E_m^{v}
\end{bmatrix}
\overset{\R{eq:model}}{\mapsto} 
\tilde{\mathbf{A}}_m \!:=\!
\begin{bmatrix}
\tilde{A}_m^{h} \\
\tilde{A}_m^{v}
\end{bmatrix},\qquad
\mathbf{B}_m \!:=\!
\begin{bmatrix}
 E_m^{h} \\
 b E_m^{v}
\end{bmatrix}
\overset{\R{eq:model}}{\mapsto}
\tilde{\mathbf{B}}_m \!:=\!
\begin{bmatrix}
\tilde{B}_m^{h} \\
\tilde{B}_m^{v}
\end{bmatrix}, 
\end{align}
where $m=1,\ldots,M$, $b:=\E^{\beta\I}$ for a fixed $\beta\in(0,2\pi)$ and $E_{m}^{h},E_{m}^{v},\tilde{A}_m^{h},\tilde{A}_m^{v},\tilde{B}_m^{h},\tilde{B}_m^{v}$ are some scalar-valued functions.
To make clear the motivation to sample according to \eqref{eq:cal_funs_vf}, observe that if $\{E_{m}^{h}\}_{m=1}^{M}$ and $\{ E_{m}^{v}\}_{m=1}^{M}$ are representation systems for $F^h$ and $F^v$ respectively, then we have that
$\{\mathbf{A}_m - \mathbf{B}_m\}_{m=1}^{M}$  and $\{\mathbf{A}_m - b^* \mathbf{B}_m\}_{m=1}^{M}$ are representation systems for $\begin{bmatrix} 0,   F^{v} \end{bmatrix}^{\top}$ and $\begin{bmatrix} F^{h},  0 \end{bmatrix}^{\top}$ respectively. In other words,  
\begin{equation}\label{eq:rep_sys_vf}
\bigl\{\mathbf{A}_m-\mathbf{B}_m,\ \mathbf{A}_m-b^*\mathbf{B}_m \ :\  m=1,\ldots,M\bigr\}
\end{equation}
can be used to represent the complex vector-valued function $\mathbf{F}$.

Mimicking  the reasoning of the previous subsection,
we proceed by recovering approximations of $F^h$ and $F^v$ with respect to some desired representation systems $\{H_k^h\}_{k=1}^{K}$ and $\{H_k^v\}_{k=1}^{K}$. Namely we aim to recover
$$
F^h_{K}(\mathbf{x}):=\sum_{k=1}^{K} f_k^h H_k^h(\mathbf{x}), \qquad F^v_{K}(\mathbf{x}):=\sum_{k=1}^{K} f_k^v H_k^v(\mathbf{x}), \qquad \mathbf{x}\in S,
$$
where we first write these representation systems in terms of the calibration functions as
$H^h_k(\mathbf{x})=\sum_{m=1}^{M} h_{m,k}^h E_m^h(\mathbf{x})+\delta^h_k(\mathbf{x})$ and $H^v_k(\mathbf{x})=\sum_{m=1}^{M} h_{m,k}^v E_m^v(\mathbf{x})+\delta^v_k(\mathbf{x})$,
for some coefficients $f_k^h,f_k^v,h_m^h,h_m^v\in\mathbb{C}$ and some error terms $\delta^h_k,\delta^v_k$. It thus follows that
\begin{align}\label{eq:f_ex_pol}
\mathbf{F}(\mathbf{x}) \!=\! 
\sum_{k=1}^{K}\!  \sum_{m=1}^{M}\!  f_k^h h_{m,k}^h \begin{bmatrix}  E_m^h(\mathbf{x}) \\ 0 \end{bmatrix} \!+\! 
\sum_{k=1}^{K}\!  \sum_{m=1}^{M}\!  f_k^v h_{m,k}^v \begin{bmatrix} 0 \\ E_m^v(\mathbf{x}) \end{bmatrix}
\!+\!  \sum_{k=1}^{K} \!  \begin{bmatrix} f_k^h\delta_k^h(\mathbf{x}) \\ f_k^v\delta_k^v(\mathbf{x}) \end{bmatrix}
\!+\!  \begin{bmatrix}
F^h(\mathbf{x})\!-\!F^h_{K}(\mathbf{x})\\
F^v(\mathbf{x})\!-\!F^v_{K}(\mathbf{x})
\end{bmatrix}\!.
\end{align}
Since
\begin{align*}
 \begin{bmatrix} 
  E_m^h(\mathbf{x}) \\
  0
 \end{bmatrix} = \frac{1}{1-b^*} \left( \mathbf{A}_m(\mathbf{x}) - b^* \mathbf{B}_m(\mathbf{x}) \right), \qquad
 \begin{bmatrix} 
 0 \\
 E_m^v(\mathbf{x})
 \end{bmatrix} = \frac{1}{1-b} \left( \mathbf{A}_m(\mathbf{x}) -  \mathbf{B}_m (\mathbf{x}) \right),
\end{align*} 
by applying \R{eq:model} to \R{eq:f_ex_pol}, we obtain
\begin{align*}
\tilde{\mathbf{F}}(\cdot) \approx \frac{1}{1-b^*} \sum_{k=1}^{K} \sum_{m=1}^{M} f_k^h h_{m,k}^h  \left( \tilde{\mathbf{A}}_m(\cdot) - b^* \tilde{\mathbf{B}}_m (\cdot) \right) + \frac{1}{1-b} \sum_{k=1}^{K} \sum_{m=1}^{M} f_k^v h_{m,k}^v \left( \tilde{\mathbf{A}}_m(\cdot) - \tilde{\mathbf{B}}_m (\cdot) \right),
\end{align*}
provided the two last terms in \R{eq:f_ex_pol} are small or they become small after applying \R{eq:model}.
By using the pointwise measurements from \R{eq:samples} and \R{eq:cal_funs_vf}, this leads to the following linear system
\begin{equation}\label{eq:pol_system}
 \mathbf{g} =  \mathbf{E} \mathbf{H} \mathbf{f} + \boldsymbol{\varepsilon}, \qquad \mathbf{E}:=
\begin{bmatrix}
  a^*\left(\mathbf{A} - b^*\mathbf{B}\right)  & a\left(\mathbf{A} - \mathbf{B}\right)
 \end{bmatrix}, \ \mathbf{H}:=
 \begin{bmatrix}
   \mathbf{H}^h & 0 \\
   0 &  \mathbf{H}^v
 \end{bmatrix}
\end{equation}
where $a:=1/(1-b)\in\mathbb{C}$, $\mathbf{H}^h\in\mathbb{C}^{M\times K}$ is such that its $(m,k)$-th entry is $h^h_{m,k}$, $\mathbf{H}^v\in\mathbb{C}^{M\times K}$ is such that its $(m,k)$-th entry is $h^v_{m,k}$, $\boldsymbol{\varepsilon}\in\mathbb{C}^{2N}$ is the error term, and $\mathbf{g}\in\mathbb{C}^{2N}$, $\mathbf{A},\mathbf{B}\in\mathbb{C}^{2N\times M}$, $\mathbf{f}\in\mathbb{C}^{K}$ are defined as
\begin{equation*}
\mathbf{g} \!:=\!
\begin{bmatrix}
\tilde F^h(\mathbf{y}_1)\\
 \tilde F^v(\mathbf{y}_1)\\
 \vdots,\\
 \tilde F^h(\mathbf{y}_{N})\\
 \tilde F^v(\mathbf{y}_{N})
\end{bmatrix}\!\!,\
\mathbf{A} \!:=\! \begin{bmatrix}
\tilde{A}_1^h(\mathbf{y}_1) & \ldots  & \tilde{A}_{M}^h(\mathbf{y}_1) \\
\tilde{A}_1^v(\mathbf{y}_1) & \ldots  & \tilde{A}_{M}^v(\mathbf{y}_1) \\
    &   \ddots & \\
\tilde{A}_1^h(\mathbf{y}_N) & \ldots  & \tilde{A}_{M}^h(\mathbf{y}_N) \\
\tilde{A}_1^v(\mathbf{y}_N) & \ldots  & \tilde{A}_{M}^v(\mathbf{y}_N) \\
\end{bmatrix}\!\!,\
\mathbf{B} \!:=\! \begin{bmatrix}
\tilde{B}_1^{h}(\mathbf{y}_1) & \ldots  & \tilde{B}_{M}^{h}(\mathbf{y}_1) \\
\tilde{B}_1^{v}(\mathbf{y}_1) & \ldots  & \tilde{B}_{M}^{v}(\mathbf{y}_1) \\
    &   \ddots & \\
\tilde{B}_1^{h}(\mathbf{y}_N) & \ldots  & \tilde{B}_{M}^{h}(\mathbf{y}_N) \\
\tilde{B}_1^{v}(\mathbf{y}_N) & \ldots  & \tilde{B}_{M}^{v}(\mathbf{y}_N) \\
\end{bmatrix}\!\!,\
\mathbf{f} \!:=\!
\begin{bmatrix}
 f^h_1 \\
 \vdots \\
 f^h_{K}\\
 f^v_1 \\
 \vdots \\
 f^v_{K}
\end{bmatrix}\!\!.
\end{equation*}
We propose to solve the linear system \R{eq:pol_system} in a similar manner to that used in \R{eq:reg}. Finally, once \R{eq:pol_system} is solved for the coefficients
$\bar{\mathbf{f}}=
\begin{bmatrix}
 \bar f^h_1,
 \ldots,
 \bar f^h_{K},
 \bar f^v_1,
 \ldots,
 \bar f^v_{K}
\end{bmatrix}^{\top} \in \mathbb{C}^{2K}$, we can define the reconstructions of $F^h$ and $F^v$ as
\begin{equation*}
\bar F_K^h(\cdot):=\sum_{m=1}^{K} \bar f_k^h H_k^h(\cdot), \qquad \bar F_K^v(\cdot):=\sum_{m=1}^{K} \bar f_k^v H_k^v(\cdot).
\end{equation*}
Observe that using $\ell_1$-regularisation in this case imposes sparsity in the reconstructions $\bar F_K^h$ and $\bar F_K^v$ with respect to $\{H_{k}^{h}\}_{k=1}^{K}$ and $\{ H_{k}^{v}\}_{k=1}^{K}$, respectively. Also similarly as before, we note that
matrices $\mathbf{H}^h$ and $\mathbf{H}^v$ are identities when the reconstruction functions $H_k^{h}, H_k^v$ and the calibration functions $E_m^{h}, E_m^v$ are the same, otherwise they can be computed approximately from the equations 
$$
\langle H_k^h , E_{m'}^h \rangle \approx \sum_{m=1}^{M} h^h_{m,k} \langle E_m^h , E_{m'}^h \rangle, \quad \langle H_k^v , E_{m'}^v \rangle\approx \sum_{m=1}^{M} h^v_{m,k} \langle E_m^v , E_{m'}^v \rangle, \quad m'=1,\ldots,M.
$$

Finally, we note that in order to reduce the impact of noise it may be possible to include measurements of additional phase-shifts of the calibration functions. In particular, 
in addition to the calibration inputs $\mathbf{A}_m$ and $\mathbf{B}_m$ given in \R{eq:cal_funs_vf}, it may also be possible to measure  
\begin{align*}
\mathbf{C}_m :=
\begin{bmatrix}
 E_m^{h} \\
 c E_m^{v}
\end{bmatrix}\overset{\R{eq:model}}{\mapsto} 
\tilde{\mathbf{C}}_m :=
\begin{bmatrix}
\tilde{C}_m^{h} \\
\tilde{C}_m^{v}
\end{bmatrix},
\end{align*}
where
$c\neq b$ ensures that $\mathbf{C}_m\neq \mathbf{B}_m$ and,
as we will see shortly below, $c$ must be carefully chosen so that $b+c\neq2$. 
With this extra information in hand, rather than using \R{eq:rep_sys_vf}, the following functions are used to represent the complex vector fields
\begin{align}\label{eq:rep_sys_vf_A_B_C}
\bigl\{ a^* \bigl( \mathbf{A}_m - \tfrac{b^*}{2} \mathbf{B}_m - \tfrac{c^*}{2} \mathbf{C}_m \bigr), \
a \bigl( \mathbf{E}_m  - \tfrac{1}{2} \mathbf{B}_m - \tfrac{1}{2} \mathbf{C}_m \bigr) \ :\ m=1,\ldots,M \bigr\},
\end{align} 
where
$a := 1/\left(1-b/2 - c/2 \right)$ is finite given that $c$ was chosen above in such a way that $b+c\neq2$.
We then proceed as above but in place of \R{eq:pol_system} obtain 
\begin{equation}\label{eq:pol_system2}
 \mathbf{g} =  \mathbf{E} \mathbf{H} \mathbf{f} + \boldsymbol{\varepsilon}, \quad  \mathbf{E}:=
 \begin{bmatrix}
  a^*\left(\mathbf{A} \!-\! \frac{b^*}{2} \mathbf{B} \!-\! \frac{c^*}{2} \mathbf{C}\right)  & a\left(\mathbf{A} \!-\! \frac{1}{2} \mathbf{B} \!-\! \frac{1}{2} \mathbf{C}\right) 
 \end{bmatrix}, \  \mathbf{H}:=
 \begin{bmatrix}
   \mathbf{H}^h & 0 \\
   0 &  \mathbf{H}^v
 \end{bmatrix},
\end{equation}
where we now have the additional matrix $\mathbf{C}\in\mathbb{C}^{2N\times M}$ containing the outputs $\tilde{\mathbf{C}}_m$. As we will see below in Section~\ref{s:numerics}, augmenting the calibration data in such a way is indeed an effective manner to decrease the influence of the measurement noise on the reconstruction.

\section{Fourier coefficients as informative features}\label{s:fourierfeatures}

Since inhomogeneities on a cellular scale caused by cancer result in increased scattering of an optical field reflected from a tumourous tissue \cite{Drezek1999}, it is expected that they also result in higher spatial frequencies of the reflected optical field. Hence, we propose that by representing such optical fields in a Fourier basis and by inspecting the corresponding Fourier coefficients it is possible to detect the increased scattering of an optical field associated with the tumourous tissue and thereby gain insight into the disease status of the tissue.

In this section, we focus exclusively on the merits of the Fourier coefficients as indicative features of increased phase scattering. By using simulated data, we show how increased changes in phase result in a slower decay of the corresponding Fourier coefficients and how this effect can be quantified. In the next section, we confirm the findings of this section on real biological data, where we use the reconstruction framework developed in Section \ref{s:framework} to recover tissue images directly in a Fourier basis and demonstrate that the recovered Fourier coefficients are indeed useful for detecting cancer.

\subsection{Fourier coefficients of a one-dimensional simulated example}

\begin{figure}[h]
	\centering
 \adjincludegraphics[width=1\textwidth,trim={{.1\width} {.00\width} {.09\width} {.0\width}},clip]{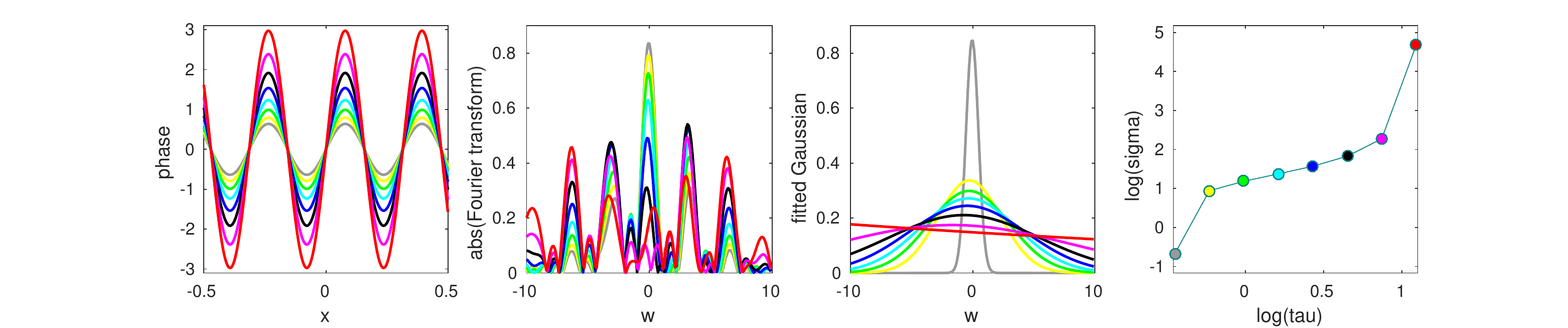} 
	\caption{Higher phase oscillations result in a slower decay of the Fourier coefficients.}
	\label{fig:1dsim}
\end{figure}

We first consider a simple one-dimensional example to illustrate the effect of increased phase oscillations on the decay of the corresponding Fourier coefficients. In this example, we measure the decay of the Fourier coefficients associated with different functions $F^{(j)}(x)=R(x)\exp({\I P^{(j)}(x)})$, $x\in I$, $j=1,\ldots,8$, with the same amplitude $R$ but with different phases $P^{(j)}$ defined on the compact interval $I:=[-1/2,1/2]$. In particular, for illustration purposes we take $R(x):=\exp(-x^2)$ and  $P^{(j)}(x):=\tau^{(j)}\sin(20x)$, where $0<\tau^{(1)}<\cdots<\tau^{(8)} < 2\pi$, so that different phase functions exhibit different degrees of oscillations. These phase functions are shown in the first panel of Figure~\ref{fig:1dsim}. 
Since the Fourier basis in the domain $I$ is given by $\{e^{2\pi\I k x}\}_{k\in\mathbb{Z}}$, then for each $F^{(j)}$ we compute its first 20 Fourier coefficients as 
$$
f^{(j)}_k:=  \int_{I} F^{(j)}(x)e^{-2\pi\I kx}\D x,
$$
where $k=-10,\ldots,9$, and approximate its Fourier transform by the classical Whittaker--Shannon interpolation formula $\sum_{k}f^{(j)}_k\mathrm{sinc}(w-k)$, $w\in\mathbb{R}$. The amplitude of the approximated Fourier transform of each $F^{(j)}$ is shown in the second panel of Figure~\ref{fig:1dsim}. 
Finally, we quantify the decay of  the Fourier coefficients  by the standard deviation $\sigma^{(j)}$ of a Gaussian function $a^{(j)}\exp(-(w-c^{(j)})^2/(2(\sigma^{(j)})^2))$ fitted to the amplitude of the approximated Fourier transform on interval $w\in[-10,10)$. The fitted Gaussian functions are shown in the third panel of Figure~\ref{fig:1dsim}.
From the forth panel of Figure~\ref{fig:1dsim}, we observe that an increased magnitude of the phase oscillations $\tau^{(j)}$ results in an increased standard deviation $\sigma^{(j)}$, which can be used to quantify the decay of the Fourier coefficients. 
It is important to note that
although in this example the zeros of the different phase functions coincide,  the same effect is observed even if this is not the case. 
Also,  if the frequency of the phase oscillation is increased while their magnitude is kept constant, then $\sigma^{(j)}$ would increase as well.

The takeaway message from this simple one-dimensional example is that representing a signal with respect to a Fourier basis is especially useful to identify variations in oscillating phase, and that the decay of the corresponding Fourier coefficients is sensitive to variations in phase scattering in a manner that can be easily identified. As we will see in the remainder of the paper, these observations remain true also in higher dimensional practical examples.

\subsection{Fourier coefficients of simulated tissue images}

\begin{figure}[htbp!]
\begin{minipage}{.75\textwidth}
	\centering
	\hspace{0.5cm} \small{Phase} \hspace{2.4cm} \small{abs(FT)} \hspace{1.8cm} \small{Fitted Gaussian} \\
1	\adjincludegraphics[scale=0.5,trim={{.05\width} {.05\width} {.05\width} {.00\width}},clip]{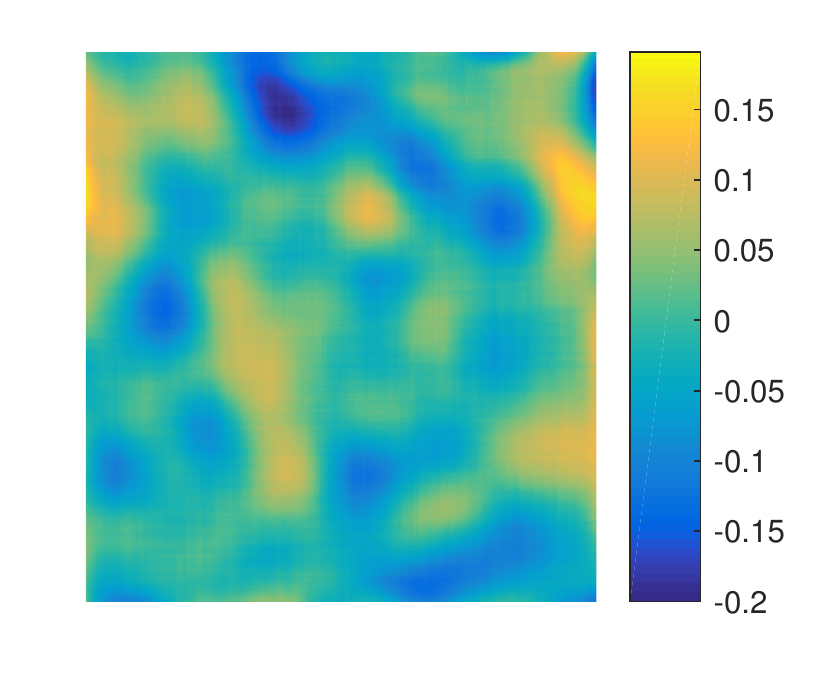}  \adjincludegraphics[scale=0.5,trim={{.05\width} {.05\width} {.05\width} {.00\width}},clip]{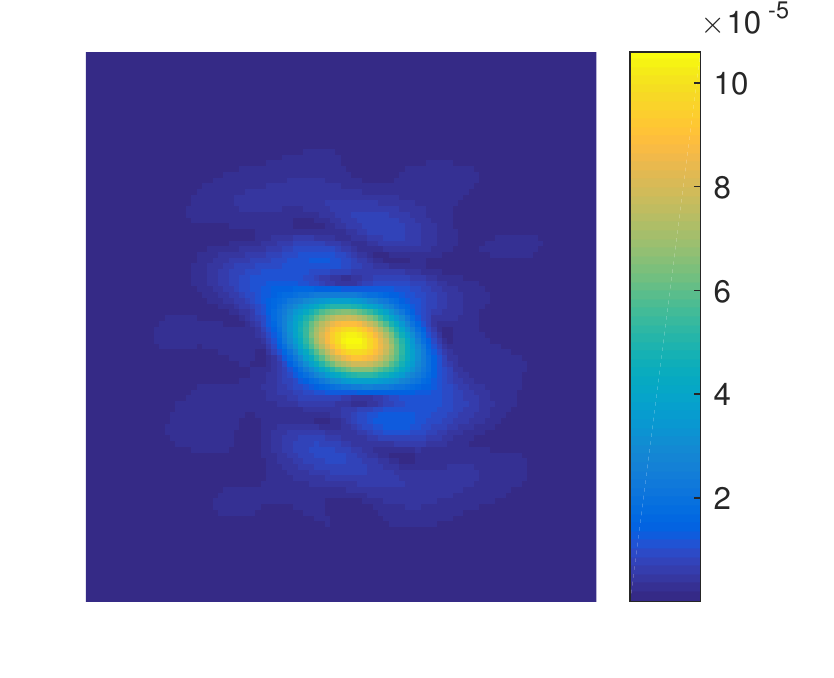} 
	\adjincludegraphics[scale=0.5,trim={{.05\width} {.05\width} {.05\width} {.00\width}},clip]{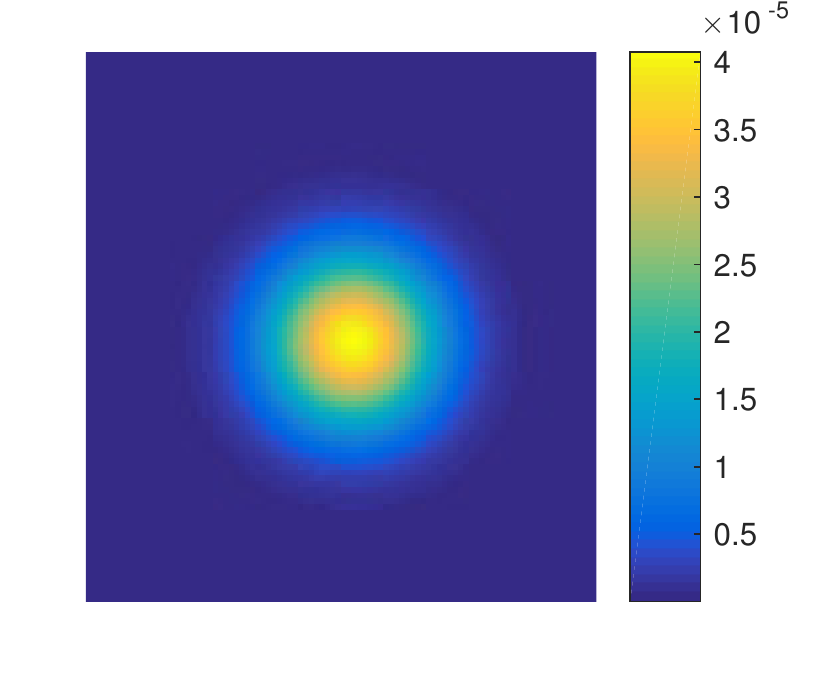}\\
2		\adjincludegraphics[scale=0.5,trim={{.05\width} {.05\width} {.05\width} {.00\width}},clip]{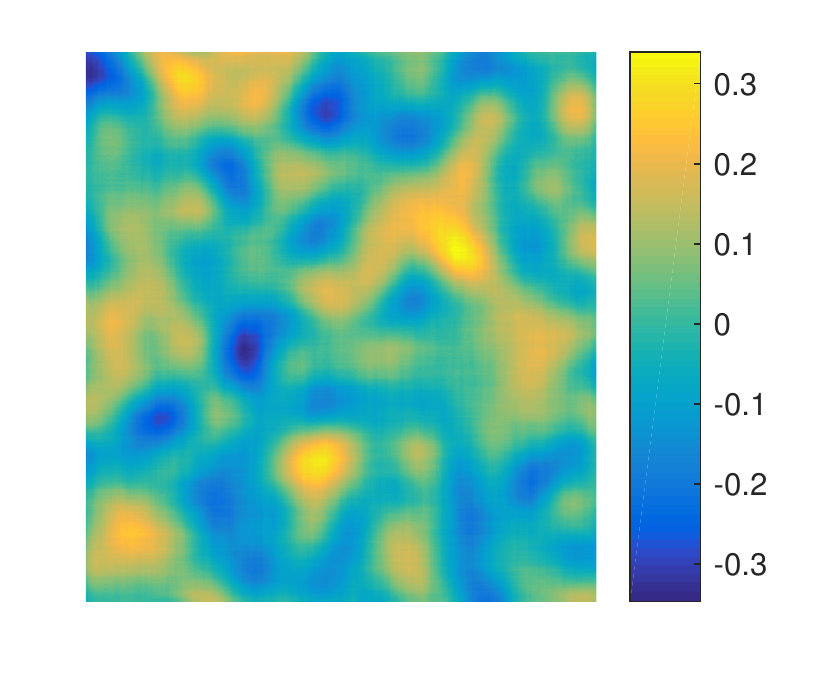}  \adjincludegraphics[scale=0.5,trim={{.05\width} {.05\width} {.05\width} {.00\width}},clip]{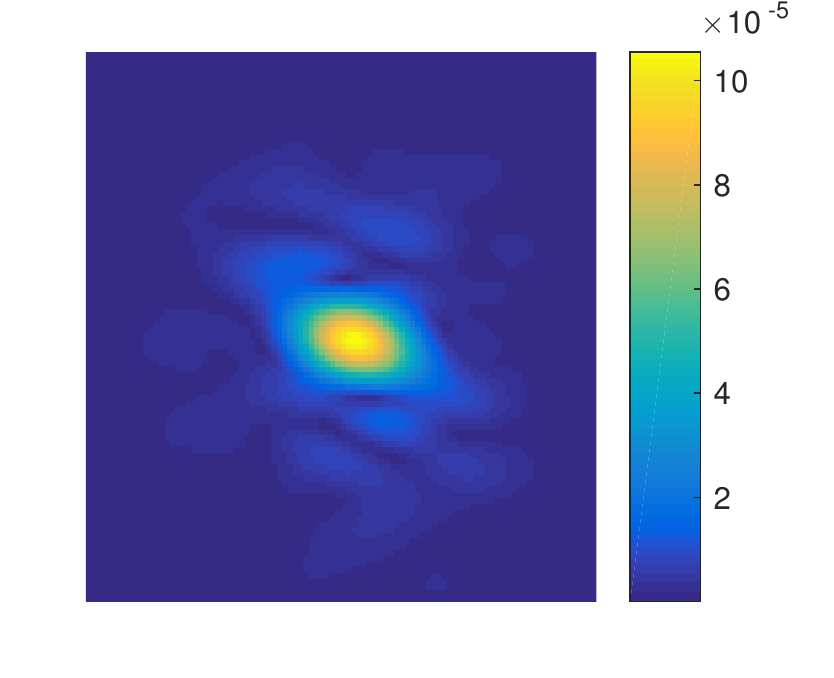} 
	\adjincludegraphics[scale=0.5,trim={{.05\width} {.05\width} {.05\width} {.00\width}},clip]{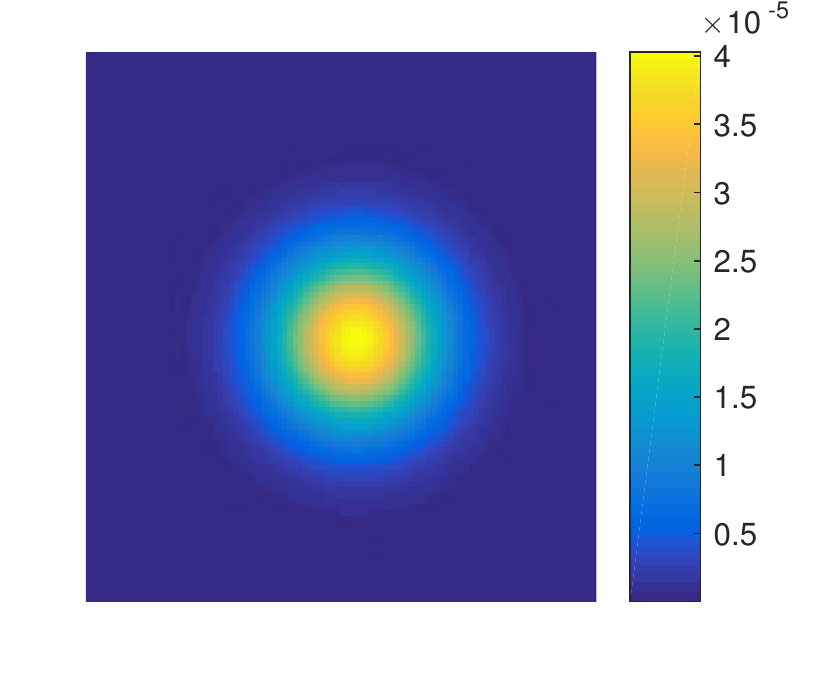}\\
3		\adjincludegraphics[scale=0.5,trim={{.05\width} {.05\width} {.05\width} {.00\width}},clip]{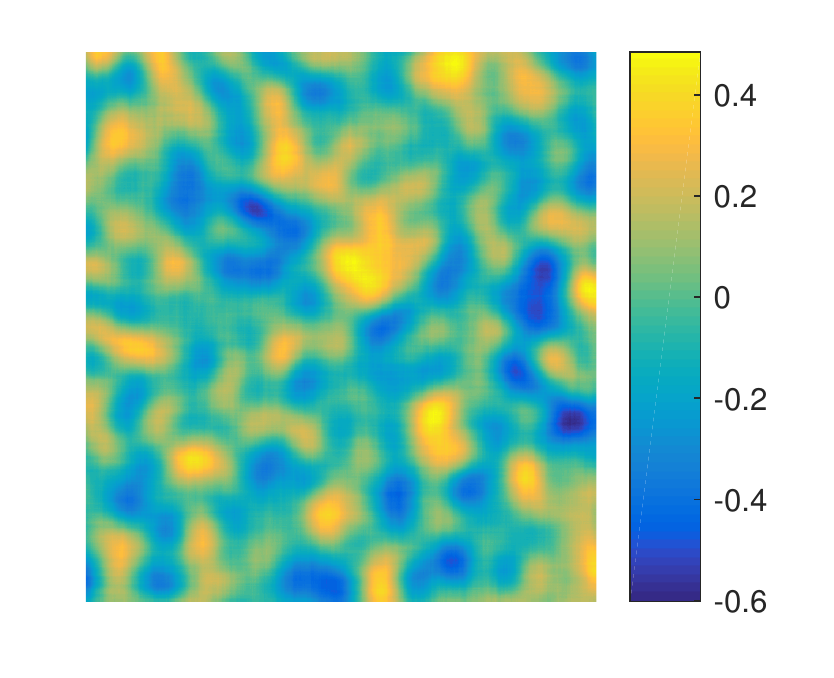}  \adjincludegraphics[scale=0.5,trim={{.05\width} {.05\width} {.05\width} {.00\width}},clip]{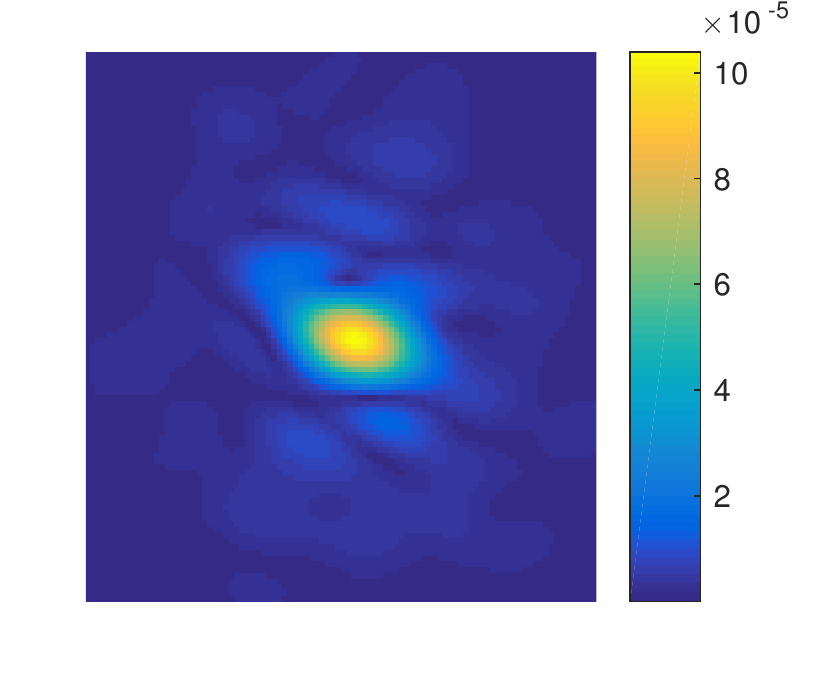} 
	\adjincludegraphics[scale=0.5,trim={{.05\width} {.05\width} {.05\width} {.00\width}},clip]{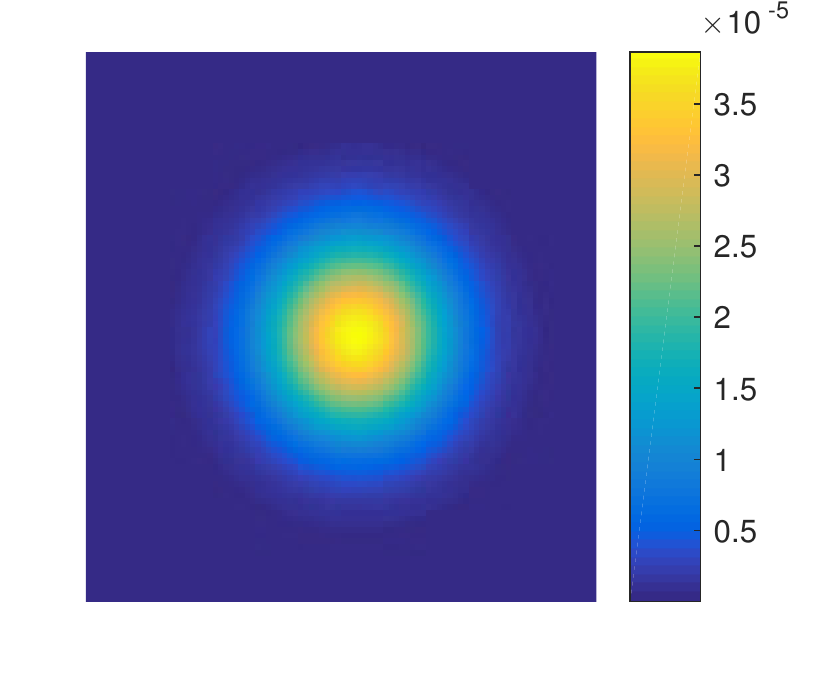}\\
4		\adjincludegraphics[scale=0.5,trim={{.05\width} {.05\width} {.05\width} {.00\width}},clip]{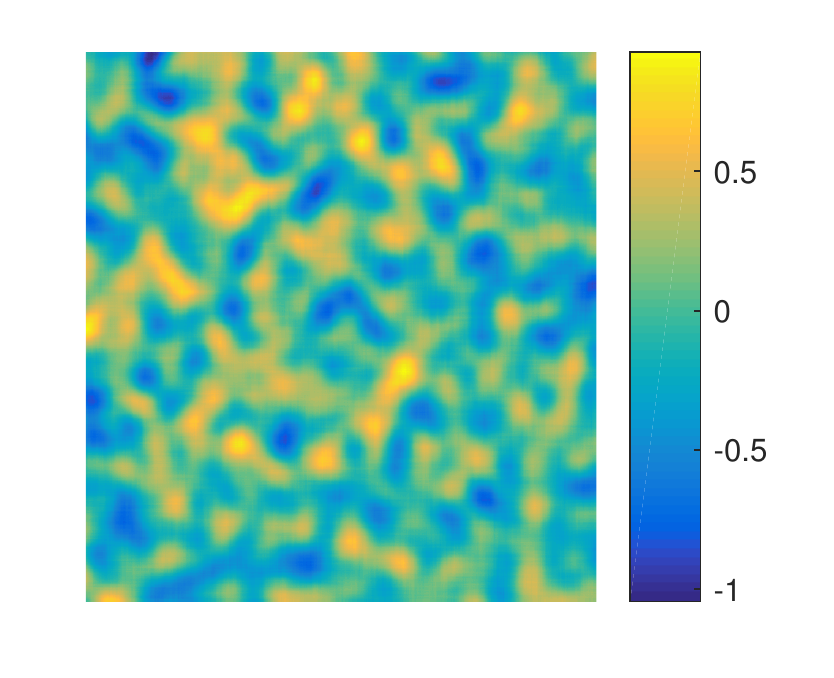}  \adjincludegraphics[scale=0.5,trim={{.05\width} {.05\width} {.05\width} {.00\width}},clip]{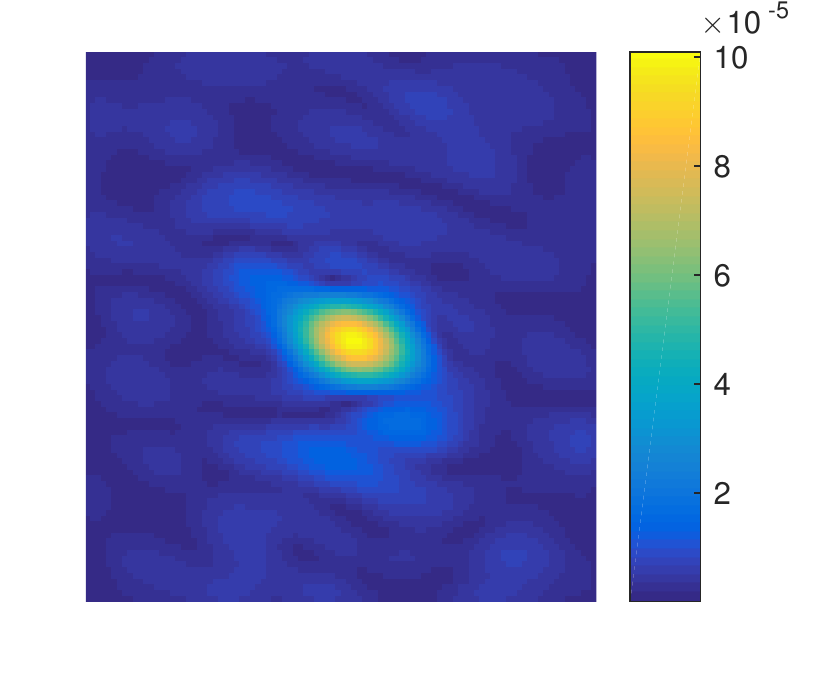} 
	\adjincludegraphics[scale=0.5,trim={{.05\width} {.05\width} {.05\width} {.00\width}},clip]{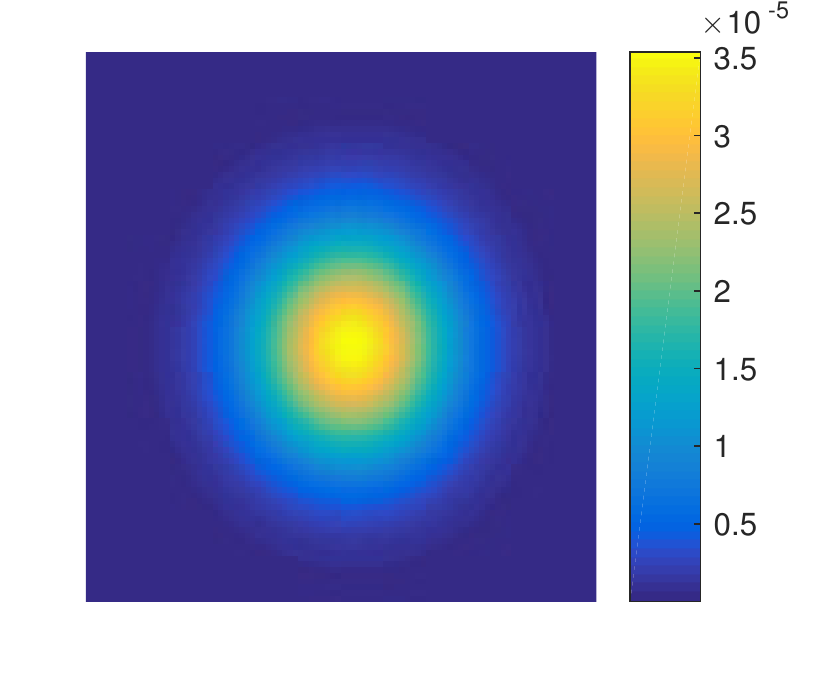}\\
5		\adjincludegraphics[scale=0.5,trim={{.05\width} {.05\width} {.05\width} {.00\width}},clip]{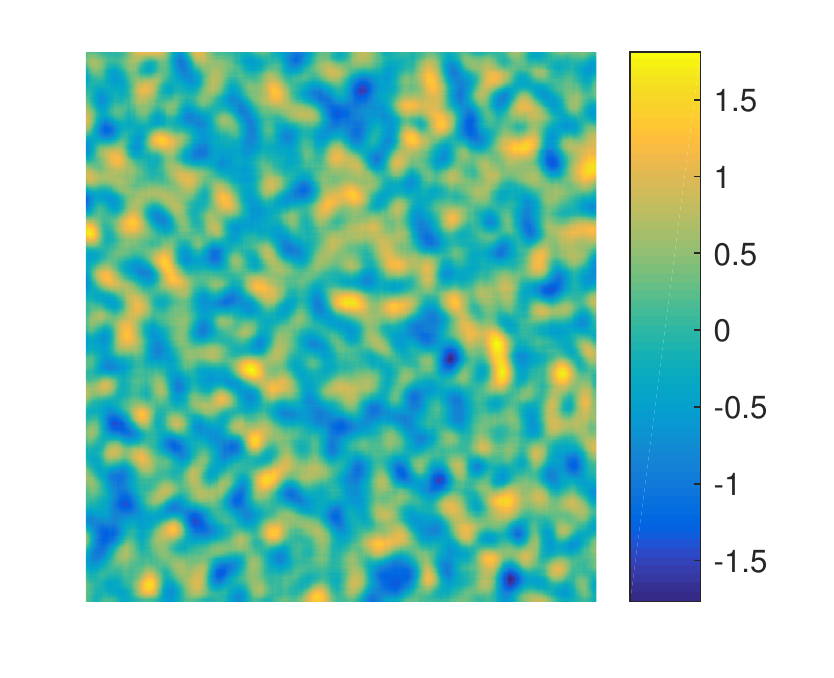}  \adjincludegraphics[scale=0.5,trim={{.05\width} {.05\width} {.05\width} {.00\width}},clip]{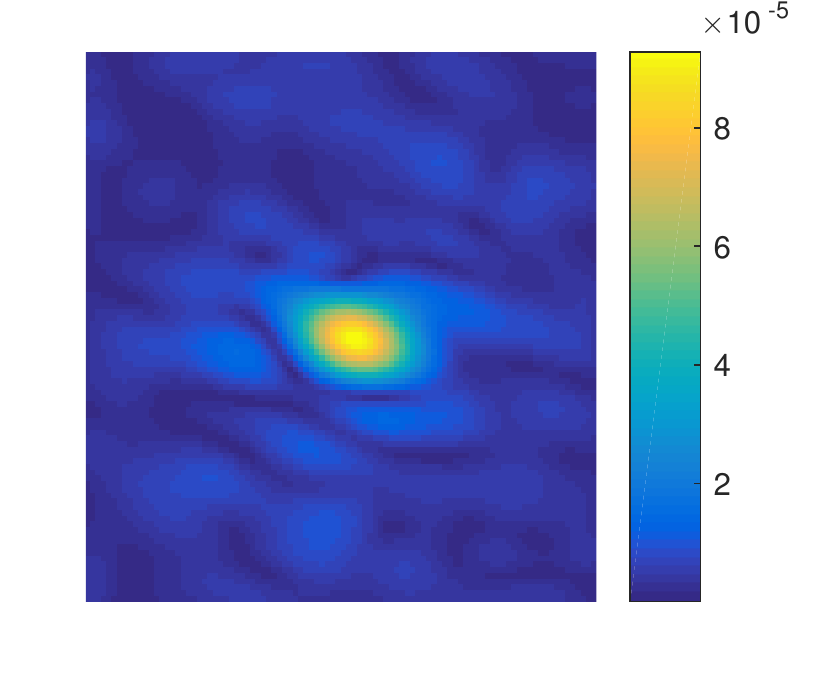} 
	\adjincludegraphics[scale=0.5,trim={{.05\width} {.05\width} {.05\width} {.00\width}},clip]{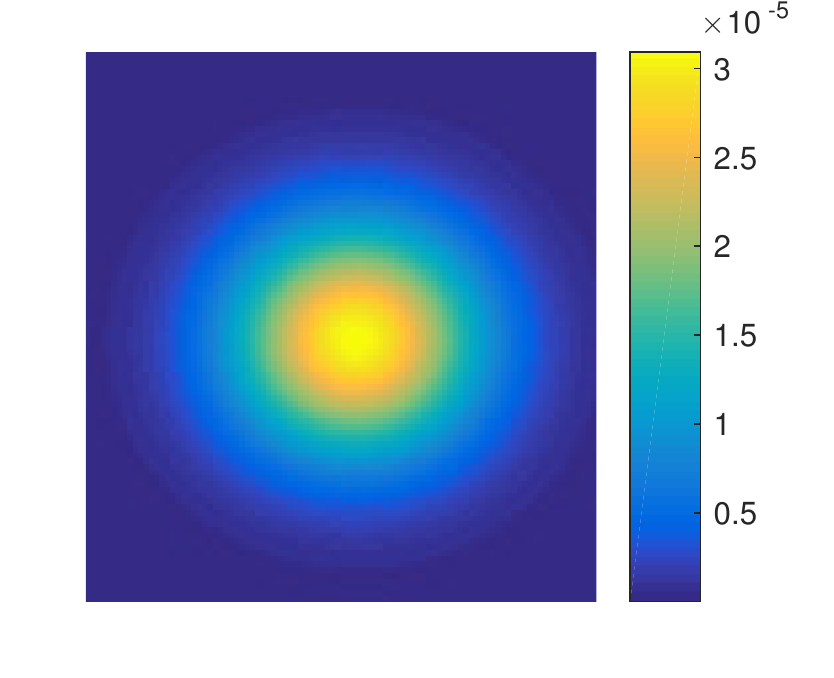}\\
6		\adjincludegraphics[scale=0.5,trim={{.05\width} {.05\width} {.05\width} {.00\width}},clip]{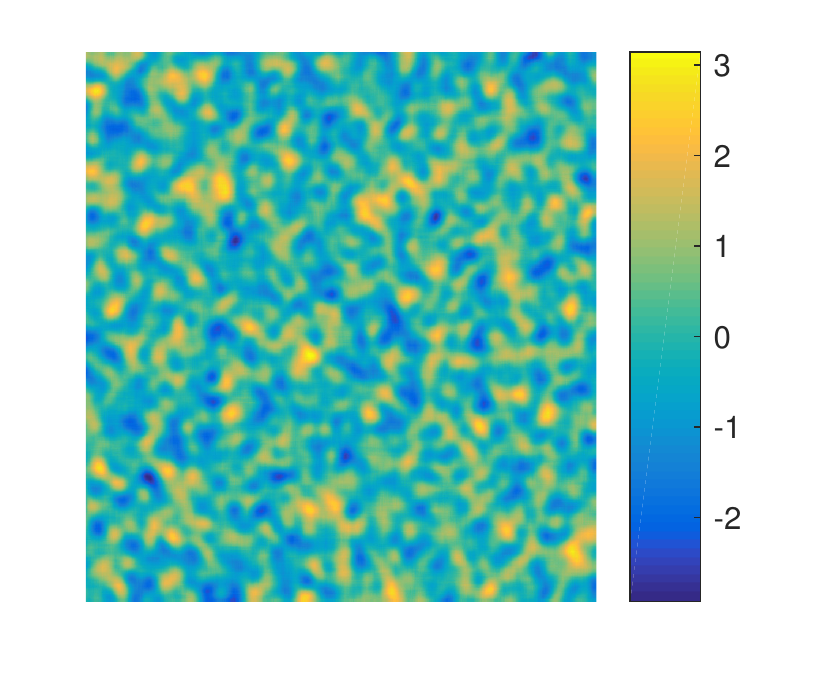}  \adjincludegraphics[scale=0.5,trim={{.05\width} {.05\width} {.05\width} {.00\width}},clip]{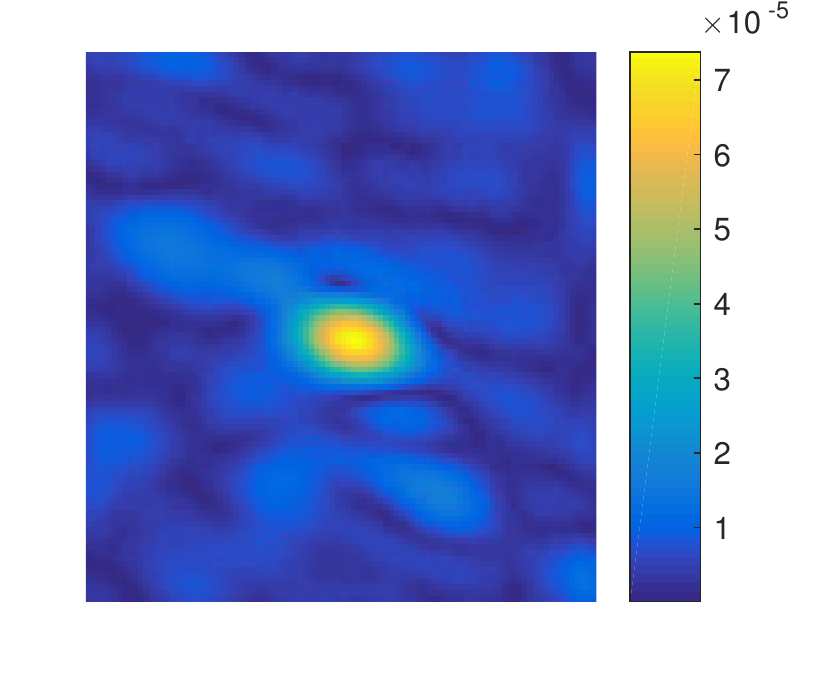} 
	\adjincludegraphics[scale=0.5,trim={{.05\width} {.05\width} {.05\width} {.00\width}},clip]{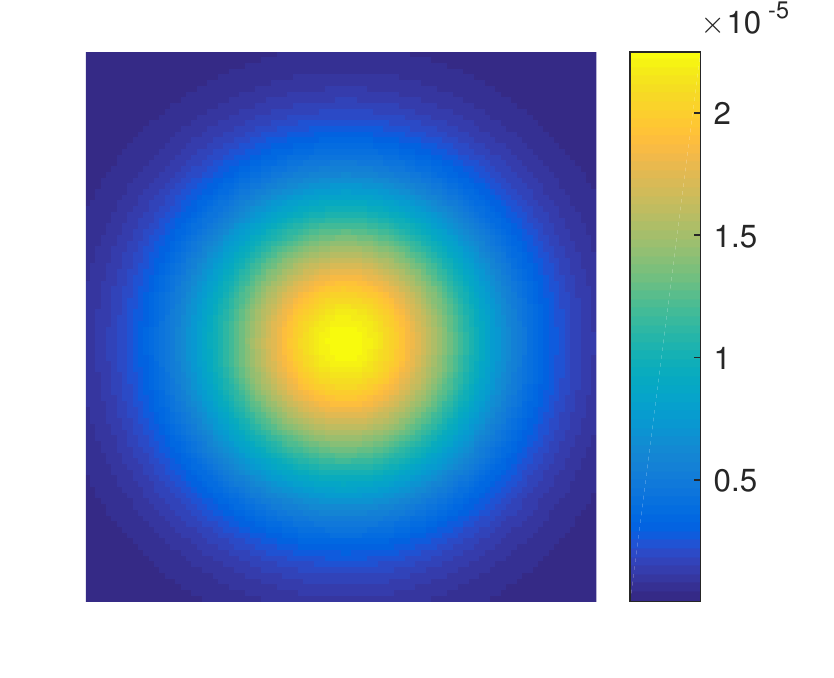}
	\end{minipage}%
\begin{minipage}[t]{.25\textwidth}
\centering
\small{Amplitude}\\
\adjincludegraphics[scale=0.5,trim={{.05\width} {.05\width} {.05\width} {.00\width}},clip]{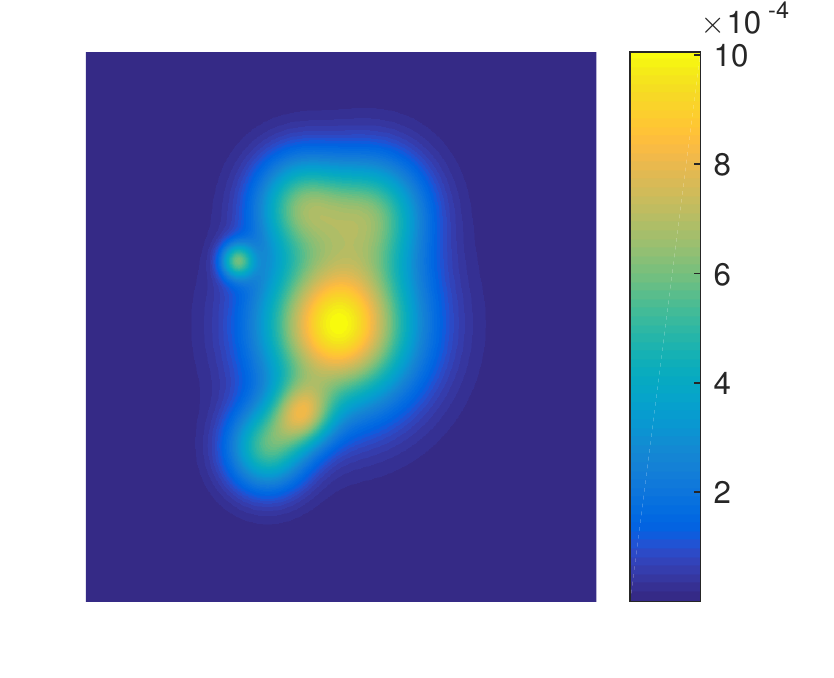}\\
\vspace{0.5cm}
\small{$\log(\sigma_1^{(j)}+\sigma_2^{(j)})$}\\
\adjincludegraphics[scale=0.49,trim={{.00\width} {.02\width} {.02\width} {.02\width}},clip]{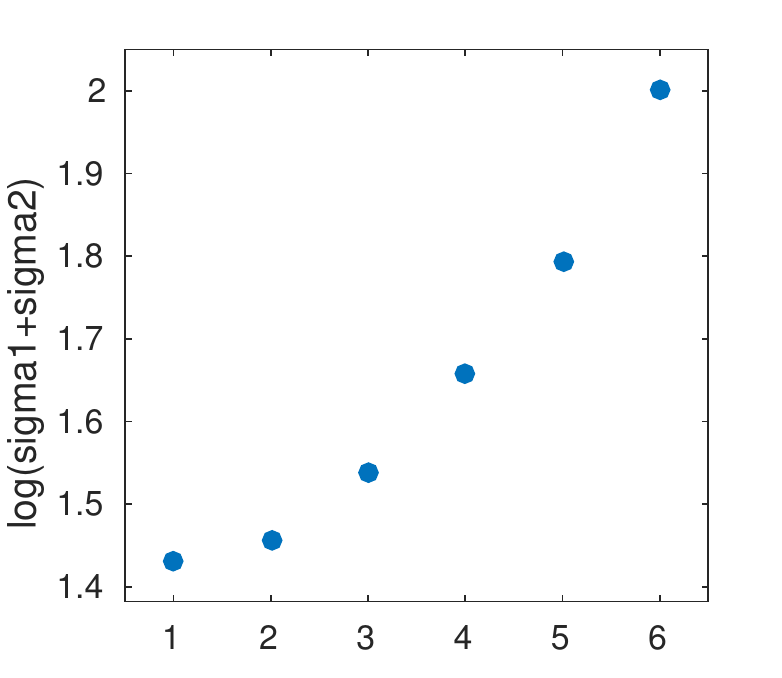}\\
{\vspace{-0.4cm}}{\hspace{3.1cm} \footnotesize{$j$}}
\end{minipage}
	\caption{Six simulated images with the same amplitude but different phase, which are generated from our model with increasing $\tau^{(j)}/\rho^{(j)}$, $j=1,\ldots,6$, so that larger $\tau^{(j)}/\rho^{(j)}$ characterises larger phase oscillations. In the scatter plot on the right, we report the sum of parameters $\sigma_1^{(j)}$ and $\sigma_2^{(j)}$ of the Gaussian fitted to the amplitude of the Fourier transform abs(FT), revealing that increased $\tau^{(j)}/\rho^{(j)}$ correlates with larger $\sigma_1^{(j)}+\sigma_2^{(j)}$.}
	\label{fig:SimPhase}
\end{figure}

We now generalise our observations to two-dimensional complex-valued functions. In particular, we perform a simulation study to demonstrate that higher and more frequent changes in phase result in a slower decay of the Fourier coefficients of the corresponding function.  For this purpose, we create a model mimicking tissue samples with a different level of phase oscillations, which we then use to generate images and compute their Fourier coefficients. 

In our model, we use randomness to achieve certain variability across different samples and two different parameters to control the degree of phase oscillations. In detail, the model that we use in our simulation study corresponds to a complex function $F(\mathbf{x}):=R(\mathbf{x})\exp(\I P(\mathbf{x}))$, $\mathbf{x}\in S$, where the original space-domain $S:=[-1/2,-1/2]^2$ is discretized into a $700\times700$ grid, while $R$ and $P$ are chosen randomly as we now describe. The phase function $P:=P^{(\tau,\rho)}$ depends on two given parameters $\tau$ and $\rho$, controlling the amplitude and the frequency of phase oscillations, respectively. Specifically, $800\times800$ pixel-values are chosen uniformly at random from $[-1,1]$,  which are then filtered by using MATLAB's function `imgaussfilt' with the smoothing parameter $\rho$. Following this step, only $700\times700$ pixels are kept by removing $50$ pixels from each boundary and such image is then rescaled so that all phase pixel-values are between $[-\tau,\tau]$, $\tau\in[0,\pi]$. The amplitude function $R$ is selected as the sum of $\exp(-50\|\mathbf{x}\|_2^2)/1000$ and five additional Gaussian functions $\exp(-\|\mathbf{x}-\mathbf{c}\|_2^2/d)/2000$ with randomly chosen parameters $\mathbf{c}$ and $d$.

In Figure~\ref{fig:SimPhase} we demonstrate how changing phase parameters $\tau$ and $\rho$ while keeping amplitude fixed changes  the decay of Fourier coefficients. Specifically, we use six different values $(\rho^{(j)},\tau^{(j)})$, $j=1,\ldots,6$ to create six different functions $F^{(j)}$, where $0<\tau^{(1)}<\cdots<\tau^{(6)}\leq\pi$ and $0.025<(\rho^{(1)})^{-1}<\cdots<(\rho^{(6)})^{-1}\leq0.125$ are increasing logarithmically. 
Similarly to the one-dimensional example of Figure~\ref{fig:1dsim}, the decay of corresponding Fourier coefficients is measured by standard deviation of a Gaussian function $a\exp(-(x_1-c_1)^2/(2\sigma_1^2)-(x_2-c_2)^2/(2\sigma_2^2))$ fitted to the absolute value of the Fourier transform approximated from the first $20\times 20$ Fourier coefficients, where the Fourier coefficients $\{f_{\mathbf{k}}\}_{\mathbf{k}\in I_{400}}$ of function $F$ are computed using the formula in \R{eq:fourierceof}. 
In particular, in Figure~\ref{fig:SimPhase}, for each $F^{(j)}$ we report the sum of the standard deviations $\sigma_1^{(j)}+\sigma_2^{(j)}$ of the fitted Gaussian, thereby observing that increased phase oscillations, i.e.~increased $\tau^{(j)}/\rho^{(j)}$, results in slower decay of the corresponding Fourier coefficients, i.e.~larger $\sigma_1^{(j)}+\sigma_2^{(j)}$.

\begin{figure}[htbp!]
\begin{minipage}[c]{0.5\textwidth}
\centering
\adjincludegraphics[width=0.8\textwidth,trim={{.00\width} {.02\width} {.02\width} {.02\width}},clip]{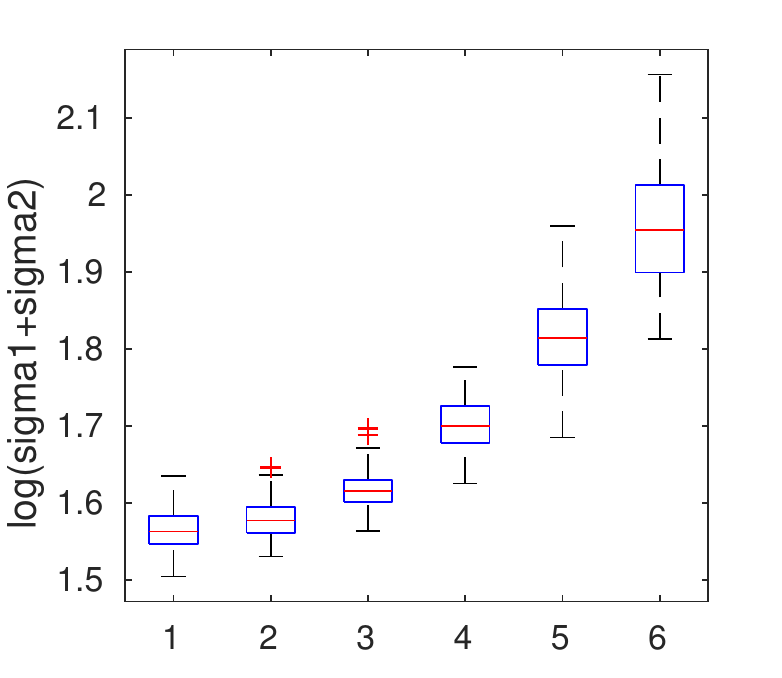} \\
{\vspace{-0.4cm}}{\hspace{5.4cm} \small{$j$}}
\end{minipage}\hfill
\begin{minipage}[c]{0.5\textwidth}
\caption{For each of the six categories, we generated a hundred images with different phase and amplitude from the tissue model with fixed parameters $\tau^{(j)}$ and $\rho^{(j)}$, $j=1,\ldots,6$. For each image, we then approximated its Fourier transform and fitted a Gaussian function whose parameter $\sigma_1^{(j)}+\sigma_2^{(j)}$ is reported. }
\label{fig:SimPhase2}
\end{minipage}
\end{figure}

Next, in Figure~\ref{fig:SimPhase2}, for each of the six different values $(\rho^{(j)},\tau^{(j)})$, $j=1,\ldots,6$, chosen as in Figure~\ref{fig:SimPhase}, we generate a hundred  different images using our model (with each image having a different phase and a different amplitude) and we report the value $\sigma_1^{(j)}+\sigma_2^{(j)}$ of the fitted Gaussian.  We observe the same trend in the decay of the Fourier coefficients in Figure~\ref{fig:SimPhase2} as in Figure~\ref{fig:SimPhase}, but now across 600 different images.

We conclude this section by noting that the features extracted from Fourier coefficients as we described in this section have three additional useful properties. First, since the amplitude of the Fourier transform is invariant to the shifts of the corresponding complex function in its space-domain, the features that we extract are invariant to the shifts of the tissue images in their space-domain. Second, the quality of the recovered phase in the space-domain is dependent on a phase unwrapping procedure and is thus highly sensitive to noise, which means that phase may bear more information in the Fourier-domain than in the space-domain. Third, once the Fourier coefficients are computed from the available measurements, each image can easily be represented in both the Fourier and the original space-domain, allowing for additional flexibility.

\section{Experimental results}\label{s:numerics}

Having established the utility of Fourier coefficients in quantifying phase scattering using simulated data in Section~\ref{s:fourierfeatures}, we now apply the reconstruction framework developed in Section \ref{s:framework} to measurements obtained experimentally using the prototype fibre endoscope developed in \cite{Gordon2018}, which can measure optical phase, polarisation and amplitude. In Section~\ref{ss:holographic}, we first demonstrate the recovery of a synthetic holographic image with a known ground-truth that can be used for validation. Next, in Section~\ref{ss:biological}, we apply our reconstruction framework to biological images of tissue samples taken from mice and demonstrate that reconstruction with respect to a Fourier basis can be used as a diagnostic indicator of early tumorigenesis.

\subsection{Reconstruction of a synthetic holographic image}\label{ss:holographic}

\begin{figure}[H]
	\centering
   \begin{tabular}{ccccc}
   \multicolumn{2}{c}{\small{Holographic image at the distal end}} & & \multicolumn{2}{c}{\small{Raw output at the proximal end}} \\
	\small{$\text{abs}(F^h)$} & \small{$\text{phase}(F^h)$} & & \small{$\text{abs}(\tilde F^h)$} & \small{$\text{phase}(\tilde F^h)$} \\
	\includegraphics[scale=0.38]{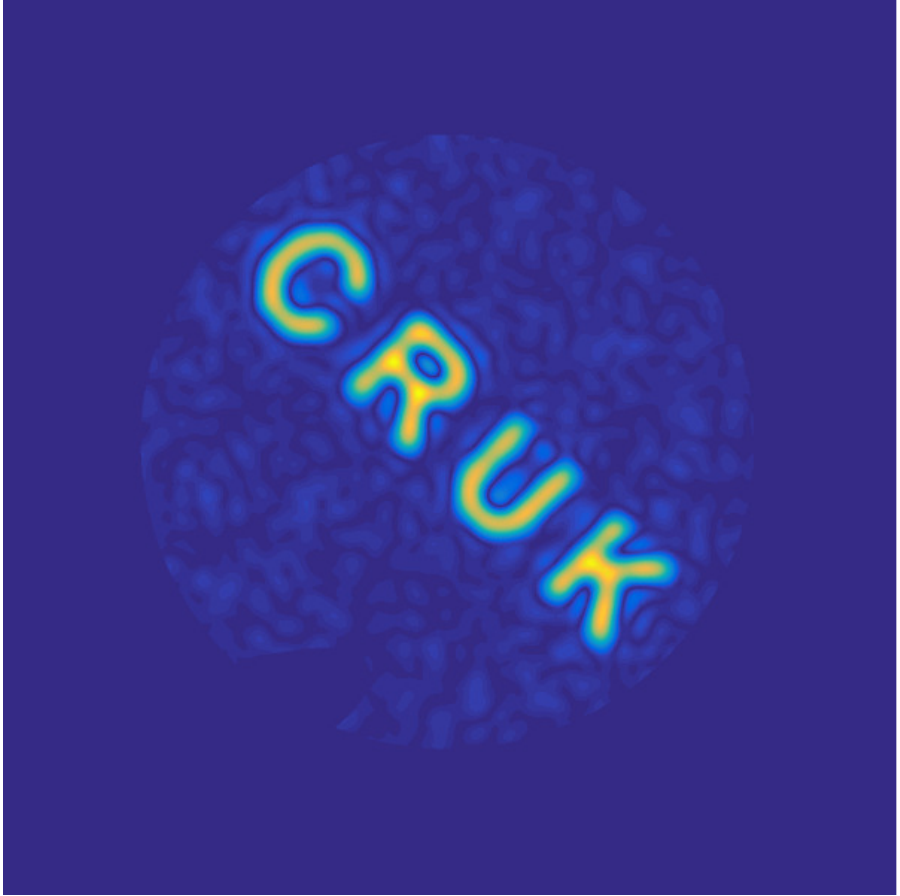} & \includegraphics[scale=0.38]{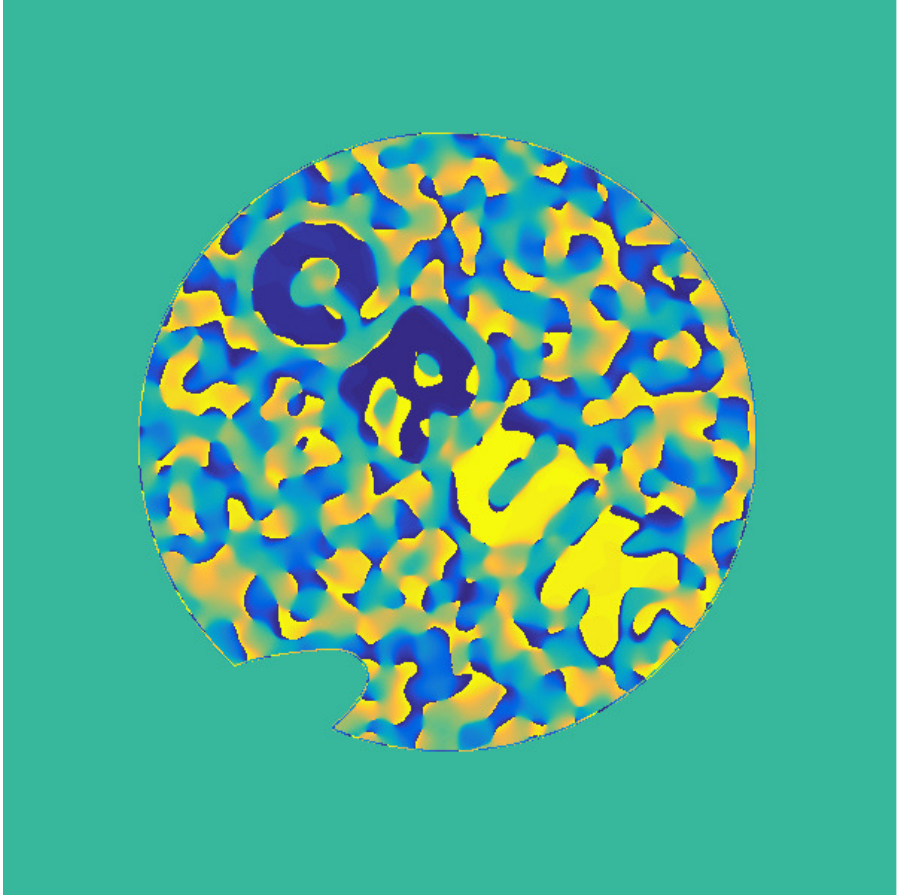} & \hspace{0.01cm} &
	\includegraphics[scale=0.38]{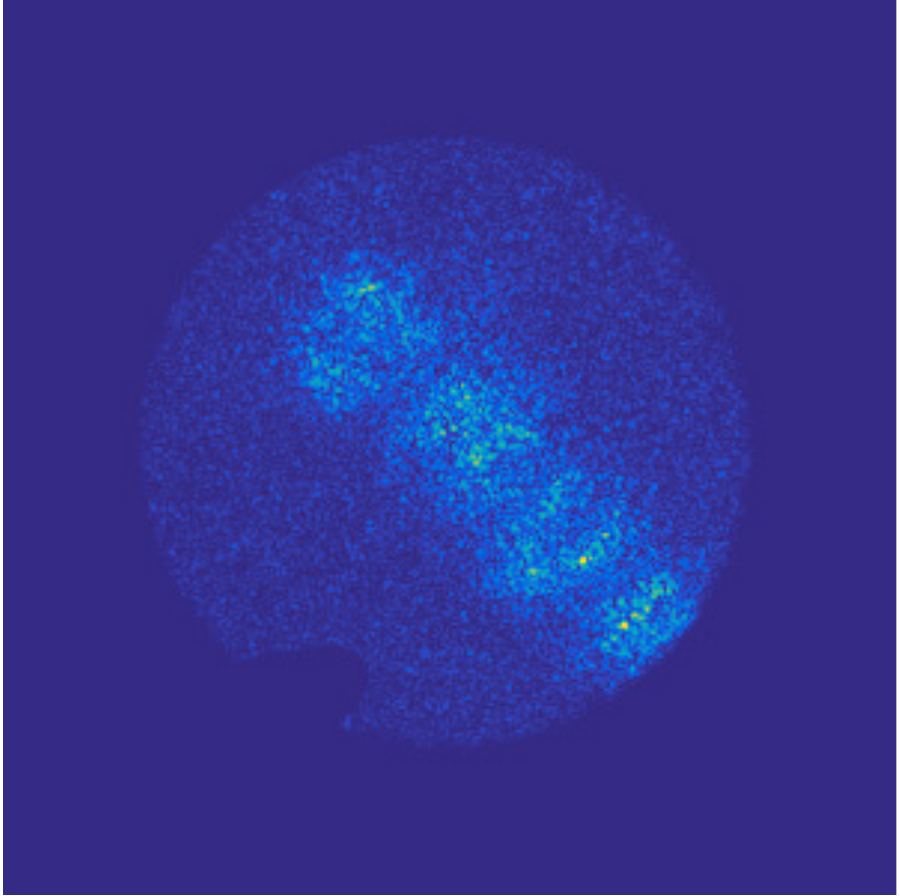} & \includegraphics[scale=0.38]{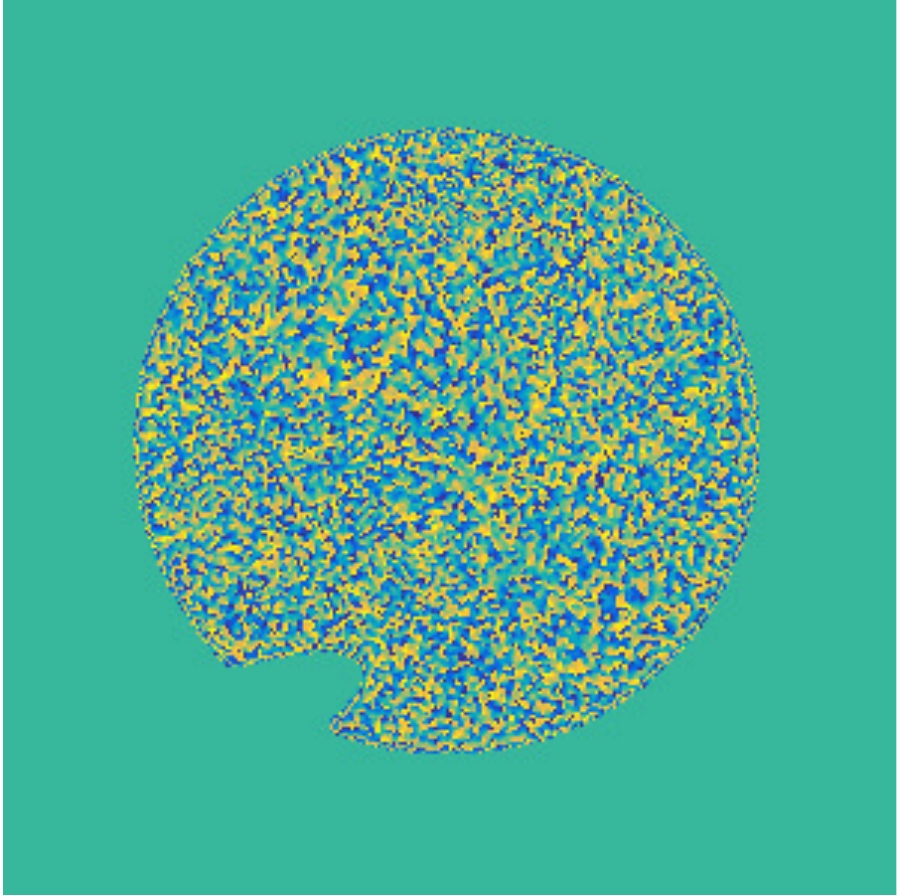}\\
	\small{$\text{abs}(F^v)$} & \small{$\text{phase}(F^v)$} & & \small{$\text{abs}(\tilde F^v)$} & \small{$\text{phase}(\tilde F^v)$} \\
	\includegraphics[scale=0.38]{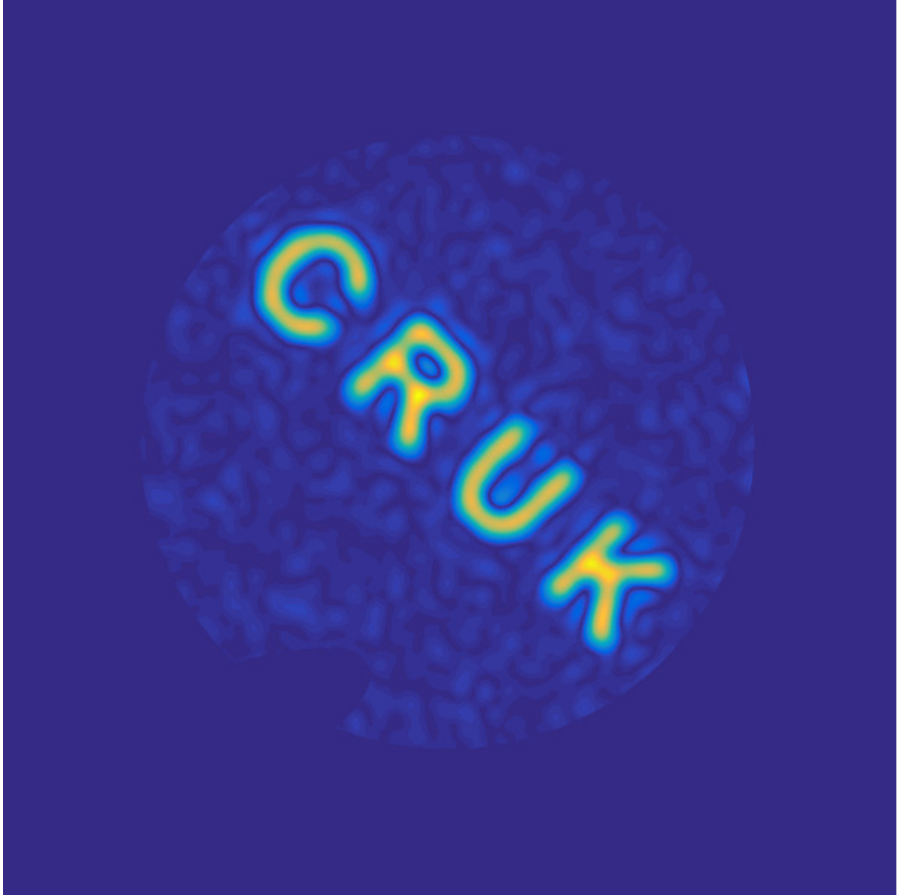} & \includegraphics[scale=0.38]{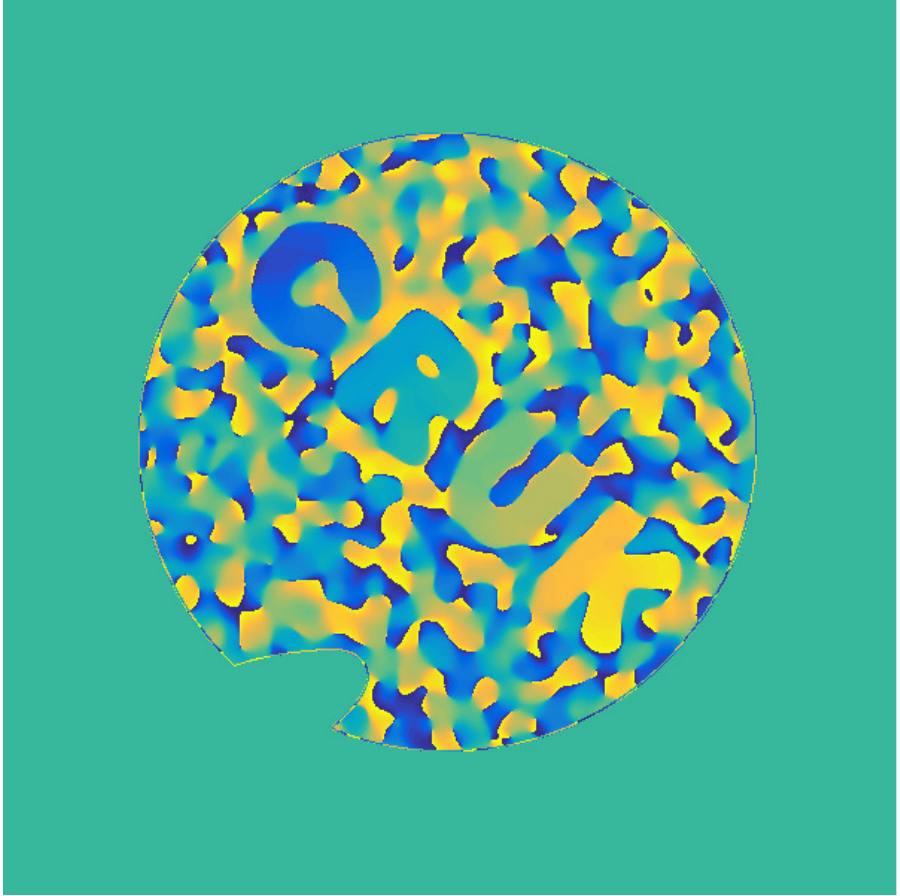} & \hspace{0.01cm} &
	\includegraphics[scale=0.38]{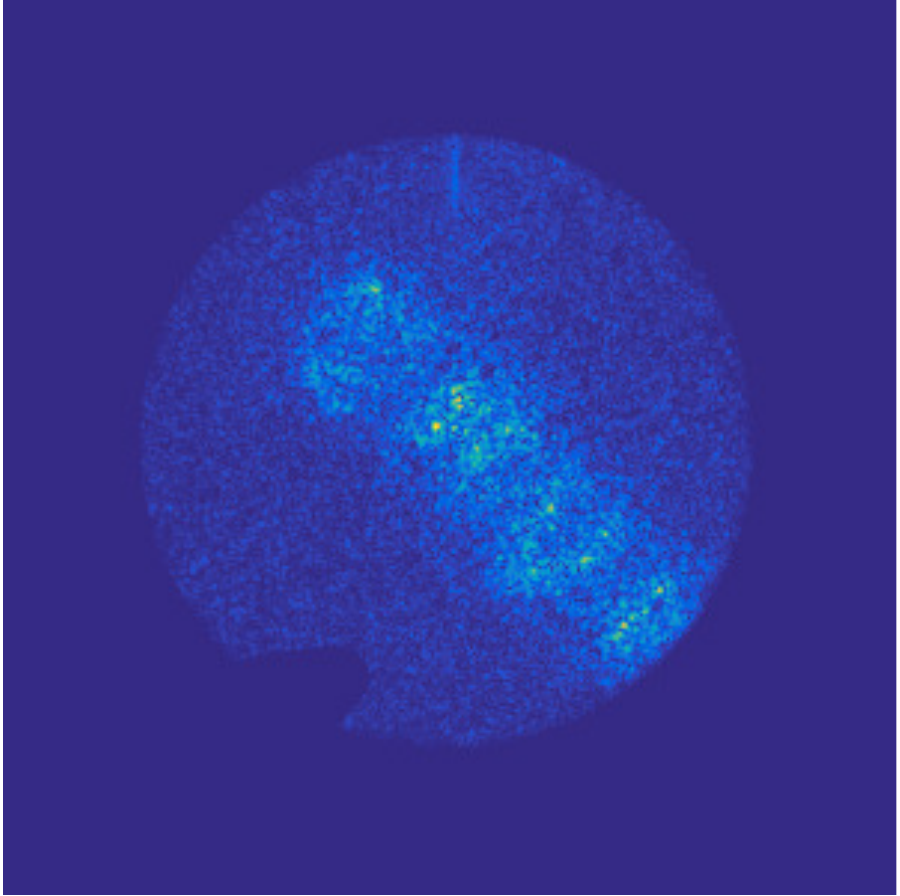} & \includegraphics[scale=0.38]{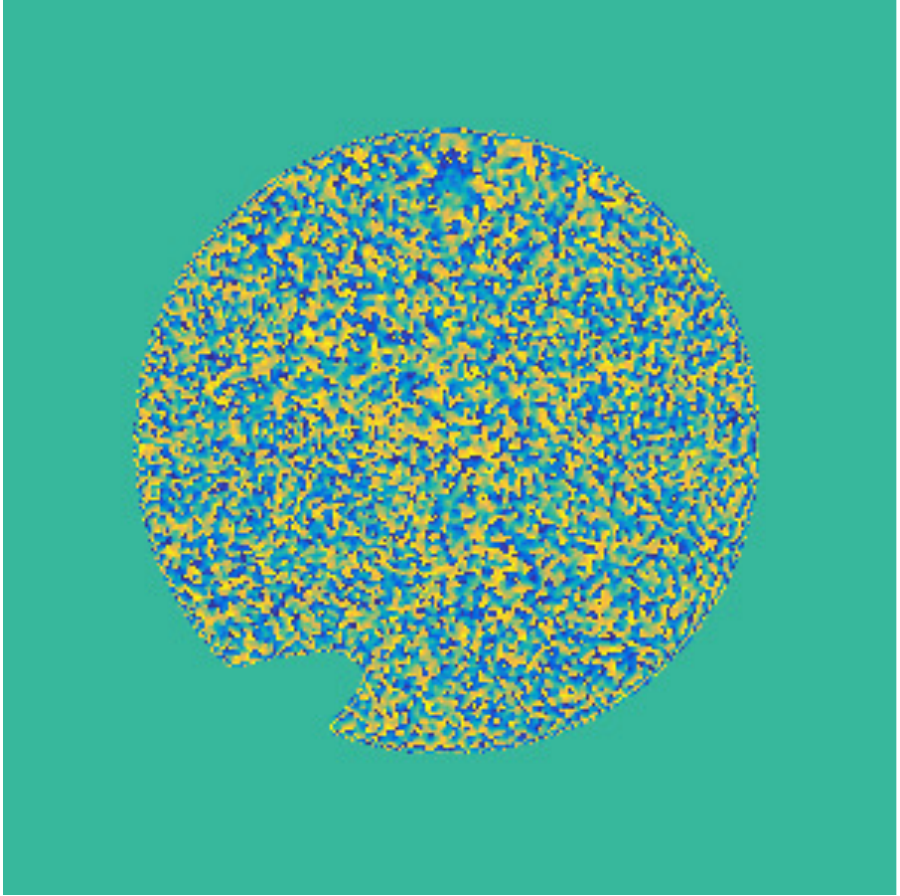}
	\end{tabular}\vspace{-0.1cm}
	\caption{Amplitude and phase of the horizontal and vertical polarisation of the ground-truth synthetic holographic image at the distal end and the corresponding output at the proximal end.}
	\label{fig:CRUKdata}
\end{figure}

To demonstrate our reconstruction algorithm on experimental data we first reconstruct a synthetic holographic image.  Since in this case we have access to the ground-truth image at the distal end of the fibre, we can visually assess the quality of our proposed imaging methodology. 
Specifically, in this subsection,
using the raw output of the synthetic holographic image shown in Figure~\ref{fig:CRUKdata}, we  test our general reconstruction framework in combination with different representation systems as well as different regularisation terms.

\begin{figure}[htb!]
	\centering
   \begin{tabular}{ccccc}
   \multicolumn{2}{c}{\small{Distal end}} & & \multicolumn{2}{c}{\small{Proximal end}} \\
	\small{$\text{abs}(E_m^h)$} & \small{$\text{phase}(E_m^h)$} & & \small{$\text{abs}(\tilde A_m^h)$} & \small{$\text{phase}(\tilde A_m^h)$} \\
	\includegraphics[scale=0.37]{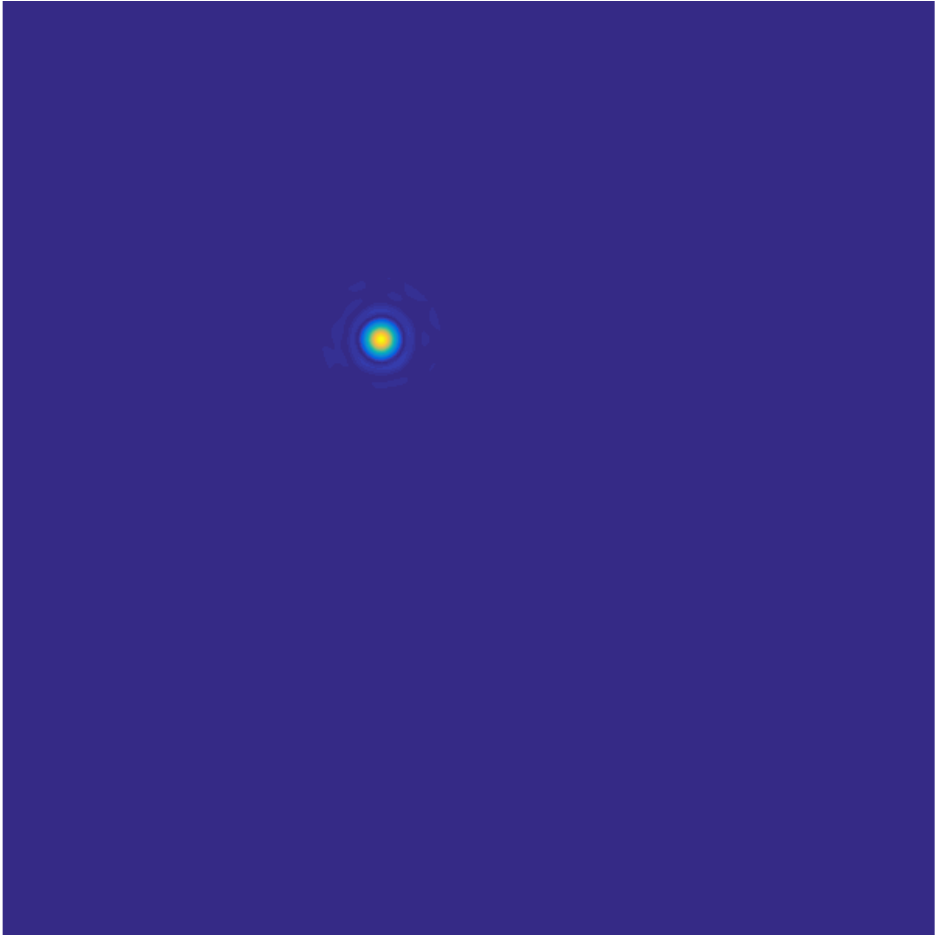} & \includegraphics[scale=0.37]{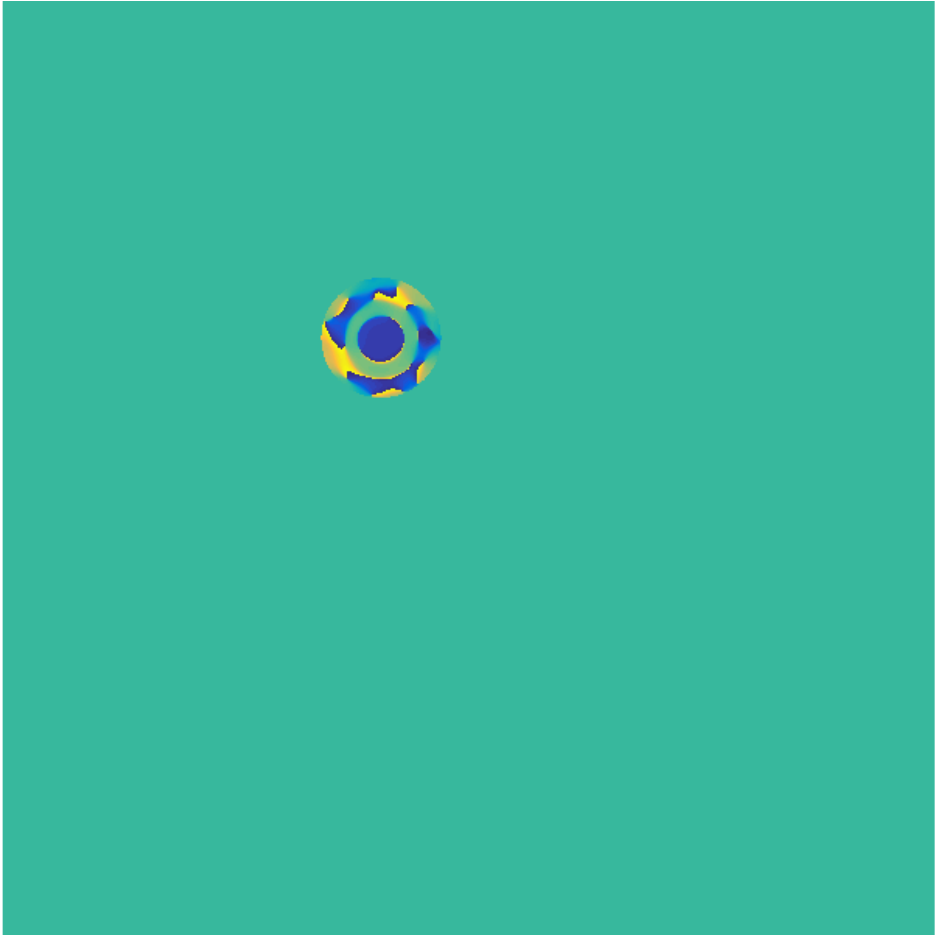} & \hspace{0.01cm} &
	\includegraphics[scale=0.37]{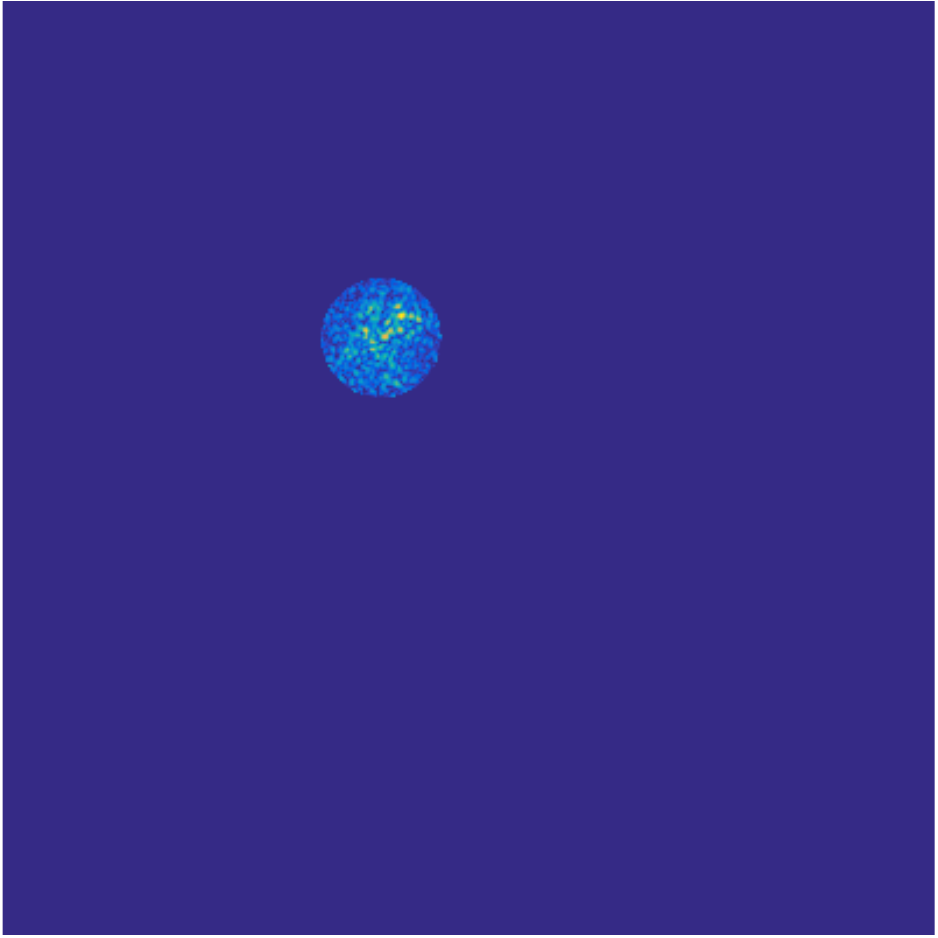} & \includegraphics[scale=0.37]{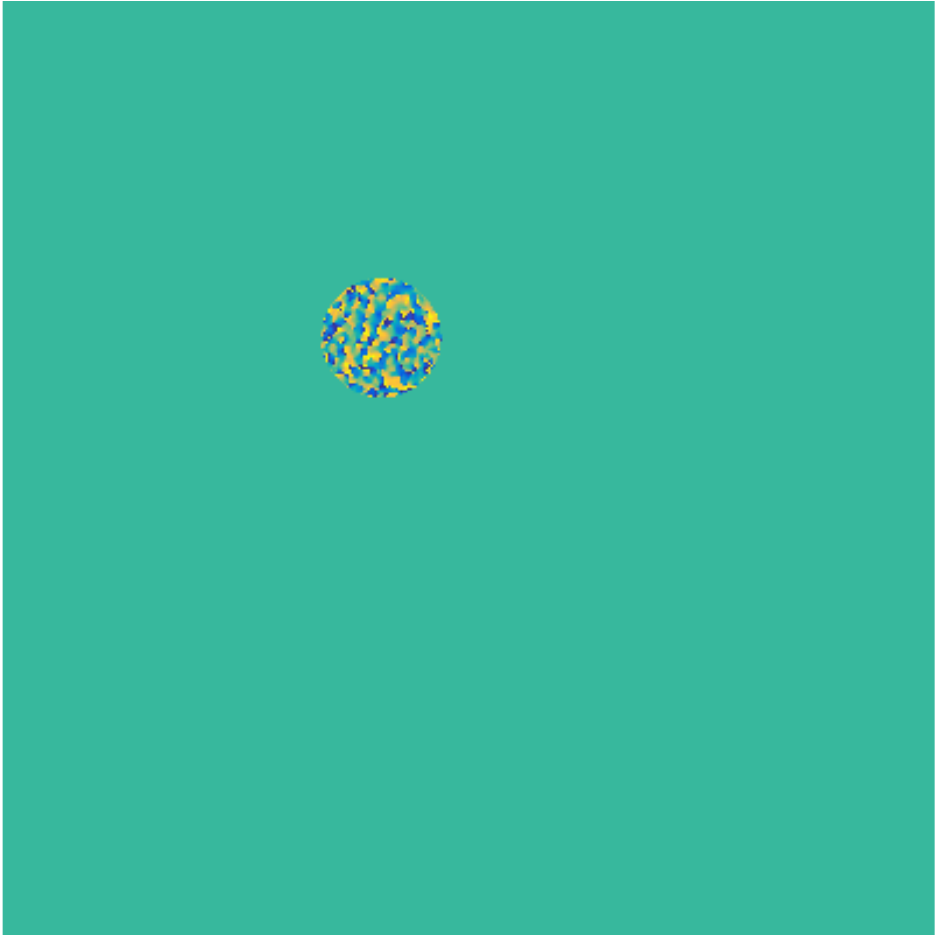}
	\end{tabular}
	\caption{Amplitude and phase of the horizontal  polarisation of one calibration input ($\mathbf{A}_m$ of Eq.~\eqref{eq:rep_sys_vf_A_B_C}) at the distal end and at the proximal end. 
	}
	\label{fig:CRUKcalibration}
\end{figure}

The fibre is calibrated using input and output pairs such as those shown in Figure \ref{fig:CRUKcalibration}.  For the purposes of reconstruction, each calibration input and output is separated from the others by evaluating each of them only over a circular region around the centre of the corresponding Gaussian-like spot. In particular, the calibration inputs in Figure \ref{fig:CRUKcalibration} are evaluated on a grid with resolution $1200\times1200$ and they are translated to $M=936$ different locations across the input imaging plane. Each output is evaluated at $N=34973$ different pixels at the output imaging plane. Considering the two polarisation states, the dimension of the system matrix in \R{eq:pol_system2} is therefore $1872\times69946$.

\begin{figure}[htb!]
	\centering
	\begin{tabular}{ccccc}
	& \multicolumn{2}{c}{Horizontal polarisation} & \multicolumn{2}{c}{Vertical polarisation}   \\
	& Amplitude & Phase & Amplitude & Phase \\
	\rotatebox{90}{\hspace{0.35cm} naive approach} &\includegraphics[scale=0.55]{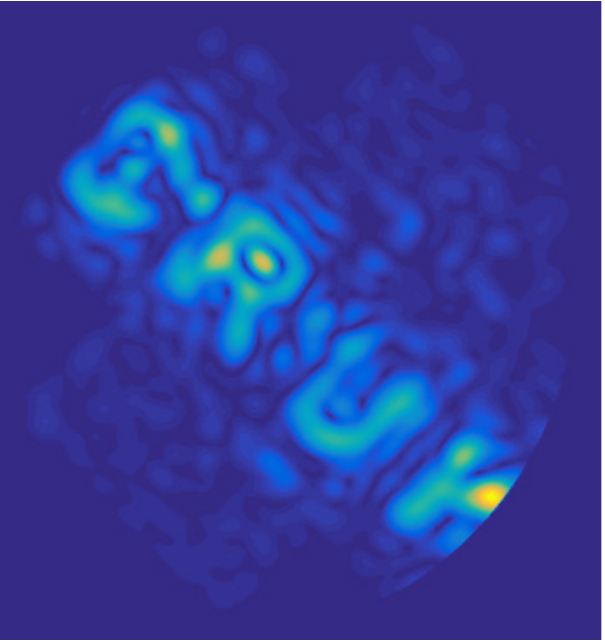} & \includegraphics[scale=0.55]{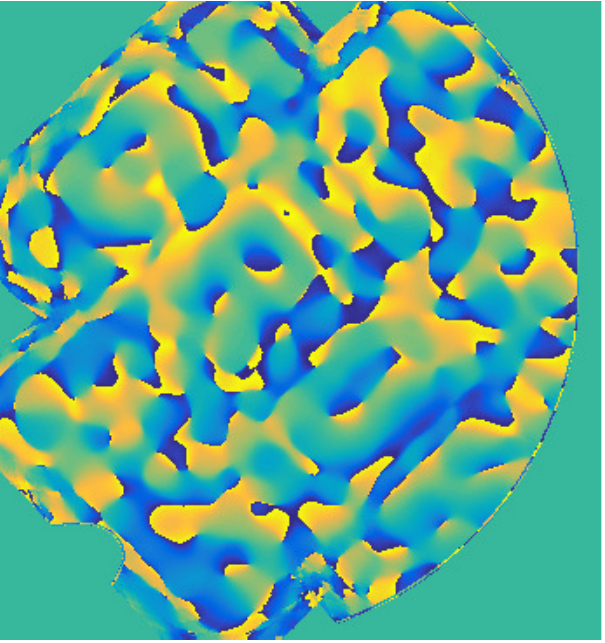} & \includegraphics[scale=0.55]{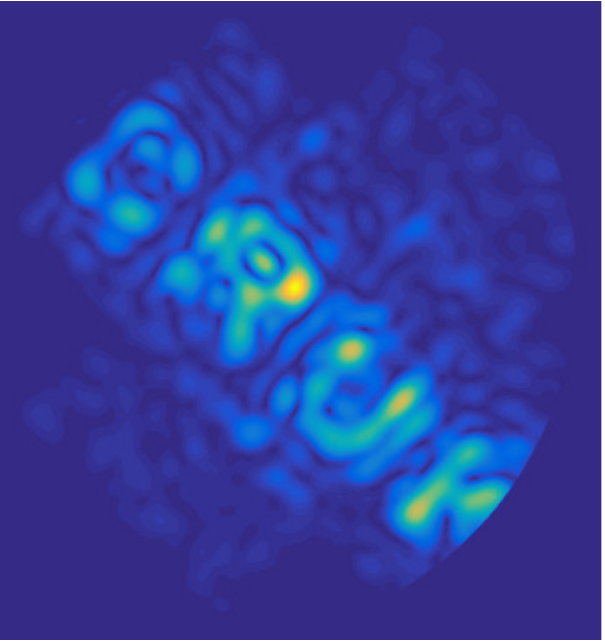} & \includegraphics[scale=0.55]{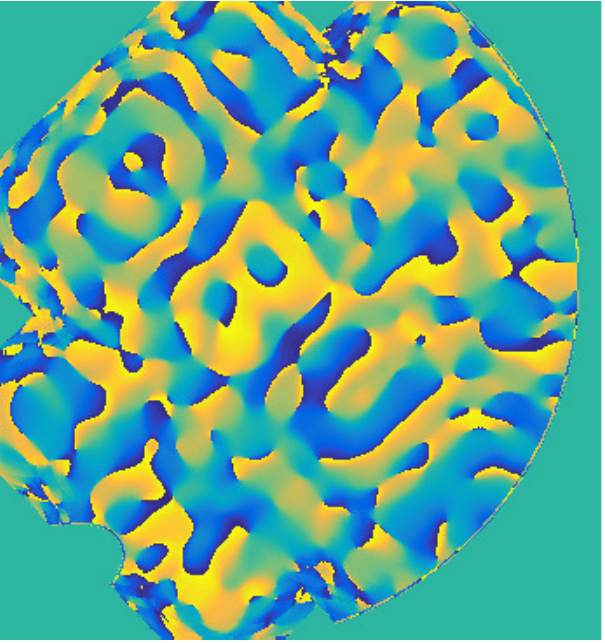} \\
	\rotatebox{90}{\hspace{0.55cm} least-squares} &\includegraphics[scale=0.55]{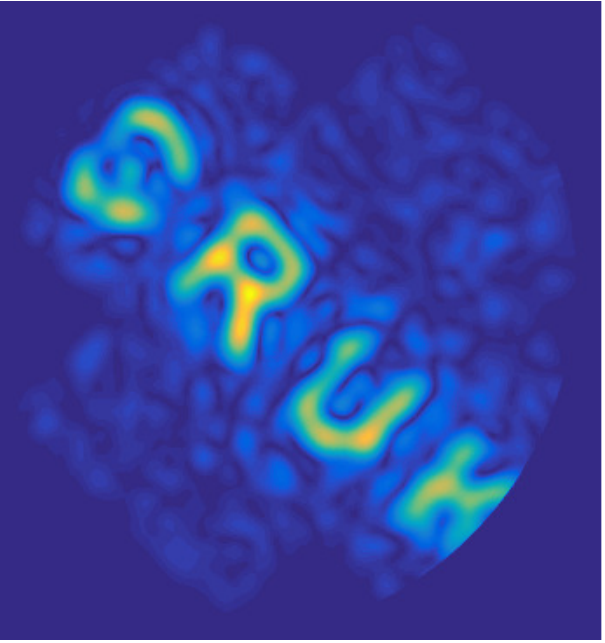} & \includegraphics[scale=0.55]{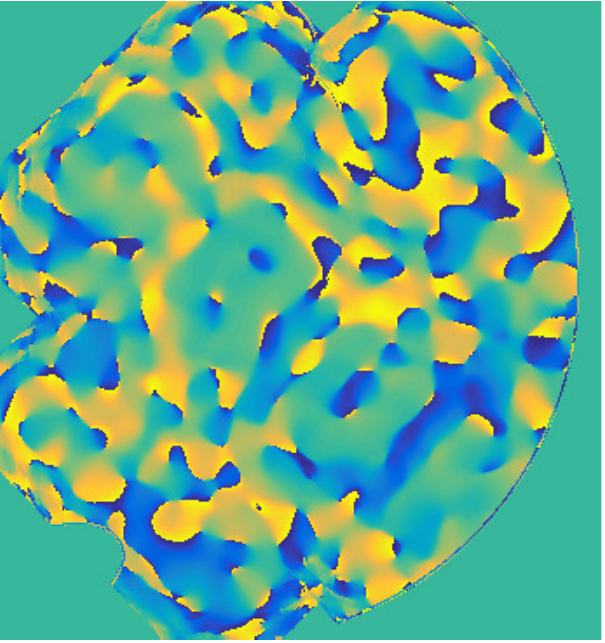} & \includegraphics[scale=0.55]{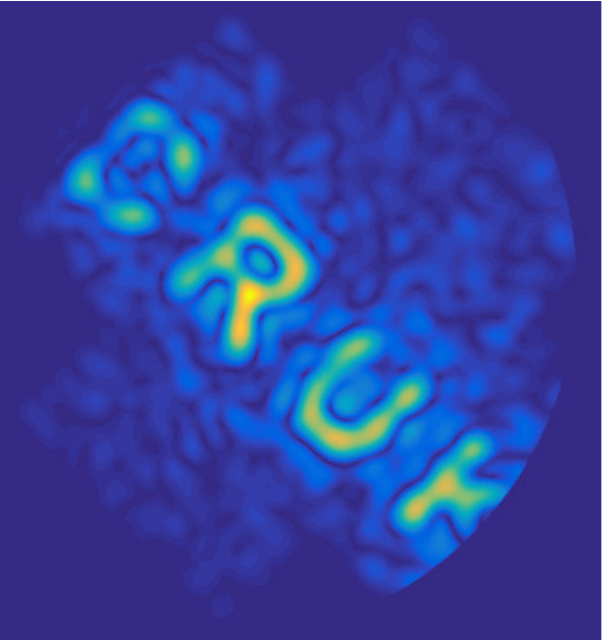} & \includegraphics[scale=0.55]{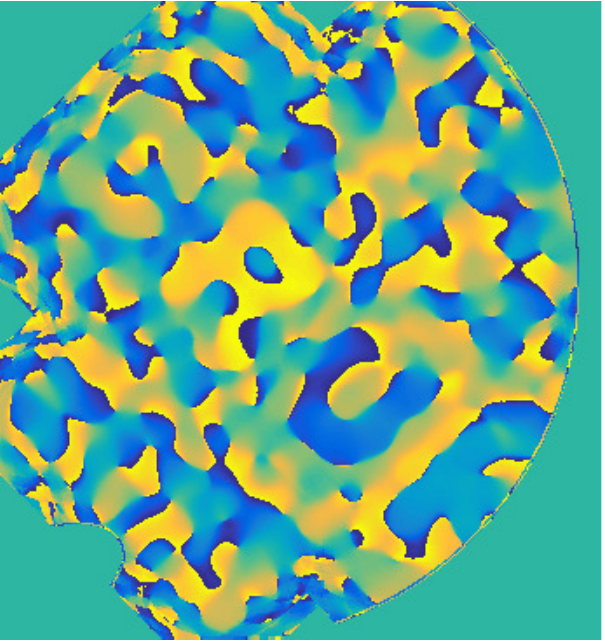} \\
	\rotatebox{90}{\hspace{0.35cm} $\ell_2$-regularisation} &\includegraphics[scale=0.55]{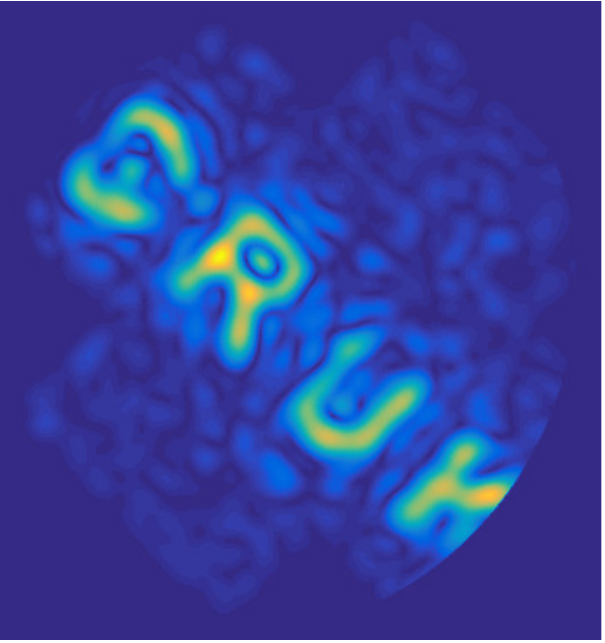} & \includegraphics[scale=0.55]{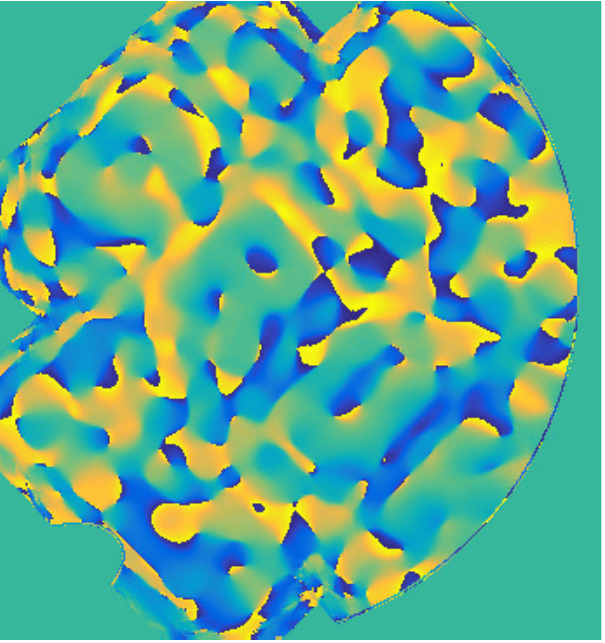} & \includegraphics[scale=0.55]{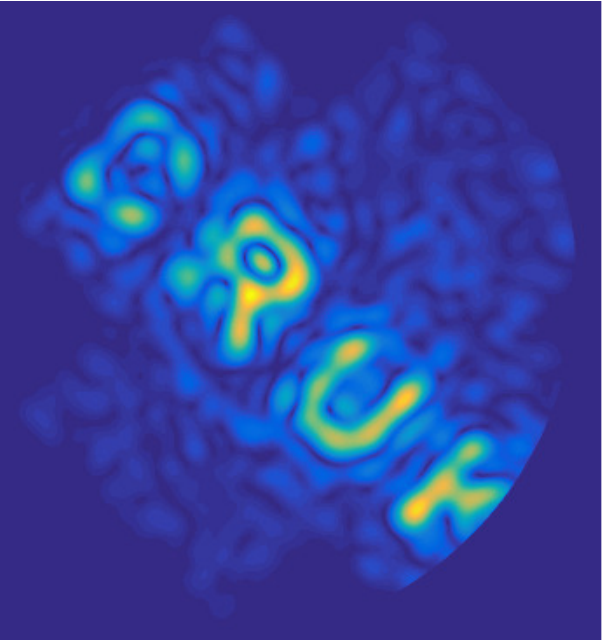} & \includegraphics[scale=0.55]{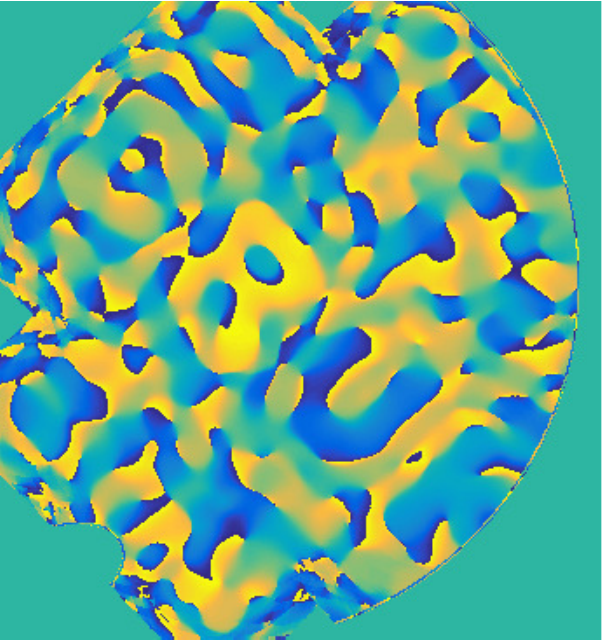} \\
	\rotatebox{90}{\hspace{0.35cm} $\ell_1$-regularisation} &\includegraphics[scale=0.55]{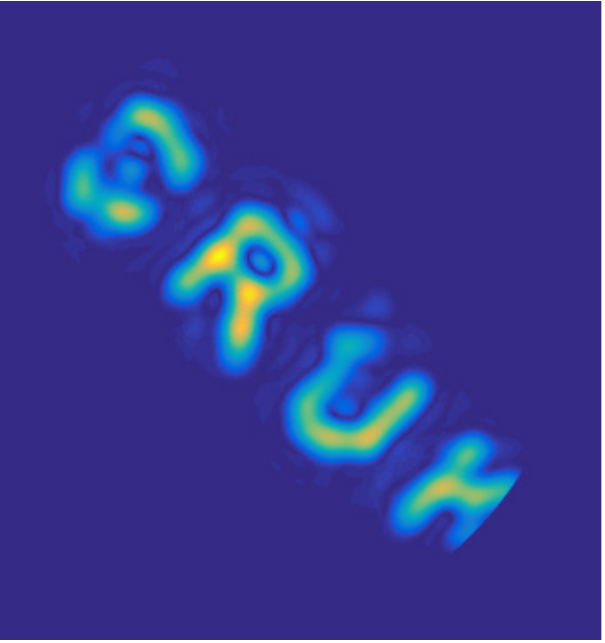} & \includegraphics[scale=0.55]{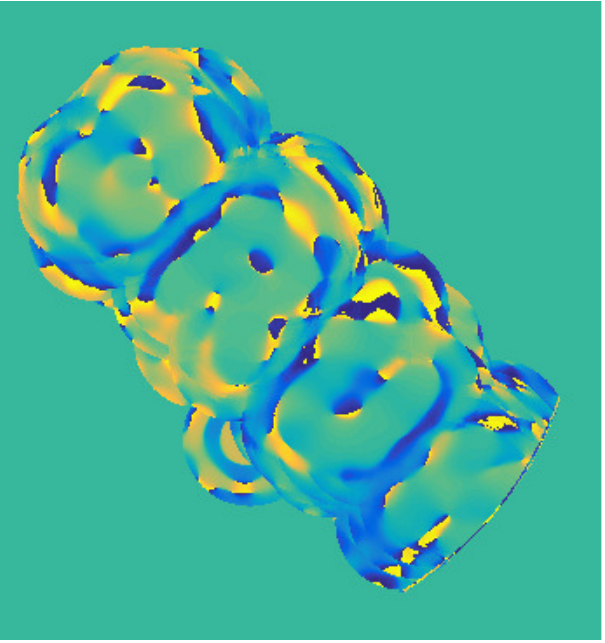} & \includegraphics[scale=0.55]{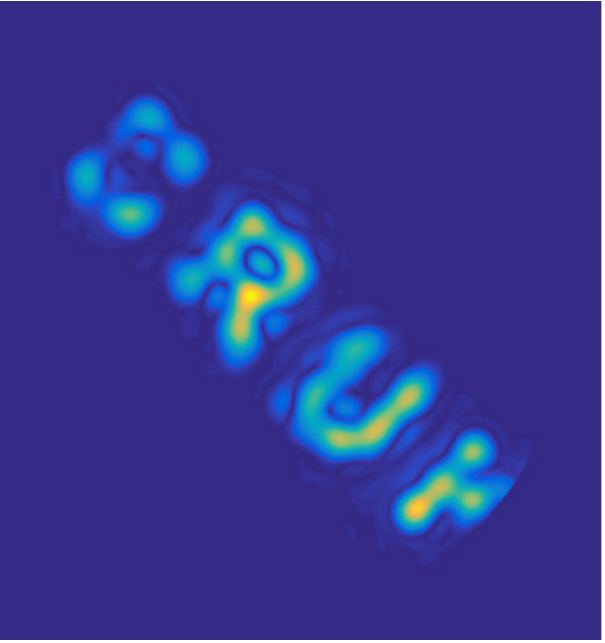} & \includegraphics[scale=0.55]{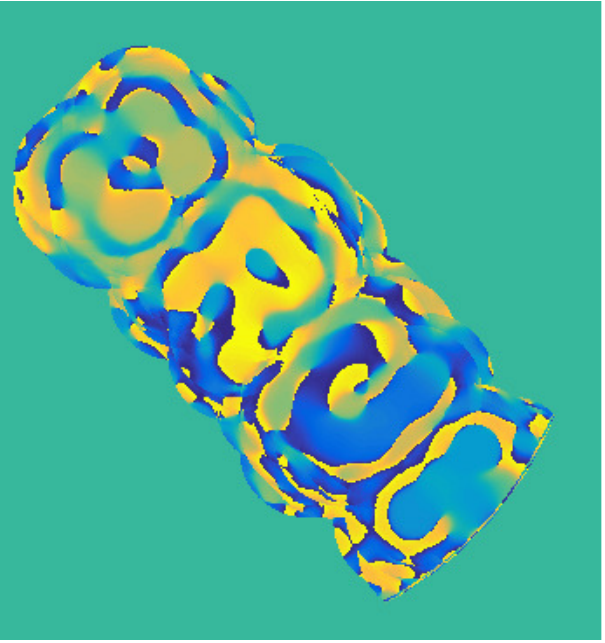} \\
	\end{tabular}
	\caption{Reconstructed amplitude and phase of the horizontal and vertical polarisations of the holographic image from Figure~\ref{fig:CRUKdata} with respect to the calibration functions such as those in Figure~\ref{fig:CRUKcalibration}, using naive, least-squares, $\ell_2$ and $\ell_1$ approaches. The regularisation parameter in the $\ell_2$ and $\ell_1$-regularisation is $\lambda=0.3$ and  $\lambda=0.267$, respectively.}
	\label{fig:CRUKrec}
\end{figure}

In Figure~\ref{fig:CRUKrec}, we recover the amplitude and the phase of the horizontal and vertical polarisations of the holographic image from raw endoscopic measurements using different inversion techniques while reconstructing with respect to the calibration coefficients. In particular, we solve \eqref{eq:pol_system2} where $\mathbf{H}^h=\mathbf{H}^v=\mathbf{I}$, by inverting the linear system in four different ways:
\begin{enumerate}[label=\arabic*.]
\item
the naive inversion $\bar{\mathbf{f}} := \mathbf{E}^* \mathbf{g} $, which corresponds to the principle of phase conjugation in that it assumes $\mathbf{E}^*\mathbf{E}=\mathbf{I}$,
\item
the least-squares approach $\bar{\mathbf{f}} :=(\mathbf{E}^* \mathbf{E} )^{-1} \mathbf{E}^* \mathbf{g}$,
\item
the $\ell_2$-regularisation $\bar{\mathbf{f}} := (\mathbf{E}^* \mathbf{E} + \lambda \mathbf{I} )^{-1} \mathbf{E}^* \mathbf{g}$, and,
\item
the $\ell_1$-regularisation $\bar{\mathbf{f}} := \argmin_{\mathbf{f}\in\mathbb{C}^{2M} } \| \mathbf{g} - \mathbf{E} \mathbf{f} \|_2+\lambda\|\mathbf{f}\|_1$ using the iterative solver \cite{spgl12007}.
\end{enumerate}
In particular, we see in Figure~\ref{fig:CRUKrec} that $\ell_1$-regularisation performs well when compared with the other approaches. In fact, since our holographic image is sparse with respect to the calibration inputs -- namely, since $F^h$ and $F^v$ are sparse with respect to $\{E_{m}^{h}\}_{m=1}^{M/2}$ and $\{ E_{m}^{v}\}_{m=1}^{M/2}$ -- $\ell_1$-regularisation successfully removes significant noise while preserving the image details. The amount of noise that is removed by $\ell_1$-regularisation depends on the strength of the regularisation parameter $\lambda$.

\begin{figure}[htb!]
	\centering
	\begin{tabular}{ccccc}
	& \multicolumn{2}{c}{Fourier exponentials} & \multicolumn{2}{c}{DB4 wavelets}   \\
	& Amplitude & Phase & Amplitude & Phase \\
	\rotatebox{90}{\hspace{0.2cm} naive approach} &\includegraphics[scale=0.55]{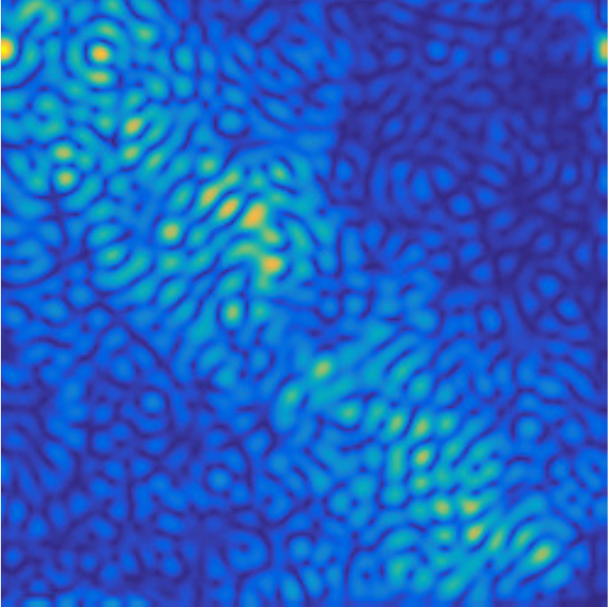} & \includegraphics[scale=0.55]{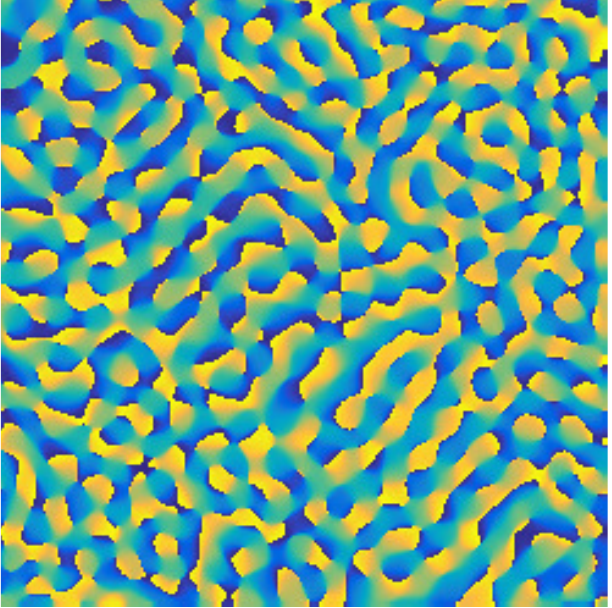} & \includegraphics[scale=0.55]{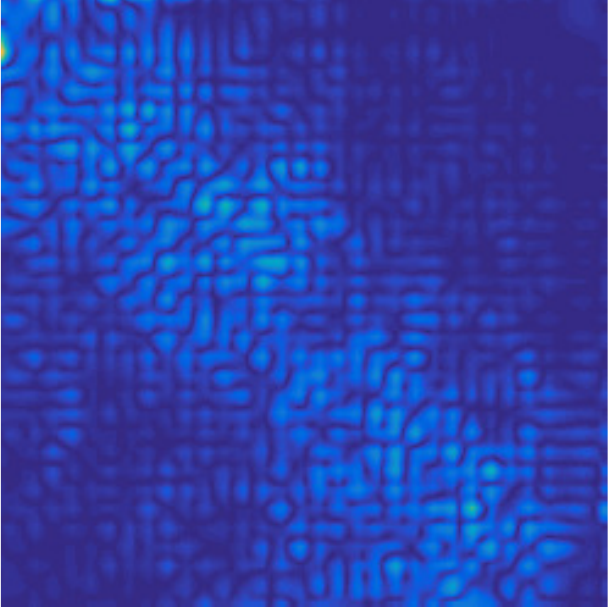} & \includegraphics[scale=0.55]{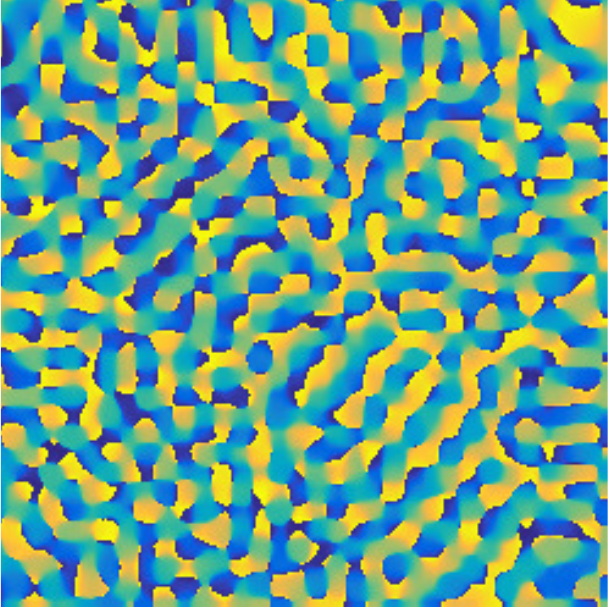}  \\
	\rotatebox{90}{\hspace{0.5cm} least-squares} &\includegraphics[scale=0.55]{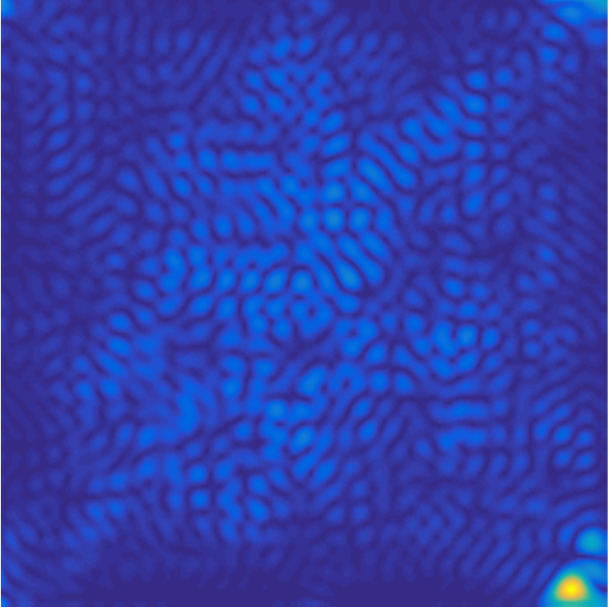} & \includegraphics[scale=0.55]{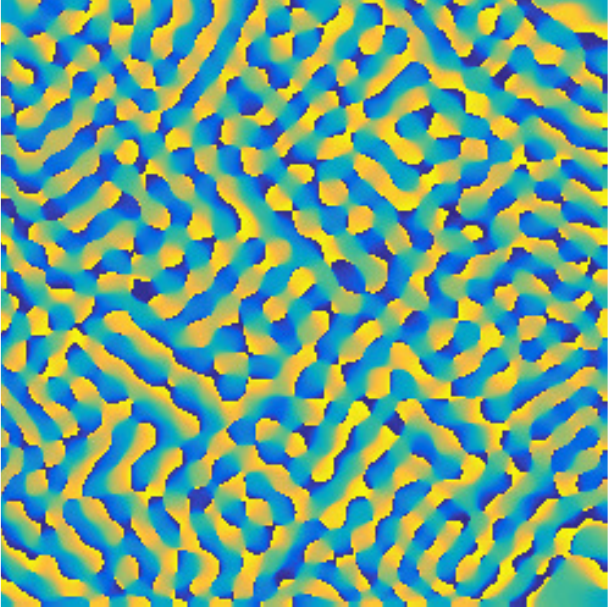} & \includegraphics[scale=0.55]{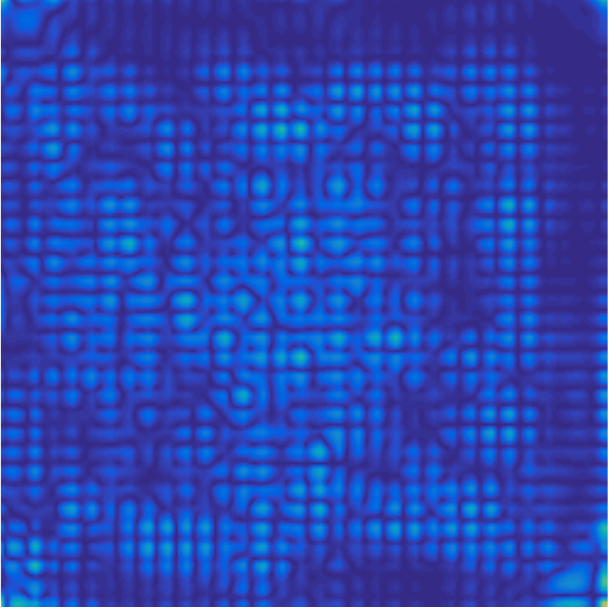} & \includegraphics[scale=0.55]{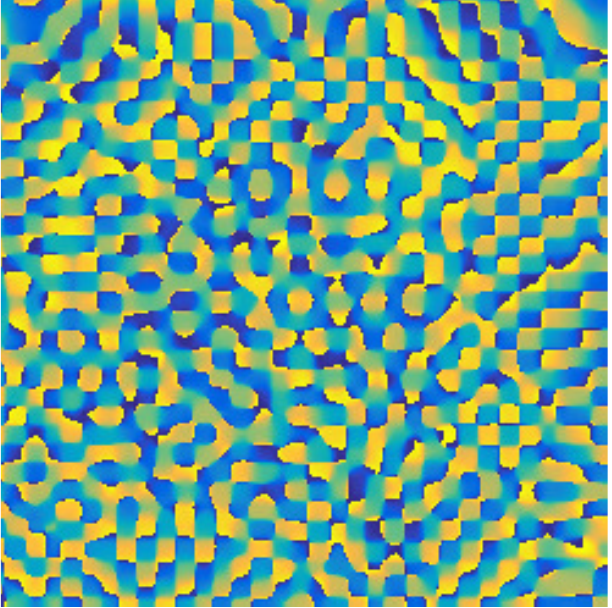}  \\
	\rotatebox{90}{\hspace{0.1cm} $\ell_2$-regularisation} &\includegraphics[scale=0.55]{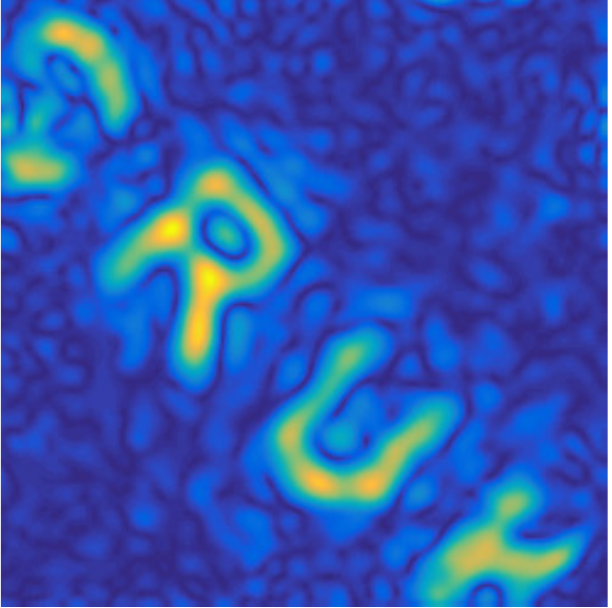} & \includegraphics[scale=0.55]{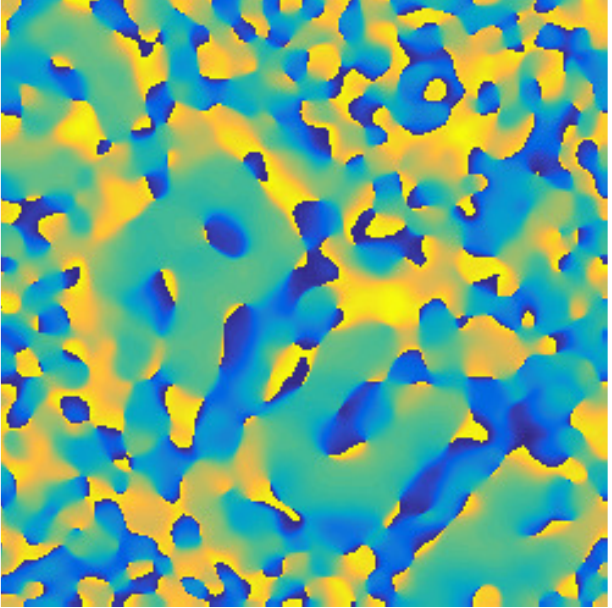} & \includegraphics[scale=0.55]{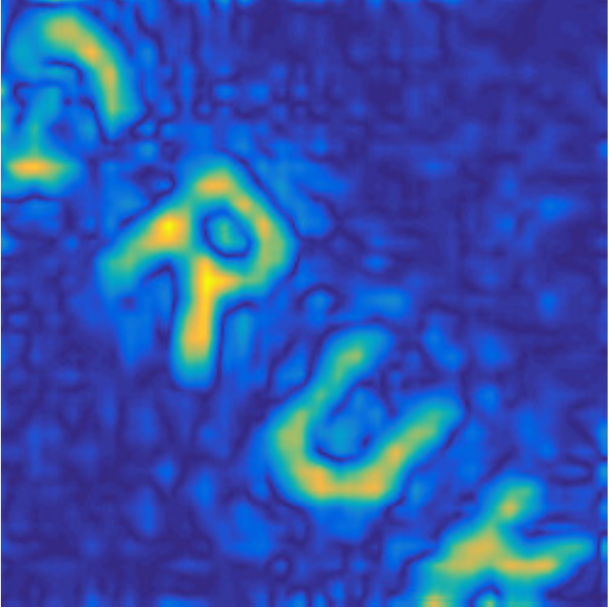} & \includegraphics[scale=0.55]{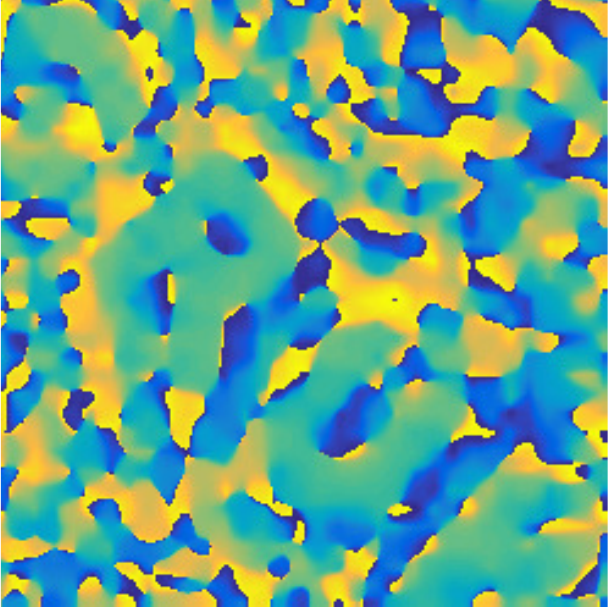}  \\ 
	\rotatebox{90}{\hspace{0.1cm} $\ell_1$-regularisation} &\includegraphics[scale=0.55]{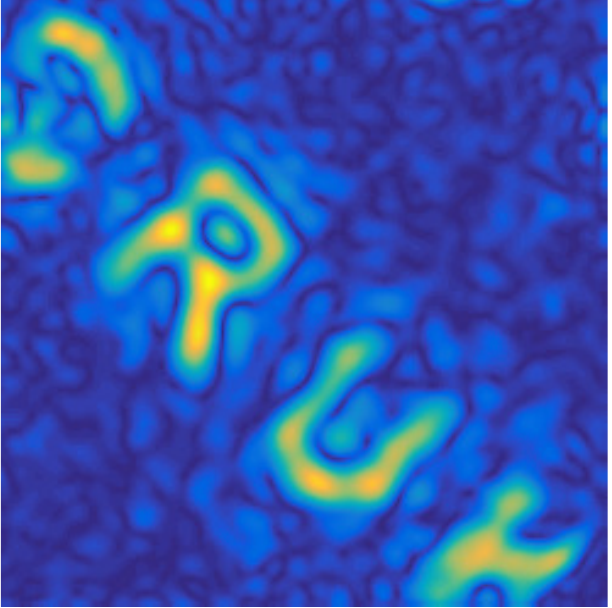} & \includegraphics[scale=0.55]{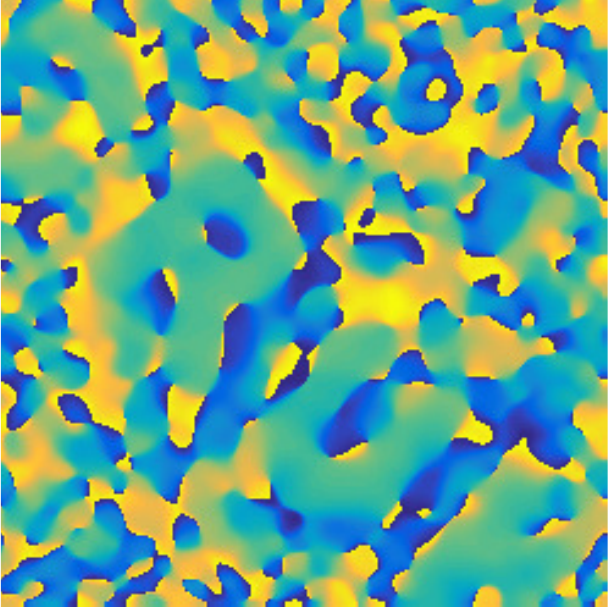} & \includegraphics[scale=0.55]{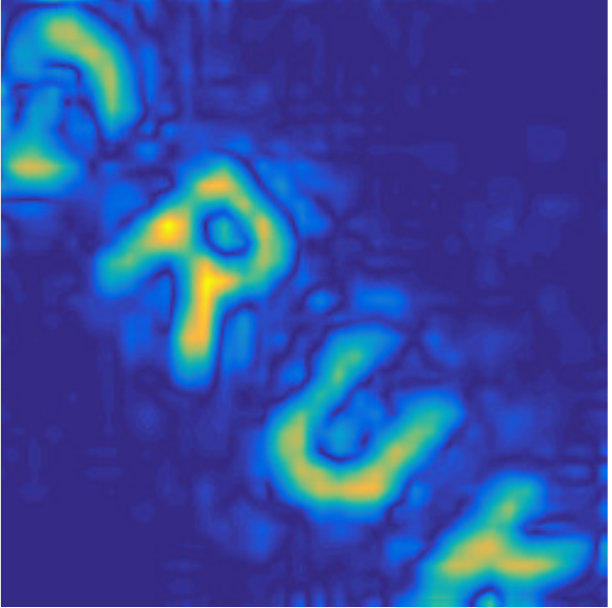} & \includegraphics[scale=0.55]{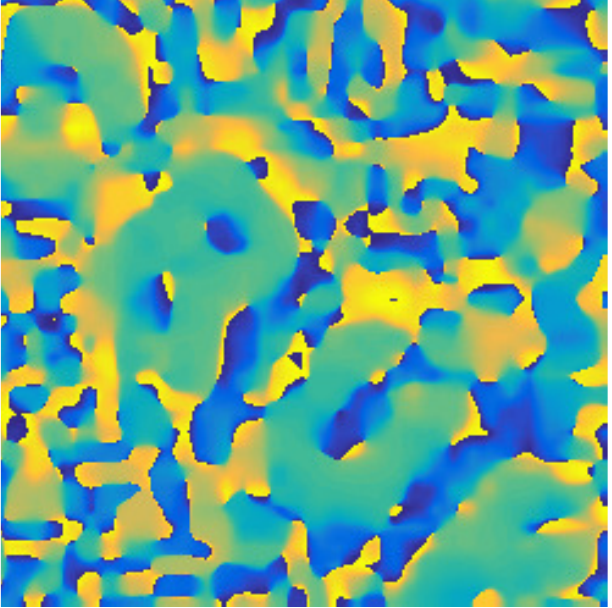} \\ 
	\end{tabular}
	\caption{Reconstructed images (horizontal polarisation) with respect to the different bases using different inversion approaches. We used $32\times32$ Fourier exponentials and $32\times32$ DB4 wavelets. Regularisation parameter in the $\ell_2$ and $\ell_1$-regularisation is $\lambda=10$ and  $\lambda=0.25$, respectively.}
	\label{fig:CRUKbases}
\end{figure}

Finally, in Figure~\ref{fig:CRUKbases}, we reconstruct the holographic image with respect to different representation systems, namely we solve \eqref{eq:pol_system2}  where both $\mathbf{H}^h$ and $\mathbf{H}^v$ correspond to the representation system of a Fourier or a wavelet basis with cardinality $K=1024$. Specifically,
\begin{enumerate}[label=(\roman*)]
\item
in the Fourier case, we choose both $\{H_k^h\}_{k=1}^{K}$ and $\{H_k^v\}_{k=1}^{K}$ to be $\{\exp(2\pi\I (k_1 x_1 + k_2 x_2)): k_1,k_2=-\sqrt{K}/2,\ldots,\sqrt{K}/2-1\}$, where $\mathbf{x}=(x_1,x_2)\in[-1/2,1/2]^2$, whereas
\item
in the wavelet case, we choose tensor-products of $\sqrt{K}$ one-dimensional boundary-corrected Daubechies wavelets with four vanishing moments (DB4) from \cite{Cohen1993}.
\end{enumerate}
We can observe in Figure~\ref{fig:CRUKbases} that least-squares fails to give a useful estimate due to the ill-conditioning of the system matrix, conveying that it is crucial to use the regularisation term.  Although the least-squares could still be used in the case where $K\ll M$, small $K$ does not necessarily lead to a good approximation of the image, and so to achieve the desired resolution one would need to increase the number of calibration measurements $M$. This is undesirable as it would incur additional experimental time. On the other hand, given that our holographic image is sparse with respect to compactly-supported wavelets, we see from  Figure~\ref{fig:CRUKbases} that $\ell_1$-regularisation performs quite well in combination with DB4 wavelets even though $K>M$. 

\subsection{Reconstruction and analysis of biological images}\label{ss:biological}

We now apply the framework developed in Sections~\ref{s:framework}--\ref{s:fourierfeatures} to reconstruct and analyse images of real-world biological samples. We imaged ex vivo samples of mouse oesophagus from healthy controls and from carcinogen treated animals with induced oesophageal tumours (lesions) using the model presented in \cite{Alcolea2014}. More concretely, 3 control mice (6 healthy areas analysed) and 6 mice with induced tumours (6 distinct lesions analysed) were used.  Each sample was segmented into areas of healthy and lesion tissue by an expert using DAPI fluorescence images.

For clarity in the examples below, we index different areas by $n=1,\ldots,12$, where the first six are healthy and the rest are lesions. Due to the limited field of view of the endoscope ($\sim 200 \mu \mathrm{m}$) relative to the sample size ($\sim 2$mm), each of the 9 mice produce 6--20 individual images corresponding to different areas of the same sample that may overlay by up to 15\%.  We therefore also introduce index $i$ to denote individual sub-images within a larger area on a given sample. Each individual sample in the data set thus has associated index $(n,i)$, where $n=1,\ldots,12$ and $i=1,\ldots, I_n$, for $I_n$ in the range 6--20.



\begin{figure}[h]
	\centering
	{\small Space domain} \hspace{4.5cm} {\small Fourier domain} \\
	\hspace{0.7cm} {\small Amplitude} \hspace{1.cm}{\small Unwrapped phase} \hspace{1.1cm} {\small Amplitude} \hspace{1.2cm} {\small Fitted Gaussian}\\
	\rotatebox{90}{\hspace{0.8cm} \footnotesize{HEALTHY}}
	\adjincludegraphics[scale=0.6]{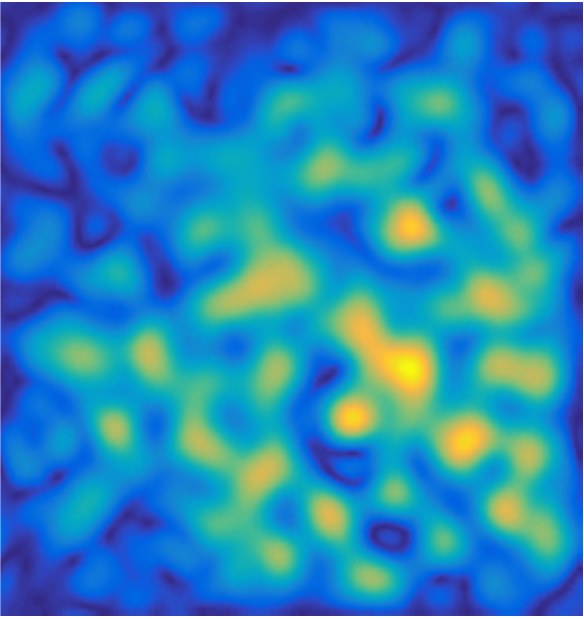}\adjincludegraphics[scale=0.6]{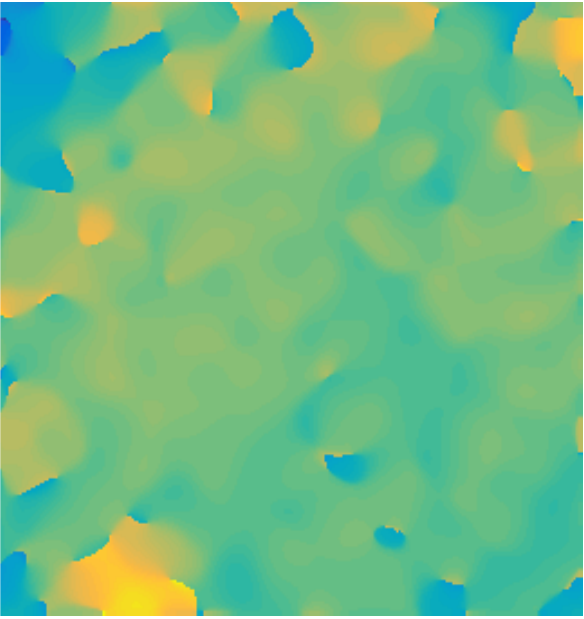}\adjincludegraphics[scale=0.6]{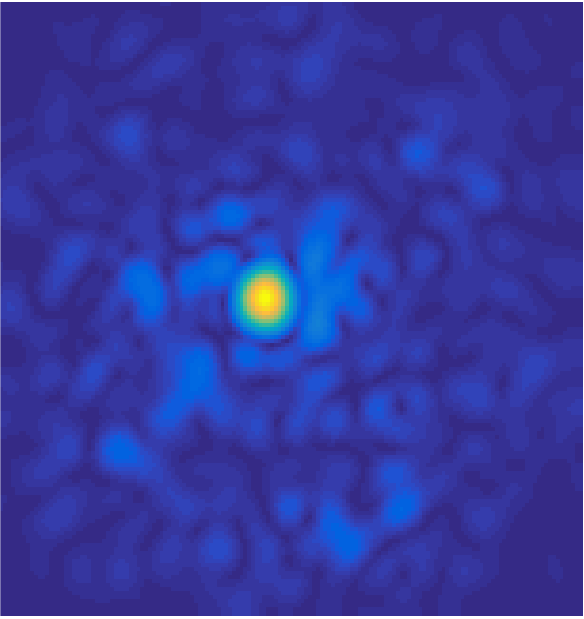}\adjincludegraphics[scale=0.6]{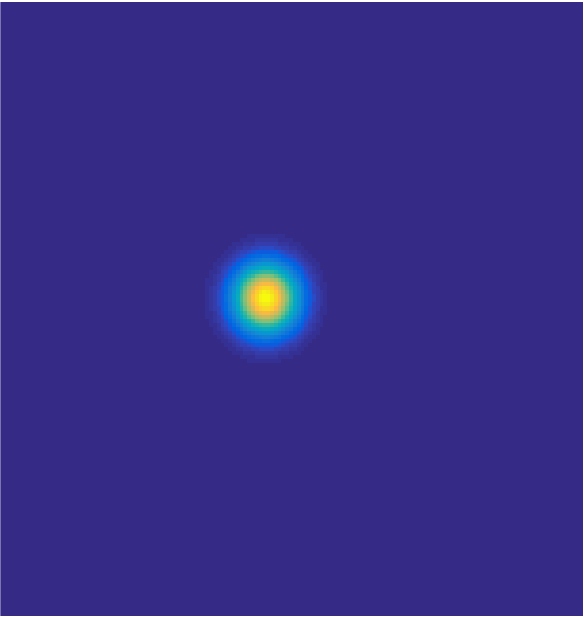}\\
	\rotatebox{90}{\hspace{1.1cm} \footnotesize{LESION}}
	\adjincludegraphics[scale=0.6]{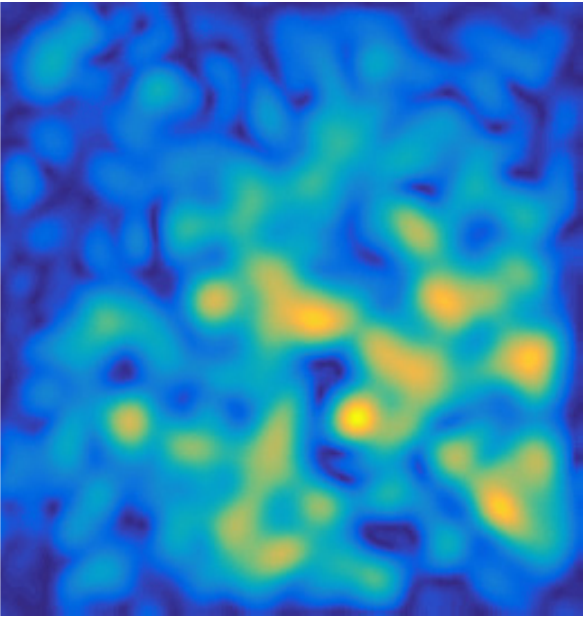}\adjincludegraphics[scale=0.6]{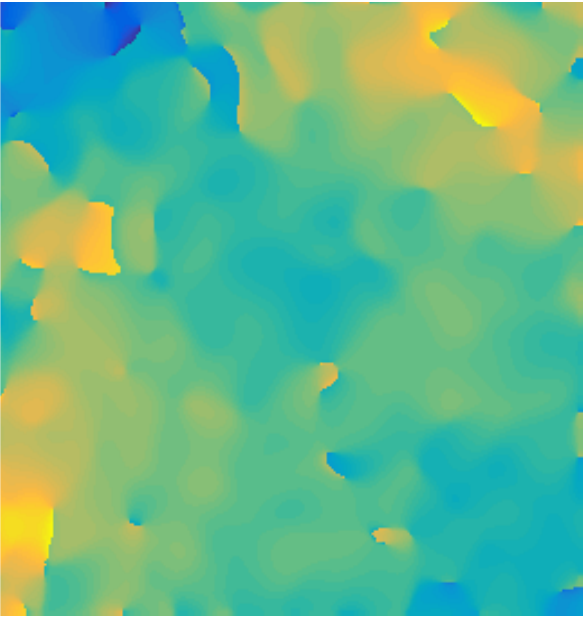}\adjincludegraphics[scale=0.6]{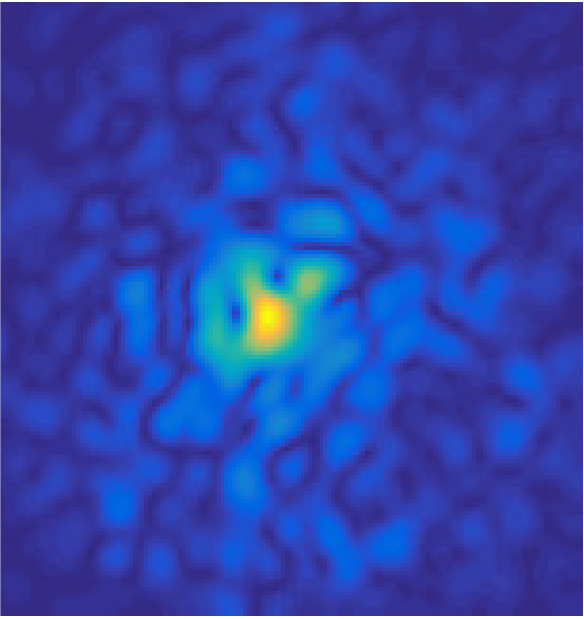}\adjincludegraphics[scale=0.6]{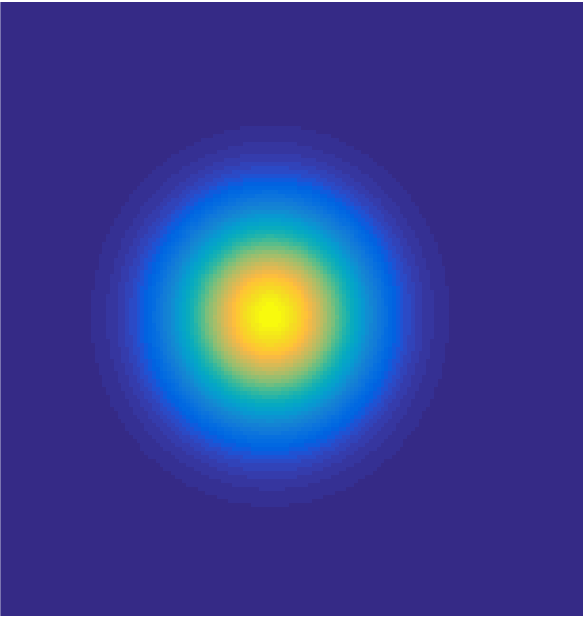}
	\caption{Healthy and lesion tissue reconstructed with respect to $K=400$ Fourier coefficients by solving the linear system \R{eq:pol_system} with $\ell_1$-regularisation  with $\lambda=0.25$. Only the horizontal polarisation in shown due to space limitations.}
	\label{fig:Real_HvL_fourier_example}
\end{figure}

In Figure~\ref{fig:Real_HvL_fourier_example}, we show the reconstruction of
the horizontal polarisation of one healthy image indexed as $(1,1)$ (first row) and one lesion image indexed as $(7,1)$ (second row), in both the space-domain (first two columns) and the Fourier-domain (last two columns). Specifically, we reconstruct $K=400$ Fourier coefficients per polarisation by solving the linear system \R{eq:pol_system} with $\ell_1$-regularisation. We then expand these coefficients with respect to Fourier-exponentials to obtain images in the space-domain, as well as, with respect to sinc-functions to obtain images in the Fourier-domain. In the space-domain, we show the amplitude (first column) and unwrapped phase (second column) of the reconstructed image, where for the unwrapping we used an efficient 2D phase unwrapper\footnote{Available on-line at \href{https://www.ljmu.ac.uk/research/centres-and-institutes/faculty-of-engineering-and-technology-research-institute/geri/phase-unwrapping}{www.ljmu.ac.uk/research/centres-and-institutes/faculty-of-engineering-and-technology-research-institute/geri/phase-unwrapping}}. In the Fourier domain, we show the amplitude (third column) of the reconstructed Fourier transform as well as a Gaussian fit to the amplitude of the Fourier transform (fourth column), where we used the procedure explained in Section~\ref{s:fourierfeatures}.
While the difference between the healthy and the lesion sample is not so apparent from the amplitude and phase in the space-domain, the difference becomes more pronounced in the Fourier domain; specifically, it can be observed that the Fourier coefficients decay slower in the lesion than in the healthy tissue, where the decay is quantified by the standard deviation of the fitted Gaussian.

\begin{figure}[h!]
(a)\\
\rotatebox{90}{\hspace{0.3cm} \footnotesize{HEALTHY}}
	\adjincludegraphics[scale=0.61]{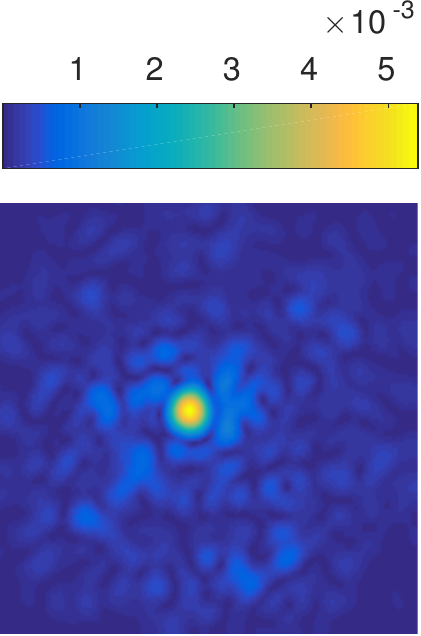}\adjincludegraphics[scale=0.61]{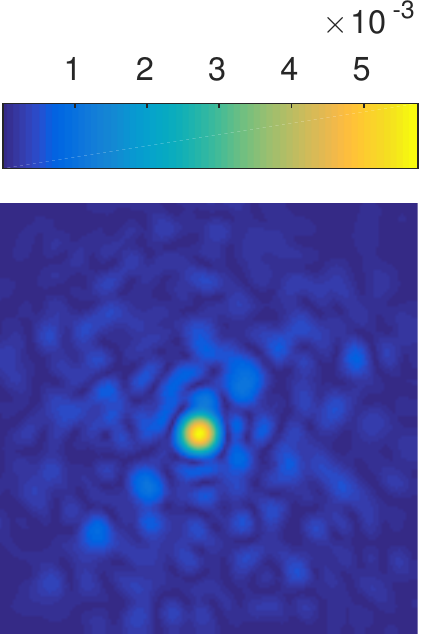}\adjincludegraphics[scale=0.61]{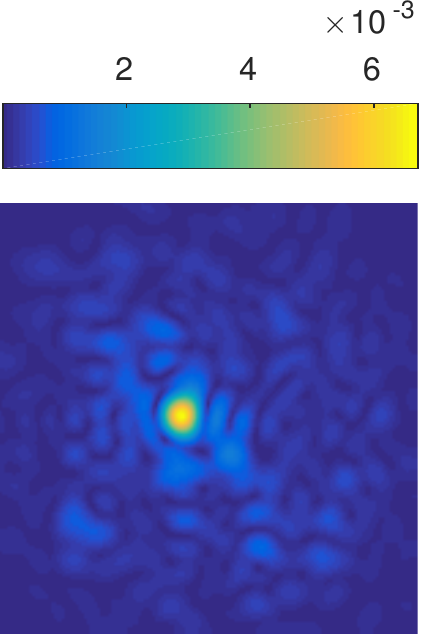}\adjincludegraphics[scale=0.61]{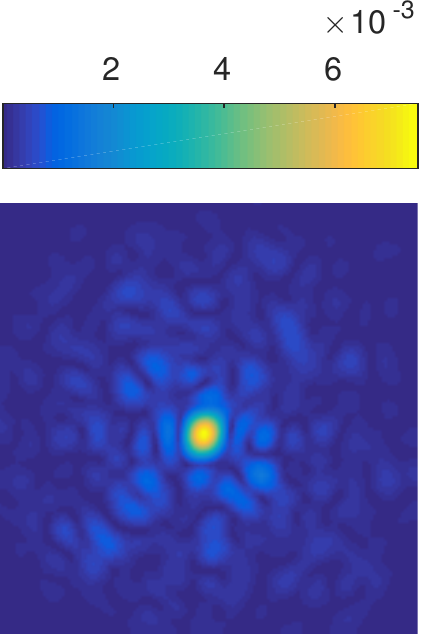}\adjincludegraphics[scale=0.61]{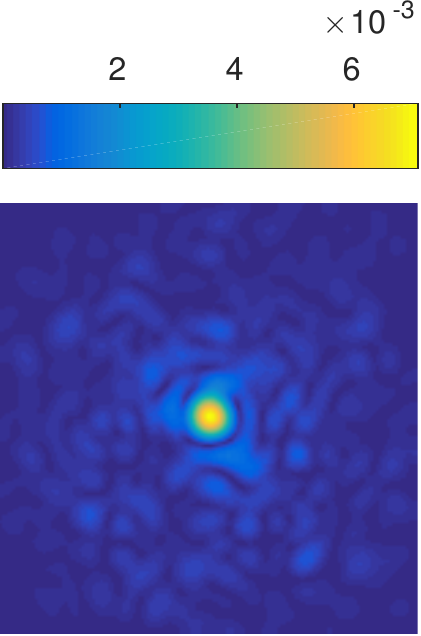}\adjincludegraphics[scale=0.61]{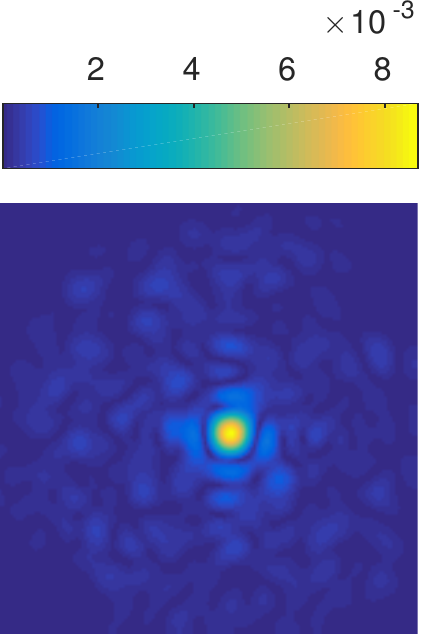}\\
\rotatebox{90}{\hspace{2.1cm} \footnotesize{LESION}}
\adjincludegraphics[scale=0.61]{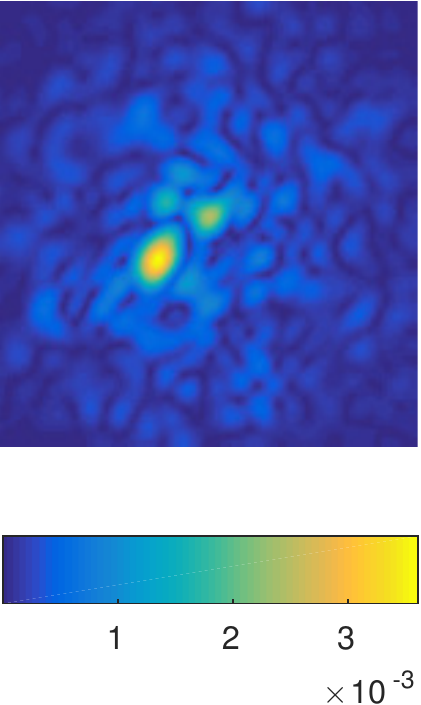}\adjincludegraphics[scale=0.61]{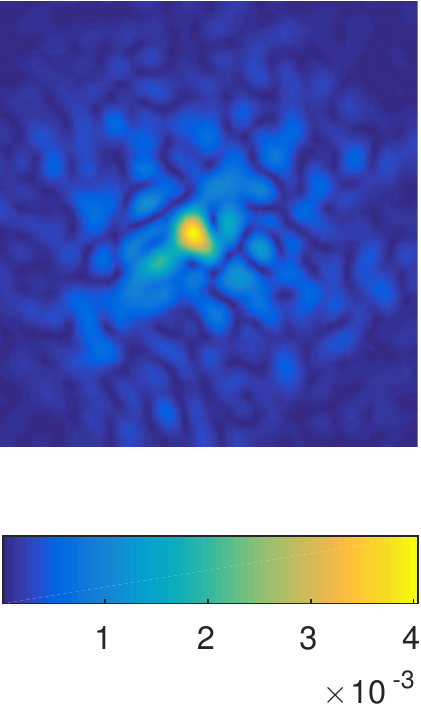}\adjincludegraphics[scale=0.61]{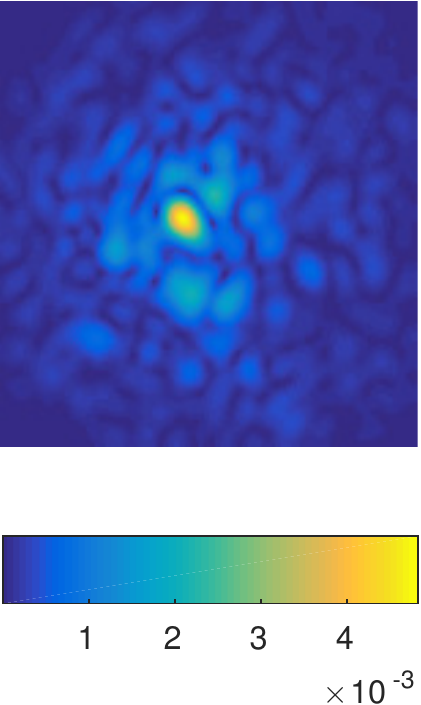}\adjincludegraphics[scale=0.61]{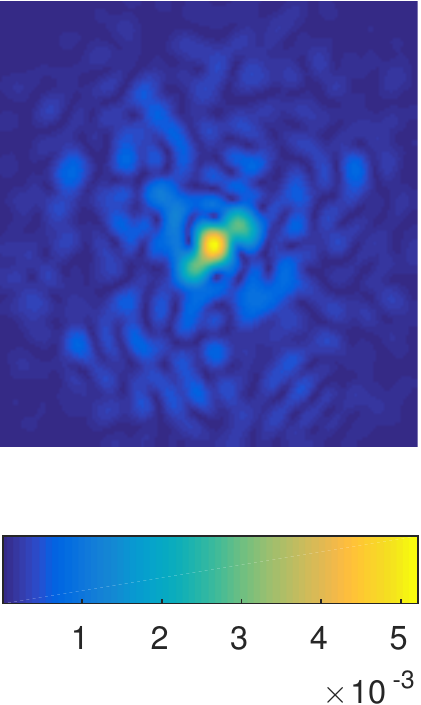}\adjincludegraphics[scale=0.61]{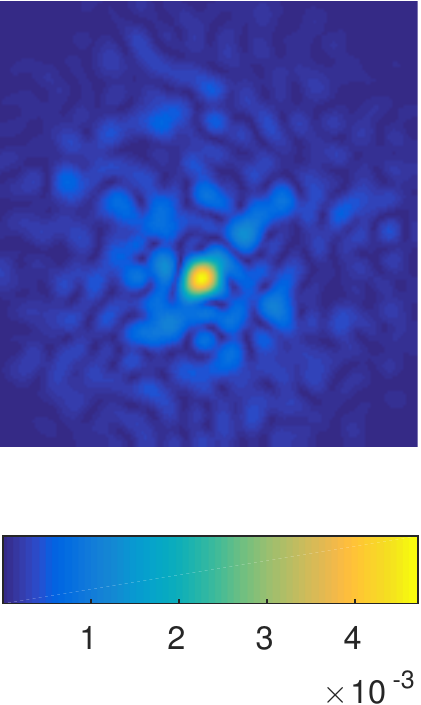}\adjincludegraphics[scale=0.61]{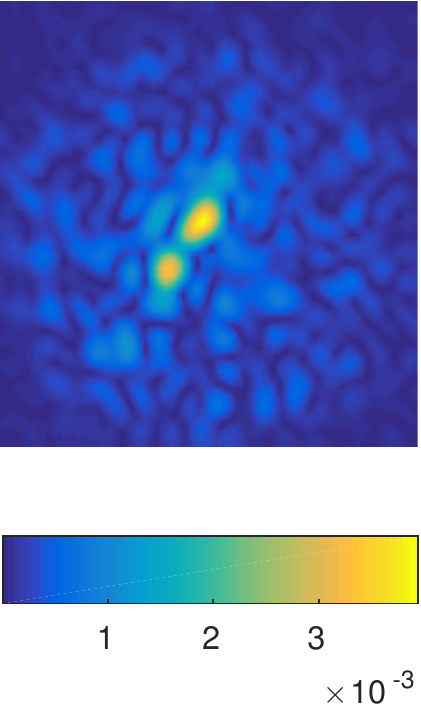}\\
(b)\\
\rotatebox{90}{\hspace{0.3cm} \footnotesize{HEALTHY}}
	\adjincludegraphics[scale=0.61]{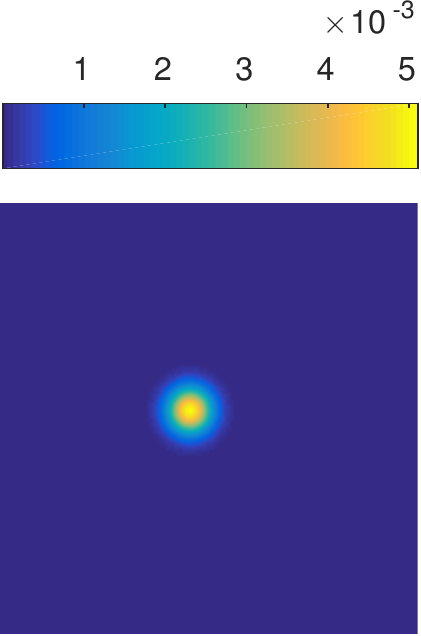}\adjincludegraphics[scale=0.61]{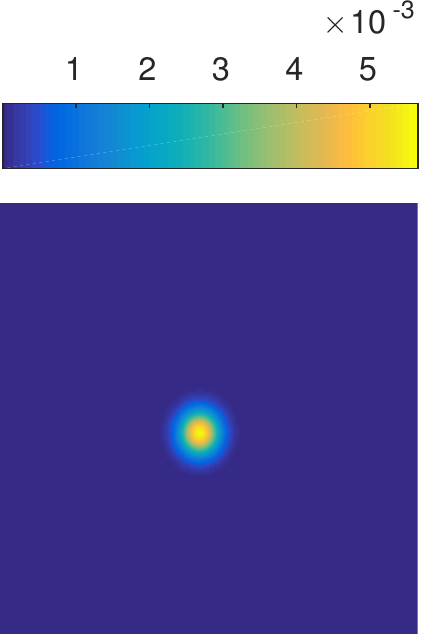}\adjincludegraphics[scale=0.61]{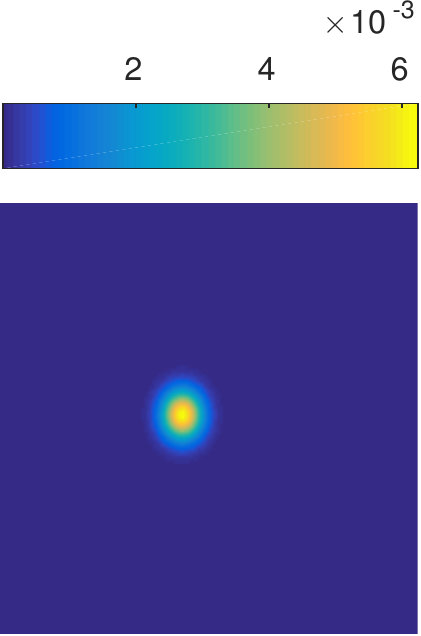}\adjincludegraphics[scale=0.61]{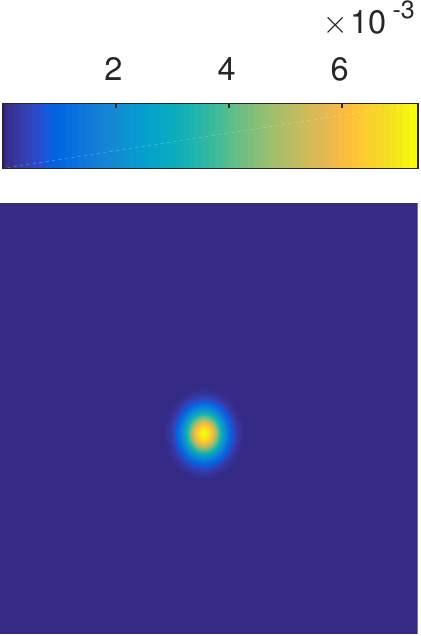}\adjincludegraphics[scale=0.61]{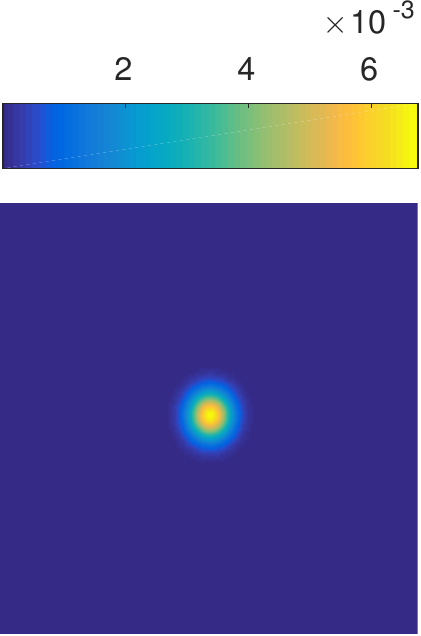}\adjincludegraphics[scale=0.61]{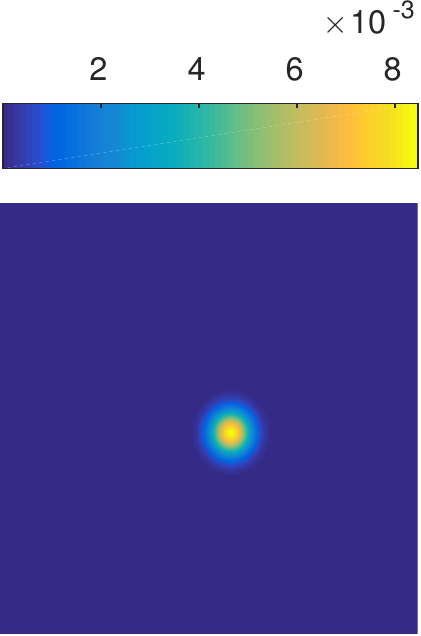}\\
	\rotatebox{90}{\hspace{2.1cm} \footnotesize{LESION}}
	\adjincludegraphics[scale=0.61]{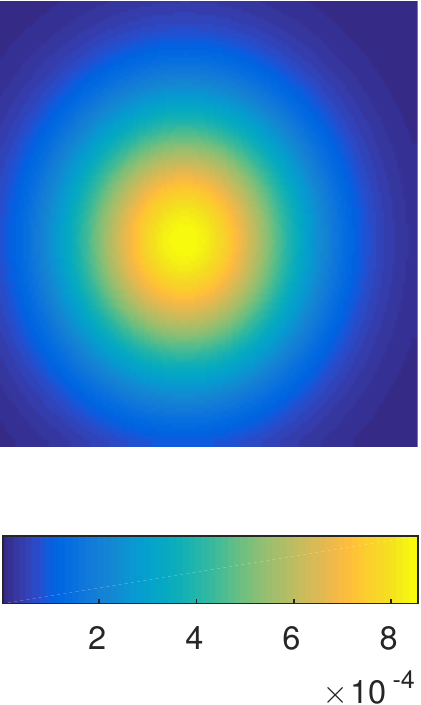}\adjincludegraphics[scale=0.61]{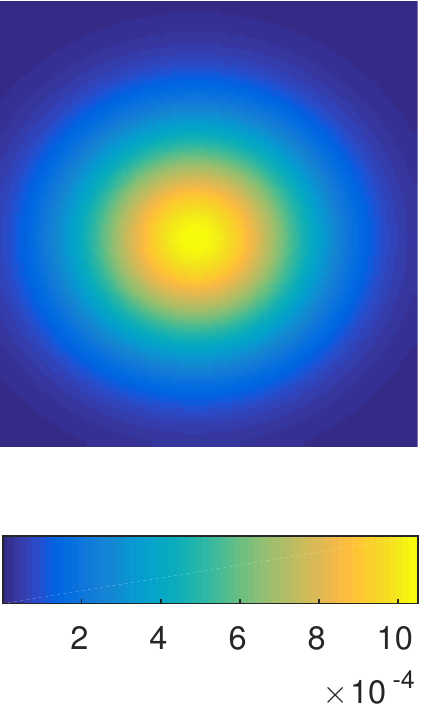}\adjincludegraphics[scale=0.61]{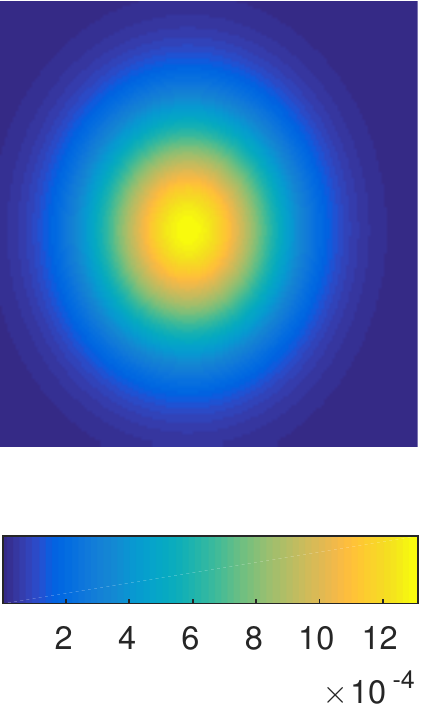}\adjincludegraphics[scale=0.61]{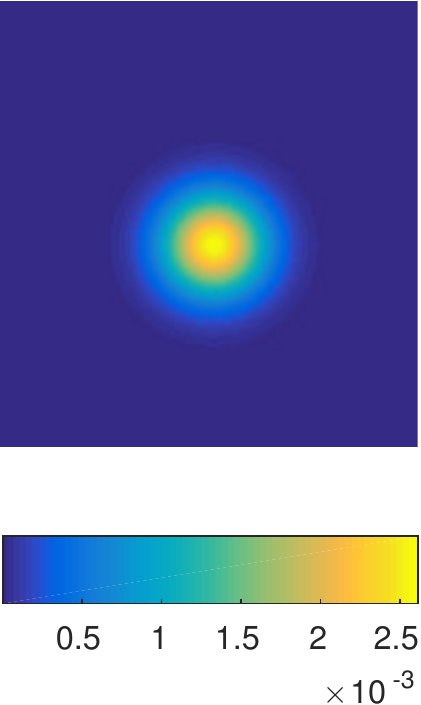}\adjincludegraphics[scale=0.61]{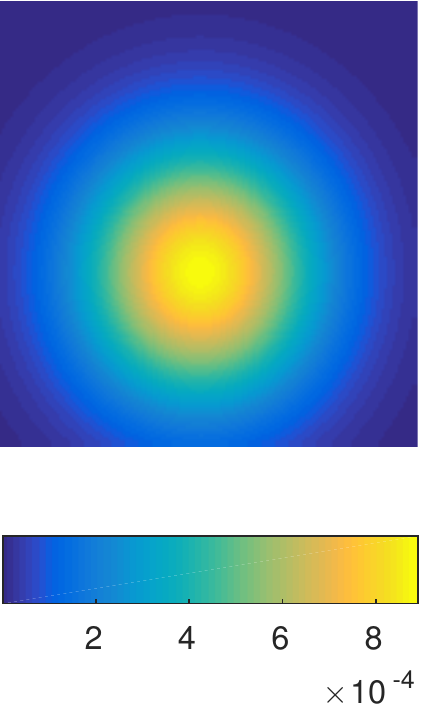}\adjincludegraphics[scale=0.61]{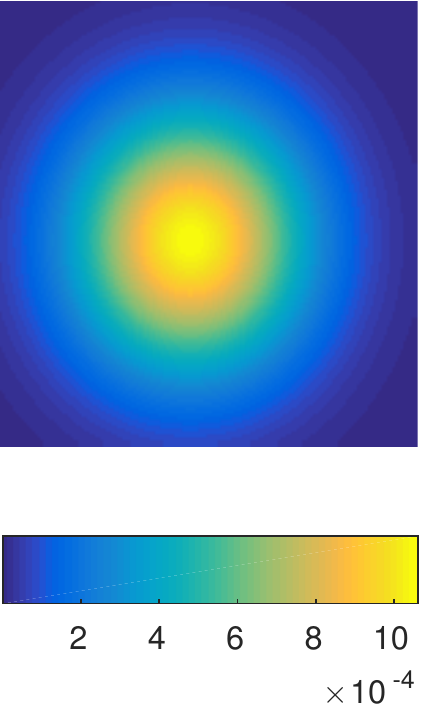}
\caption{(a): Amplitude of the Fourier transform of different healthy and lesion tissues, computed from the Fourier coefficients of the horizontal polarisation. (b): Gaussian function fitted to the amplitude of the Fourier transforms shown in (a). }
	\label{fig:Real_HvL_fourier_Gauss}
\end{figure}

To see if the standard deviation of the fitted Gaussian  can be used to discriminate between healthy and lesion samples in general, we present in Figure~\ref{fig:Real_HvL_fourier_Gauss} the variation of the Fourier-domain information across different tissue samples.
In particular, for each reconstructed image $(n,i)$, $n=1,\ldots,12$, we show the amplitude of the Fourier transform  and the associated Gaussian function.
We can see in Figure~\ref{fig:Real_HvL_fourier_Gauss} that the discriminative behaviour of the standard deviation of the fitted Gaussian between one particular pair of healthy and lesion samples seen in Figure~\ref{fig:Real_HvL_fourier_example}, holds more generally throughout the dataset.

\begin{figure}[h!]
	\centering
	\adjincludegraphics[scale=0.75,trim={{.0\width} {.0\width} {.05\width} {.00\width}},clip]{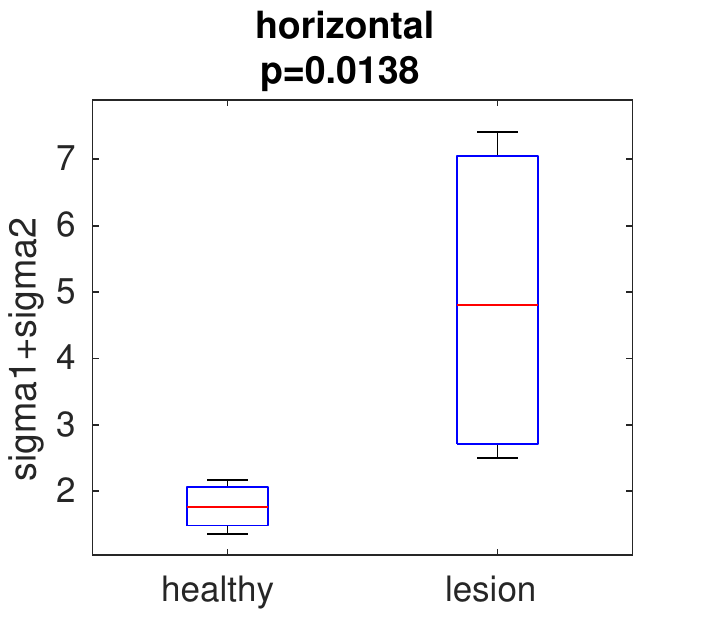} \adjincludegraphics[scale=0.75,trim={{.0\width} {.0\width} {.05\width} {.00\width}},clip]{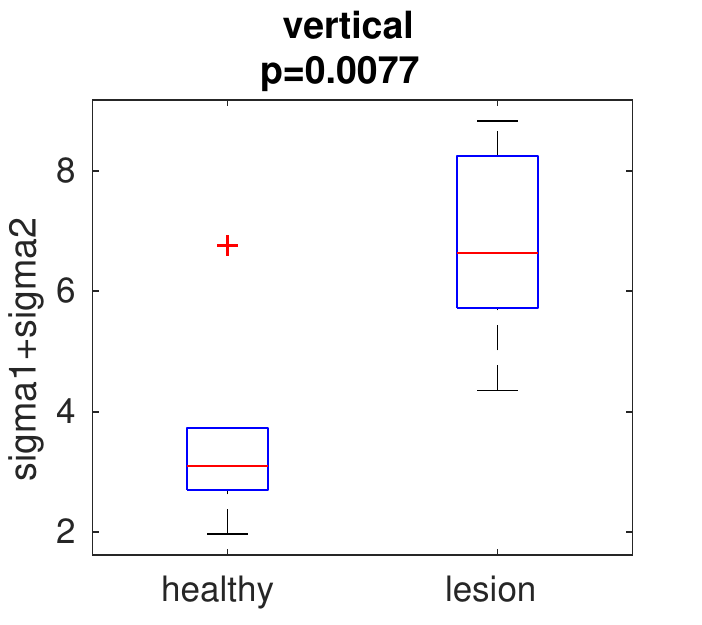} \adjincludegraphics[scale=0.75,trim={{.0\width} {.0\width} {.05\width} {.00\width}},clip]{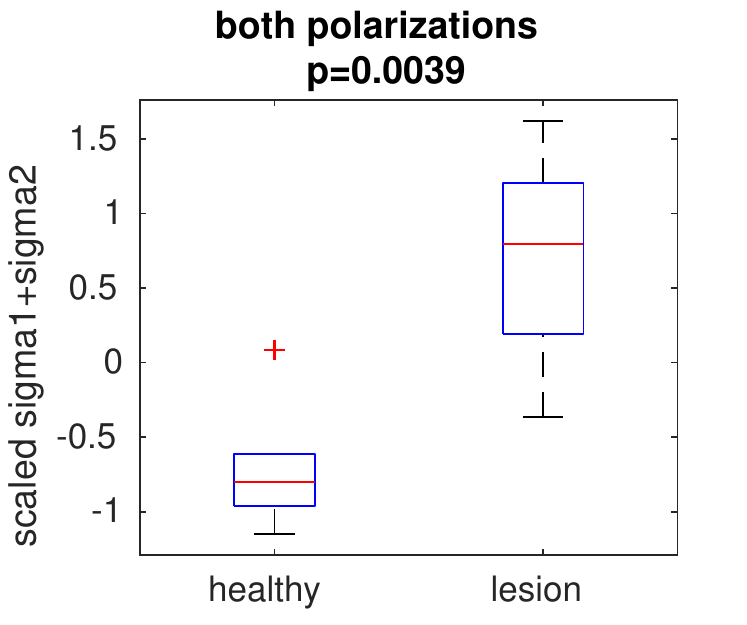}
	\caption{Standard deviation $\sigma_1^{(n)}+\sigma_2^{(n)}$, $n=1,\ldots,12$, of a Gaussian function fitted to the amplitude of the Fourier transform of 6 healthy and 6 lesion areas across our sample set of mouse oesophagus tissue.  In the box-plot on the right, a two-dimensional feature corresponding to two different polarisations is used to characterise each sample, and thus, prior to the t-test, each sample is rescaled by the mean and standard deviation of the total of 12 samples.}
	\label{fig:RealData_boxplot}
\end{figure}

Finally, in Figure~\ref{fig:RealData_boxplot} we perform a statistical test which confirms that the standard deviation of the fitted Gaussians is an informative feature to distinguish between healthy tissues and lesion tissues. In particular, 
for each individual sample $(n,i)$ in the data set, we compute $\sigma_1^{(n,i)}+\sigma_2^{(n,i)}$ of the fitted Gaussian with parameters $\sigma_1^{(n,i)}$ and $\sigma_2^{(n,i)}$. For each tissue sample $n=1,\ldots,12$, we then compute the average value $\sigma_1^{(n)}+\sigma_2^{(n)}:=I_n^{-1}\sum_{i=1}^{I_n}(\sigma_1^{(n,i)}+\sigma_2^{(n,i)})$ and show a box-plot of these twelve values $\sigma_1^{(n)}+\sigma_2^{(n)}$ grouped according to their class label `healthy' or `lesion', for each polarisation as well as for both polarisations combined. We also compute the $p$-value of Welch's t-test \cite{Welch1947}, showing the significant difference in the decay of Fourier coefficients between healthy and lesion samples.

We may conclude that the flexibility of the approach developed in Section \ref{s:framework}, which permits the use of different systems for calibration and image representation, is particularly useful when reconstructing images of biological tissues. The possibility to use a Fourier basis for reconstruction provides the opportunity to investigate images in both the space-domain and the Fourier domain, and thus investigate the degree to which the reconstructed Fourier coefficients decay.  We showed that this decay quantified by the  standard deviation of the fitted Gaussian is a feature with a discriminative power between healthy tissues and lesions, which in the future, in conjunction with a larger data set, could be used to build an automated classifier for distinguishing healthy and lesion samples.

\section{Discussion and future research}\label{s:discussion}

Towards the development and practical use of fibre endoscopes, the main contributions in this paper are two-fold. Firstly, we demonstrated that a Fourier basis yields a diagnostically relevant representation of the optical field reflected from a tissue, using both simulated and experimental real-world data. Secondly, we provided a general reconstruction algorithm that through regularisation can stably recover such representation directly from the calibration measurements
and the measurements of the output optical field transmitted through a fibre, 
where the system used for calibration is allowed to be different and thus more efficient than a Fourier basis.

Nevertheless several open problems remain. One possible direction for future research, which would require a significantly larger number of biological samples to be tested, relates to learning an `optimal' dictionary (alternative to Fourier) as a means to minimise the classification error between healthy and lesion tissues.
More importantly, further work is required to achieve real-time imaging through fibre endoscopes in realistic clinical settings. Specifically, future research is needed to lift the time-independence assumption present in the kernel of the linear model which in everyday clinical use varies across time with bending and temperature. In practice, the time-independence assumption means that the calibration measurements need to be taken often and under similar bending and temperature conditions as when sampling the output optical field,
which is difficult to achieve in realistic clinical deployments. Although recently there have been some initial steps in addressing this issue, as for example \cite{Gu2015}, the development of a clinically-feasible recovery procedure that accounts for significant fibre changes remains an important open problem.

\paragraph{Acknowledgements:} MG and SEB were supported by an EPSRC grant EP/N014588/1 for the centre for Mathematical and Statistical Analysis of Multimodal Clinical Imaging. GSDG and SEB received funding from CRUK (C47594/A16267, C14303/A17197, C47594/A21102) and a pump-priming award from the Cancer Research UK Cambridge Centre Early Detection Programme (A20976). The research of FR was performed in part while he was at University of Cambridge, and it was funded in part by the European Union's Horizon 2020 research and innovation programme under the Marie Sk{\l}odowska-Curie grant agreement No 655282 and in part by the FCT grant SFRH/BPD/118714/2016. The research of AGCPR was conducted during January--June 2017, while visiting University of Cambridge.

\end{document}